\documentclass[11pt]{amsart} %

\usepackage{amsmath, amsthm}
\usepackage{amssymb, mathtools}
\usepackage{extarrows}
\usepackage{mathabx}

\usepackage{relsize}

\usepackage[dvipsnames]{xcolor}

\usepackage{xr}
\usepackage{bigfoot}

\usepackage{enumitem}

\usepackage{todonotes}

\usepackage{accents}

\usepackage{hyperref}

\usepackage{tikz-cd}

\usepackage{endnotes}

\usepackage{booktabs}
\usepackage{tabulary}

\usepackage{soul}
\usepackage{scalerel}

\usepackage[utf8]{inputenc}
\usepackage[greek,english]{babel}
\usepackage{lmodern}

\usepackage{graphicx}

\usepackage{pict2e}

\newtheorem{THEOREM}{Theorem}[section]

\newtheorem{Conclusion}[THEOREM]{Conclusion}

\newtheorem{Theorem}[THEOREM]{Theorem}
\newenvironment{theorem}{\begin{Theorem}}{\end{Theorem}}

\newtheorem{Lemma}[THEOREM]{Lemma}
\newenvironment{lemma}{\begin{Lemma}}{\end{Lemma}}

{\theoremstyle{definition}
\newtheorem{notation}[THEOREM]{Notation}
\newtheorem{definition}[THEOREM]{Definition}

\newtheorem{remark}[THEOREM]{Remark}

\newtheorem*{convention}{Convention}
\newtheorem*{note-nono}{Note}
\newtheorem{note-no}{Note}
\newtheorem{sublemma}[THEOREM]{Sublemma}
}

\newtheorem{Claim}[THEOREM]{Claim}

\newtheorem{Subclaim}[THEOREM]{Subclaim}

\newtheorem{Corollary}[THEOREM]{Corollary}
\newenvironment{corollary}{\begin{Corollary}}{\end{Corollary}}

\newtheorem{Proposition}[THEOREM]{Proposition}
\newenvironment{proposition}{\begin{Proposition}}{\end{Proposition}}

\newtheoremstyle{myrmkstyle}
  {12pt}
  {3pt}
  {\upshape\color{blue}}
  {}
  {\bfseries\color{blue}}
  {.}
  {.5em}
  {}

  \theoremstyle{myrmkstyle}

\newtheoremstyle{myrmkstyleblk}
  {12pt}
  {0pt}
  {\upshape\color{black}}
  {}
  {\bfseries\color{black}}
  {.}
  {.5em}
  {}

  \theoremstyle{myrmkstyleblk}
  \newtheorem*{rmkk}{Remark}

  \newtheorem*{triteob}{Trite Observation}

\newtheoremstyle{privatermkstyle}
  {12pt}
  {3pt}
  {\upshape\color{blue}}
  {}
  {\bfseries\color{blue}}
  {}
  {0em}
  {}

  \theoremstyle{privatermkstyle}

\newtheoremstyle{privatermkstylered}
  {3pt}
  {3pt}
  {\upshape\color{red}}
  {}
  {\bfseries\color{red}}
  {}
  {0em}
  {}

  \theoremstyle{privatermkstylered}

\newenvironment{eqn}{\begin{equation*}}{\end{equation*}}
\newenvironment{itmz}{\begin{itemize}[itemsep=6pt,topsep=-3pt,leftmargin=0.15in]}{\end{itemize}}

\def\>{\rangle}
\def\<{\langle}

\DeclareFontEncoding{LS1}{}{}
\DeclareFontSubstitution{LS1}{stix}{m}{n}
\DeclareSymbolFont{symbols2}{LS1}{stixfrak}{m}{n}

\DeclareMathSymbol{\lcurvyangle}{\mathopen}{symbols2}{"E9}
\DeclareMathSymbol{\rcurvyangle}{\mathclose}{symbols2}{"EA}

\newcommand{\squ}{\scalebox{1.2}{$\square$}}
\newcommand{\squareof}[1]{\raise1.183pt\hbox{\scalebox{0.915}{$\Box(#1)$}}\hskip1.5pt}
\newcommand{\squaresub}[1]{\raise1.183pt\hbox{\scalebox{0.915}{$\Box_{\ssss #1}$}}\hskip1.5pt}

\newcommand{\rge}{\mathop{\mathrm{rge}}} 
\newcommand{\dom}{\mathop{\mathrm{dom}}} 
\newcommand{\rel}{\mathop{\mathrm{rel}}} 
\newcommand{\fun}{\mathop{\mathrm{fun}}} 
\newcommand{\FV}{\mathop{\mathrm{FV}}} 
\newcommand{\Rel}{\mathop{\mathrm{Rel}}} 
\newcommand{\Sent}{\mathop{\mathrm{Sent}}} 
\renewcommand{\lim}{\mathop{\mathrm{lim}}} 

\renewcommand{\to}{\rightarrow}

\newcommand{\ob}{\mathop{\mathrm{ob}}}
\newcommand{\lh}{\mathop{\mathrm{lh}}}

\newcommand{\gtp}{\mathop{\mathrm{gtp}}}
\newcommand{\sh}[1]{\mathop{{\mathrm{Sh}}}({#1})}
\newcommand{\psh}[1]{\mathop{\mathrm{pSh}}({#1})}
\newcommand{\sho}{\mathop{{\mathrm{Sh}}}({\Omega})}
\newcommand{\psho}{\mathop{\mathrm{pSh}}({\Omega})}
 
\newcommand{\concat}{\kern-.25pt\raise4pt\hbox{$\frown$}\kern-.25pt}
\newcommand{\on}{\upharpoonright}  

\newcommand{\Nbb}{\mathbb{N}}

\newcommand{\pret}{\sss\mathrel{\triangleleft}\sss}

\newcommand{\forces}{\Vdash}

\newcommand{\inclusion}{\hookrightarrow}
\newcommand{\Vee}{\bigvee}
\newcommand{\Wedge}{\bigwedge}

\def\sss{\hskip2pt}  
\def\ssss{\hskip1pt}
\def\tupof#1#2{\<\ssss #1\sss:\sss #2\ssss\>} 
\def\tup#1{\<\ssss #1\ssss\>} 
\def\setof#1#2{\{\ssss #1\sss:\sss #2\ssss\}} 
\def\set#1{\{\ssss #1\ssss\}} 
 
\def\gen#1{\lcurvyangle #1 \rcurvyangle} 
\renewcommand{\to}{\mathop{\parbox{.5cm}{\rightarrowfill}}}

\newcommand{\cardrule}{\hrule height.2pt}
\newcommand{\cardrulefill}{\cleaders\cardrule\hfill}
\newcommand{\cardchar}[1]{\vbox{\ialign{##\crcr
    \cardrulefill\crcr\noalign{\kern1pt\nointerlineskip}
    $\hfil\displaystyle{#1}\hfil$\crcr}}}
\newcommand{\barchar}[1]{\vbox{\ialign{##\crcr
    \cardrulefill\crcr\noalign{\kern1.5pt\nointerlineskip}
    $\hfil\displaystyle{#1}\hfil$\crcr}}}
\def\card#1{\cardchar{\cardchar{#1}}}

\newcommand{\forkindep}[2][]{%
  \mathrel{
    \mathop{
      \vcenter{
        \hbox{\oalign{\noalign{\kern-.3ex}\hfil$\vert$\hfil\cr
              \noalign{\kern-.7ex}
              $\smile$\cr\noalign{\kern-.3ex}}}
      }
    }\displaylimits^{#2}_{#1}
  }
}

\title{Some contributions to presheaf model theory}

\author{Andreas Brunner, Charles Morgan and Darllan Concei\c{c}\~ao Pinto}

\address{Departamento de Matem\'atica, Universidade Federal de Bahia,
  Avenida Milton Santos s/n, S/N, Ondina, Salvador, Bahia, Brazil. CEP: 40170-110.}
\email{andreas@dcc.ufba.br}
\email{charlesmorgan@ufba.br}
\email{darllan@ufba.br}

\begin{document}
\maketitle

\begin{abstract} This paper makes contributions to ``pure'' sheaf
  model theory, the part of model theory in which the models are
  sheaves over a complete Heyting algebra. We start by outlining
  the theory in a way we hope is readable for the non-specialist. We then
  give a careful treatment of the interpretation of terms and formulae.
  This allows us to prove various preservation results,
  including strengthenings of the results of \cite{BM14}.
  We give refinements of
  Miraglia's work on directed colimits, \cite{M88}, and an analogue of
  Tarski's theorem on the preservation of $\forall_2$-sentences under
  unions of chains. We next show various categories whose objects are (pairs of) presheaves and sheaves
  with various notions of morphism are accessible in the category theoretic
  sense. Together these ingredients allow us ultimately to prove that
  these categories are encompassed in the AECats framework for
  independence relations developed by Kamsma in \cite{K22}.
\end{abstract}

\section*{Introduction}\label{intro} 

Pure and applied model theory have advanced in consort since (at latest) 
Hrushovski’s proof of the Manin-Mumford conjecture in the 1990s.\footnote{Published in \cite{Hru}.}
Classical model theory continues to have a wide range of new, striking applications,
so wide that it would be invidious to attempt to survey them here.
On the pure side too, \emph{inter alia}, what has come to be known as \emph{neostability theory}
has thriven.\footnote{An extremely rapid
overview of one perspective on neostability can be gleaned from reading the introductions to the recent \cite{KrR},
the three papers \cite{Mut1}, \cite{Mut2} and \cite{Mut3}, or the slightly earlier \cite{d'El}.}

However, arguably even more striking has been the spectacular impact this century of
continous model theory, where truth values are taken in the real
interval $[0,1]$ as opposed to a set of size two,
$\set{T,F}$, as in classical model theory.  At this point of development its
applications are, again, too widespread to survey here.  As for
classical model theory, a `pure' part of the theory, plays an
important r\^ole both intrisically and in applications.\footnote{See
\cite{BYU} as a starting place for continuous stability theory and,
for example, \cite{Kha} or \cite{BHV} for sample excursions in
continuous neostability.}

We assert the rise of continuous model theory is arguably more striking
because it proves the key methodological point that one should not be afraid to
change one's set of ``truth values'' if doing so will help one
\emph{mathematically} in working with particular collections of
mathematical structures and in different mathematical areas. 
Using the terms algebraic and analytic loosely, while classical
model theory has proven marvelously useful in areas where one is
working with `algebraic' structures, continous model theory has in general
turned out to be more appropriate when dealing with more `analytic' structures
and areas.

This paper is a contribution to a third type of model theory,
variously known as the theory of sheaves of models, the theory of sheaf models,
model theory in sheaves or sheaf model theory. The central idea is that one takes the ``truth values''
as lying in a complete Heyting algebra, $\Omega$, and thus one is equipped to
deal with sheaves over $\Omega$. (As will be seen in \S{}\ref{prelims},
complete Heyting algebras are a common generalization of
complete Boolean algebras and the collections of open subsets of topological spaces.)
Clearly, sheaves are a more recent invention
than algebra or analysis,\footnote{See, for example, \cite{Mil},
for an outline of their history.} however over the past 80 years they have
come to play an central r\^ole in mathematics
and we believe the theory of presheaf and sheaf models can make important contributions.

We want to stress immediately that talk about ``truth values'' in
continous model theory and in sheaf model theory is \emph{not}
supposed to be taken literally.  Its aim is to provide an evocative
way for readers to key onto the crucial technical change that formulae
interpreted in models are no longer simply either true or false, but
take values in $[0,1]$ or a complete Heything algebra, respectively.

Continuous model theorists do not advocate for mathematicians to
switch from working in classical (first order) logic to \L{}ukasiewicz's logic,
Pavelka's logic or continuous logic.\footnote{See \cite{Pav},
\cite{Haj} and \cite{BYP} for introductions to these logics.} They do not pursue an agenda
to reform mathe\discretionary{-}{}{}matical argumentation, as perhaps a constructivist might.
The point, rather, is to develop tools analogous to those from classical model theory,
but more suitable for dealing with
parts of mathematics where analysis plays an essential r\^ole,
in order to work in the usual mathematical tradition.

Likewise, our point of view in this paper is to contribute to the development of a theory
analogous to classical model theory, but more suitable for dealing with parts of mathematics
where sheaves plays an essential r\^ole, in order to work in the usual mathematical tradition. 
Two things follow from this. Firstly, this theory is concerned with first order logic and
higher order logic is definitively not something about which it is concerned.\footnote{There
are generalizations of classical model theory to second order logic. See, for example, \cite{Vaa}
for an introduction to this area. In one's daydreams one might consider analogues
of these extensions. However, once coming back to reality these seem rather far off as feasilbe projects.}
Secondly, we absolutely
eschew any interest in molding day to day mathematical practice in, for example,
an intuitionistic direction. Our objectives are very different from those working, say, in
topos theory.\footnote{See, for example, \cite{Joh}.}

Sheaf model theory began in the 1970s with a substantial body of work on
sheaves of models over topological spaces or for Boolean valued models,
This strand of the area has continued to the present. 

An early example is the line of work on model completion, model
companions and certain generalizations of the Feferman-Vaught theorem
in the restricted case of sheaves over Boolean spaces
(\emph{i.e.}, those which are compact, Hausdorff and totally
disconnected), initiated by Lipschitz and Saracino, \cite{LS}, and
Carson, \cite{Carson}, and pursued by Macintyre, \cite{Mac},
\cite{Lou79}, Comer, \cite{com0}, \cite{com}, Burris and Werner
\cite{burwer}, \cite{burwer2} and Pappas, \cite{Pappas}, see also
Volger, \cite{vol1}, \cite{vol2}, \cite{vol3}, using Boolean valued
structures rather than sheaves over Boolean spaces, and Lavendhomme
and Lucas, \cite{LL}, for related results for ``mellow'' sheaves over
Hausdorff spaces without isolated points, and in a slightly different
direction by Weispfennig, \cite{Weis}, Chatzadakis, \cite{Zoe} and
Sureson, \cite{Sureson1}, \cite{Sureson2}. (Chatzadakis works in the
more general context of sheaves over locally compact, zero-dimensional
spaces. Recall that locally compact Hausdorff spaces are totally
disconnected if and only if they are zero-dimensional.) In a similar
setting is Ellerman's work, \cite{ell}, on Boolean valued
generalizations of ultraproducts, and the work of Pierobon and Viale,
\cite{PV}, on \L{}os's theorem for Boolean valued models.

Other examples include the work of 
the work of Prest, alone and variously with Puninskaya, Ralph, Ranjani and Sl\'avik,
on sheaves of modules over (arbitrary) topological spaces, see
\cite{PPR}, \cite{Prest-Ralph1}, \cite{Prest-Ralph2}, \cite{PreRaj},  \cite{Prest-Slavik}, much of which is
surveyed in \cite{Prest}, and
the work of Caicedo and others on sheaves of models over (arbitrary) topological spaces,
particularly with regard to generic model theorems, see Caicedo, 
\cite{Caicedo1}, \cite{Caicedo2}, Sette and Caicedo, \cite{SetCai}, 
Forero, \cite{Forero}, Montoya, \cite{Montoya}, Benavides, \cite{Bena}, and Ochoa and Villaveces, \cite{ochvil}.

However, perhaps unlike in the case of continuous model theory, relatively
early on Fourman and Scott (\cite{FS}) gave a very appealing framework
within which to carry out sheaf model theory over complete Heyting
algebras.  (This framework was also studied independently by Higgs
(\cite{Higgs1}, \cite{Higgs2}), however only limited parts
of \cite{Higgs1} appeared in the published literature in
\cite{Higgs2}, according to Higgs's account there.)

We focus on this Fourman-Scott-Higgs framework here. One reason we do so is because of its elegance. Other motivations
are more concerned with practicality. At present we believe this framework
gives the best combination available of generality and the demonstrability of results
analogous to those of classical model theory and continuous model theory.
For example, in the framework one has versions of:
the downward L\"owenheim-Skolem theorem, due to Miraglia, \cite{M88},
back-and-forth arguments/Fra\"isse constructions, due to Brunner, \cite{B00},
the omitting types theorem, Brunner and Miraglia, \cite{BM04},
the generic model theorem, due to Beer\'io, \cite{Beerio},
\L{}os's theorem, due to Miraglia, \cite{M90}, and Aratake, \cite{Ara}, 
and Robinson's method of diagrams, Brunner and Miraglia, \cite{BM14}.
See also the survey \cite{B16}.

  We should note we are not doctrinaire about the superiority of the Fourman-Scott-Higgs framework. In future
  it may be that similar substantial model theoretic results can be proved for, for example, sheaves over various types
  of quantales or over lineales. See for example the work of Mariano, with variously Alves, Mendes and Ten\'orio,
  \cite{AAM1}, \cite{AAM2}, \cite{TAM}, Reyes-Zambrano, \cite{RZ}, H\"ohle-Kubiak, \cite{HK11}, Miraglia-Solitro, \cite{MirSol},
  Coniglio-Miraglia, \cite{ConMir1}, \cite{ConMir2},
  and Resende, \cite{Resende}, \cite{Resende-lectures}, \cite{Resende-talk},
  for various analogues of the notion of $\Omega$-set for quantales, which would be
  necessary preliminaries for this type of development, and see \cite{deP} for lineales.
  However, as of yet little has been proved model theoretically in any of these directions.
  Hence our evaluation that the Fourmann-Scott framework is currently the most productive one in which to work.

  In this paper we aim to contribute to the ``pure'' side of sheaf model theory, with results in analogues
  of work from the early days of classical model theory, such as various preservation theorems, and with results in
  some first steps in analogues of neostability theory. The former are needed for the latter.

  We briefly summarize the contents of the paper. In \S{}\ref{prelims}
  and \S{}\ref{L-structures} we go over the Fourman-Scott-Higgs
  framework.  \S{}\ref{prelims} discusses presheaves and sheaves over
  complete Heything algebras. Much of this material is known but we
  develop it here freshly and explicitly.  In \S{}\ref{L-structures}
  we introduce $L$-structures in the collection of presheaves over a
  complete Heyting algebra.  We then discuss how terms and the
  realizations of formulae should be interpreted in $L$-structures. We
  this treat relatively carefully as the detail of the interpretation
  of terms are needed for our results in following section and were
  passed over somewhat rapidly in previous presentations. We then
  define, with examples, what it is for the realization of formulas to
  be forced by an $L$-structure. We next briefly discuss extending the
  language with unary connectives corresponding to elements of the
  complete Heyting algebra. Finally, we discuss various notions of
  morphism between $L$-structures. In \S{}\ref{pres_phenom_section} we
  discuss various results on the preservation of the realization of
  formulae being forced under a map between two $L$-structures being
  either an $L$-morphism or an $L$-monomorphism. These culminate in
  our being able to improve substantively
  on some of the results on Robinson's method of diagram contained in \cite{BM14}. 
  In \S{}\ref{cats_and_dir_colims} we introduce categories of interest to
  us with regards to neostability.  The objects in each are pairs
  consisting of either a presheaf and a subpresheaf, a sheaf and a
  subpresheaf, a sheaf and a subsheaf or a sheaf and a subsheaf of density at most that of the
  sheaf, where in each case the larger structure is a model of a
  particular theory. For each collection of objects there are various
  notions of morphism, each leading to a separate category. Ultimately
  the most interesting ones for this paper are those where the
  morphism is a morphism of the larger structure which is either an
  $L$-monomorphism or an elementary $L$-embedding. We prove various
  results concerning the existence of directed colimits for these
  categories, building heavily on Miraglia's results in \cite{M88}. In
  \S{}\ref{accessibility} we discuss the (category-theoretic)
  accessibility of various of our categories. Finally, in
  \S{}\ref{indep_relns} we discuss how the categories dealt with in
  \S{}\ref{accessibility} fit into Kamsma's framework of AECats,
  developed in \cite{K22}, and outline the immediate implications of
  this in terms of uniqueness for independence relations for these
  categories.  \vskip30pt

\section{Presheaves and sheaves}\label{prelims}

We start with basic definitions of partially ordered sets, complete lattices and complete Heyting algebras.

\begin{definition}\label{defn_poset} Let $P$ be a set and $\le$ a binary relation on $P$. The pair $(P,\le)$ is
  a \emph{partially ordered set} (\emph{poset}) if
    \begin{eqn}
    \begin{split}
     & \forall a \in P \sss\sss a\le a \\
     & \forall a, \sss b \in P \sss ( ( a\le b \sss\sss\& \sss\sss b \le a ) \longrightarrow a = b ) \\ 
      &    \forall a, \sss b, \sss c \in P \sss ( (a \le b \sss\sss\& \sss\sss b \le c ) \longrightarrow a \le c )
    \end{split}
    \end{eqn}
  \end{definition}

    \begin{definition}\label{defn_glb_lub} Let $(P,\le)$ be a partially ordered set and $A\subseteq P$.
      $A$ has a \emph{least upper bound} (\emph{lub}), also called the \emph{supremum} or the \emph{join} of $A$,
      if there is some $p\in P$ such that
\[
     \forall a \in A  \sss\sss a\le p \sss\hbox{ and }\sss
     \forall q\in P \sss((\forall a\in A \sss \sss a \le q) \longrightarrow p \le q )
\]
    Similarly $A$ has a \emph{greatest lower bound} (\emph{glb}), also called the \emph{infimum} or the \emph{meet} of $A$,
    if there is some $p\in P$ such that
\[
     \forall a \in A  \sss\sss p\le a\sss \hbox{ and }\sss
     \forall q\in P \sss ((\forall a\in A \sss \sss q \le p) \longrightarrow q \le p )
\]
   The least upper bound and greatest lower bound are respectively denoted, if they exist, by $\Vee A$ and $\Wedge A$.
      \end{definition}

   We have the usual infix notation for dyadic joins and meets.
    
    \begin{notation}\label{defn_vee_wedge_infix} If $(P,\le)$ is a partially ordered set and
      $A = \set{x, y}\subseteq P$ we write
      $x\vee y$ and $x \wedge y$ for $\Vee A$ and $\Wedge A$ respectively they exist.
    \end{notation}
    
    \begin{definition}\label{defn_complete_lattice} A partially ordered set $(L,\le)$ is a
      \emph{complete lattice} if every $A\subseteq L$ has both a least upper bound and a greatest lower bound.
    \end{definition}
    
    \begin{notation} If $(L,\le)$ is a complete lattice we write $\top$ and $\bot$ for
      $\Wedge \emptyset$ and $\Vee \emptyset$, respectively, its greatest and least elements.
    \end{notation}

    \begin{definition}\label{defn_cHa} A complete lattice $(\Omega,\le)$ is a \emph{complete Heyting algebra},
      or \emph{frame}, if
      \[ \forall p\in \Omega \sss\sss \forall A\subseteq \Omega
              \sss \sss \sss \sss p \wedge \Vee A = \Vee \setof{ p \wedge a}{a \in A} .\]
  \end{definition} 

    \begin{definition}\label{defn_implication_negation_in_cHa} If $(\Omega,\le)$ is a complete Heyting algebra and $p$,
      $q\in \Omega$
      we define
  \begin{eqn}
    \begin{split}
      &\hbox{$p \to q$ ($\in \Omega$) as $\Vee \setof{r \in \Omega}{r\wedge p \le q}\sss$ and
        $\sss\neg p$ ($\in \Omega$) as $p\to \bot$.}
    \end{split}
    \end{eqn}
    \end{definition}

    \begin{lemma}\label{adjunction} (Adjunction) Let $(\Omega,\le)$ be a complete Heyting algebra.
      If $p$, $q$, $r\in \Omega$ then
           $r \le p \to q$ if and only if $r\land p \le q$.
    \end{lemma}

    We make a few simple observations about identities and inequalities in complete Heyting algebras will be useful in the ensuing, particularly
    in \S{}4.
    \vskip12pt

\begin{lemma}\label{modus_ponens} (\emph{Modus ponens}.) Let $p$ and $q\in \Omega$.
Then $p\land (p\rightarrow q) = p \land q$.
\end{lemma}

\begin{proof} As $p\land q \le q$ always holds, adjunction shows $q\le p\to q$.
  Taking the meet ($\land$) with $p$ on both sides gives $p\land q \le p \land (p\to q)$.

  Conversely, as $p\to q \le p\to q$ holds, adjunction gives $p\land (p\to q) \le q$.
  Taking the meet with $p$ on both sides gives $p\land (p\to q) \le p\land q$
\end{proof}

\vskip12pt
  \begin{lemma}\label{char_equiv_in_Omega} Let $p$, $q$ and $r\in \Omega$.
    Then $p = p \land ( q \leftrightarrow r )$ if and only if $p \land q = p \land r$.
\end{lemma}

  \begin{proof} If $p = p \land ( q \leftrightarrow r )$ then $p\le ( q \leftrightarrow r )$ and, by adjunction (Lemma (\ref{adjunction}),
    $p\land q \le r$ and $p\land r \le q$. Thus
    $p\land q \le p\land r$ and $p\land r \le p \land q$. Conversely, if $p \land q = p \land r$ then $p \land q \le r$ and $p\land r \le q$, and so by adjunction
    $p\le q \to r$ and $p\le r\to q$. 
  \end{proof}

  \vskip12pt
  \begin{lemma}\label{char_equiv_d} Let $p$, $q\in \Omega$ and let $d= p\to q$.
\begin{enumerate}[align=parleft,labelsep=1cm, label=(\arabic*),itemsep=6pt,topsep=0pt]
\item The identity $p = p \land ( q \leftrightarrow d )$ always holds.
\item The identity $ p = p \land ( p \leftrightarrow d )$ is true if and only if $p\le q$.
\end{enumerate}
  \end{lemma}

  \begin{proof} \emph{(1).} By Lemma (\ref{char_equiv_in_Omega}),
    $p = p \land ( q \leftrightarrow d ) $ holds if and only if $p \land q = p \land d$ does. However, $p \land d = p \land (p\to q) = p \land q$ by
    \emph{modus ponens}. Thus $p \land q = p \land d$ always holds.

    \emph{(2).}  Similarly, by Lemma (\ref{char_equiv_in_Omega}),
    $p = p \land ( p \leftrightarrow d )$ holds if and only if  $ p = p\land p = p\land d$ does. However, as in the proof of \emph{(1)},
    $p\land d = p \land (p\to q) = p \land q$, so $p = p \land ( p \leftrightarrow d )$ holds if and only if $p = p\land q$, that is, if and only if $p\le q$.
  \end{proof}

\begin{definition} Let $\Omega$ be a complete Heyting algebra and let $p$, $q\in \Omega$.
  We say $p$ is \emph{dense in} $q$ if $p \le q \le \neg\neg p$ and $p$ is \emph{dense}
  if it is dense in $\top$: $\neg\neg p = \top$.
\end{definition}
    
We next introduce presheaves over complete Heyting algebras and the two types of products of them that
are used pervasively in the paper.\footnote{See, for example, \cite{Shaf} or \cite{Tenn},
for introductions to and examples of presheaves and sheaves over topological spaces. See
\cite{Borceux} and \cite{TvD} for introductions 
to presheaves and sheaves over complete Heyting algebras
similar to, but somewhat more effusive, than this one, in the former interwoven with an introduction
to presheaves and sheaves over topological spaces.}

\begin{definition}\label{defn_presheaf} Let
  $\Omega$ be a complete Heyting algebra.  $M = (|M|,\on^M,E^M)$ is a
  \emph{presheaf over} $\Omega$ if
        $|M|$ is a set, $\on^M:|M|\times\Omega \longrightarrow |M|$ and $E^M:|M|\longrightarrow \Omega$, 
      and for all $a\in |M|$ and $p$, $q\in \Omega$,
  \begin{eqn}
    \begin{split}
      & a\on^M (p\wedge q) =   (a\on^M p)\on^M q\\
      & a\on^M E^Ma = a \\
      & E^M(a\on^M p) = E^Ma \wedge p.
    \end{split}
    \end{eqn}where, as we will throughout, we write $E^M a$ for $E^M(a)$ -- the \emph{extent of}
  $a$ -- and $a\on^M p$ for $\on^M(a,p)$ -- \emph{the restriction of
  $a$ to $p$}.

  When clear from the context -- and it almost invariably will be clear -- we drop the superscripts and
  simply write $E$ and $\on$ for $E^M$ and $\on^M$, respectively.
\end{definition}

\vskip6pt
\begin{definition}\label{defn_extensional}
  If $\Omega$ is a complete Heyting algebra,
  a presheaf $M$ over $\Omega$ is \emph{extensional} or \emph{separated}
  if for all $a$, $b\in |M|$ and $D\subseteq\Omega$
  \[\big((\forall p \in D\sss\sss a\on p=b\on p)\sss\sss\&\sss\sss
  Ea=Eb=\bigvee D\big) \Longrightarrow a=b .\]
\end{definition}

\begin{rmkk}   From the point of view of the model theory of presheafs of first-order structures,
  which will be our focus here
  and which are introduced in Definition \ref{defn_L_str_in_psho},
  nothing is lost by considering only extensional presheaves as there
  is a canonical and well-behaved process of extensionalization of a
  presheaf of first-order structures, as described in \cite{Mir06},
  Theorem (23.21).
\end{rmkk}
Consequently we adopt the following convention.

\begin{convention} From here on all presheaves (and sheaves) are assumed to be extensional.
\end{convention}

We shall make extensive use of the $\Omega$-set structure derived from a presheaf.
\begin{definition}\label{defn_evaluation_for_presheaves}
  Let $\Omega$ be a complete Heyting algebra and
  $M$ be a presheaf over $\Omega$. Define $[.=.]_M:|M|\times
  |M|\longrightarrow \Omega$ by for $a$, $b\in |M|$ setting
    \begin{eqn}
    \begin{split}
      [a=b]_M = & \bigvee \setof{E(a\on p) \wedge E(b\on p)}{p\in \Omega \sss\sss \& \sss\sss a\on p = b\on p} \\
      \big( = & \bigvee \setof{E(a\on p) }{p\in \Omega \sss\sss \& \sss\sss a\on p = b\on p} \big).\\
   \end{split}
  \end{eqn}
     \end{definition}

\begin{lemma}\label{[a=b]_below_Ea_Eb} Let $\Omega$ be a complete Heyting algebra, 
  $M$ be a presheaf over $\Omega$ and $a$, $b\in |M|$. Then $[a=b]_M \le Ea$, $Eb$.
  \end{lemma}

\begin{proof} Immediate from the definition of $[a=b]_M$.
\end{proof}

\begin{lemma}\label{[a=b]_[b_=c]_le_[a=c]} Let $\Omega$ be a complete Heyting algebra, 
  $M$ be a presheaf over $\Omega$ and $a$, $b$ and $c\in |M|$. Then
  $[a=b]_M \land [b=c]_M \le [a = c]_M$.
\end{lemma}

\begin{proof} By definition we have \begin{eqn}
  \begin{split}
  [a=b]_M \land [b=c]_M & =  \bigvee \setof{E(b\on p) }{p\in \Omega \sss\sss \& \sss\sss a\on p = b\on p} \\ &\land
  \bigvee \setof{E(b\on q) }{q\in \Omega \sss\sss \& \sss\sss c\on q = b\on q}.
  \end{split}
  \end{eqn}
  By distributivity of $\land$ over $\bigvee$ twice this gives
  \begin{eqn}\begin{split} [a=b]_M \land [b=c]_M = &  \bigvee \setof{E(b\on p) \land E(b\on q)}{ p,\sss q\in \Omega \sss\sss \& \sss\sss a\on p = b\on p \\ & \qquad \quad  \sss\sss \& \sss\sss b\on q = c\on q } \\
      = & \bigvee \setof{E(b\on p\land q)}{ p,\sss q\in \Omega \sss\sss \& \sss\sss a\on p\land q  = \\ & \qquad \quad b\on p\land q \sss\sss \& \sss\sss b\on p\land q = c\on p\land q }  = [a=c]_M
    \end{split}
  \end{eqn}
\end{proof}

Note, we clearly cannot expect equality in Lemma (\ref{[a=b]_[b_=c]_le_[a=c]}), for example, take $a=c \ne b$.

As a first example we prove an equivalent criteria in terms of this structure for a presheaf to be existential.
\begin{lemma}\label{equiv_extensional}
  Let $\Omega$ be a complete Heyting algebra and
  $M$ a presheaf over $\Omega$. $M$ is {extensional} if and only if  for all $a, b\in |M|$
  $ [a=b]_M = Ea = Eb $ implies $a=b$.
\end{lemma}

\begin{proof} In one direction, suppose $M$ is extensional and let $a$, $b\in |M|$ be such that $ ([a=b]_M = Ea = Eb $. By Definition (\ref{defn_evaluation_for_presheaves})
  we have
  \begin{eqn}\begin{split} [a=b]_M = & \bigvee \setof{E(a\on p) }{p\in \Omega \sss\sss \& \sss\sss a\on p = b\on p}\\
      = & \bigvee \setof{p }{p\in \Omega \sss\sss \& \sss\sss p \le Ea \sss\sss \& \sss\sss a\on p = b\on p} .
  \end{split}\end{eqn}

  Define $D= \setof{p\in \Omega}{p \le Ea \sss\sss\&\sss\sss a\on p = b\on p}$. Then $[a=b]_M = \bigvee D$. By entensionality we thus have $a=b$.

  Conversely, let $a$, $b\in |M|$. By Lemma (\ref{[a=b]_below_Ea_Eb}), $[a=b]_M \le Ea$, $b$.
  Let $D\subseteq \Omega$ be such that for all $p\in D$ we have $a\on p = b\on p$ and $Ea=Eb=\bigvee D$. We thus
  have
  \begin{eqn}\begin{split}  Ea= Eb =\bigvee D & = \bigvee\setof{p \in D}{a\on p = b\on p} \\ & \le \bigvee\setof{p \le \bigvee D}{a\on p = b\on p} = [a=b]_M.
    \end{split}
  \end{eqn}
  Consequently $[a=b]_M = Ea = Eb$, and so by hypothesis $a=b$.
\end{proof}

\vskip12pt

\begin{definition}\label{defn_cartesian_prod}
  If $\Omega$ is a complete Heyting algebra and $M_0$, \dots, $M_{n-1}$ are 
  presheaves  over $\Omega$ then $|M_0|\times\dots\times |M_{n-1}|$
  is simply the \emph{Cartesian product} of the $M_i$:
  $\setof{\tup{a_0,\dots,a_{n-1}}}{\forall i < n \; a_i \in |M_i|}$. If $M$ is such
  that for all $i<n$ we have $M=M_i$ we write $|M|^n$ for the
  $n$-fold \emph{Cartesian power} of $|M|$.
\end{definition}

\begin{definition}\label{defn_presheaf_prod}
  If $\Omega$ is a complete Heyting algebra, $n\ge 1$ and
  $M_0$, \dots , $M_{n-1}$ are presheaves  over $\Omega$ then
  $M_0 \times \dots\times M_{n-1}$ is the \emph{presheaf product}   of the $M_i$,
  $M_0 \times \dots\times M_{n-1}=
  (|M_0\times\dots\times M_{n-1}|,
  \on^{M_0\times\dots\times M_{n-1}} ,
  E^{M_0\times\dots\times M_{n-1} }))$, where
  \begin{eqn}
    \begin{split}
        |M_0\times\dots\times M_{n-1}|   = & \\ 
          \setof{\tup{a_0,\dots,a_{n-1}}}{ & \forall i < n \; a_i \in |M_i|
    \sss\sss\&\sss\sss \forall i, \sss j <n\; E^{M_i}a_i = E^{M_j}a_j},
    \end{split}
  \end{eqn}if $p\in \Omega$ then
  $\tup{a_0,\dots,a_{n-1}}\on^{M_0\times\dots\times M_{n-1} } p = \tup{a_0\on^{M_0} p,\dots,a_{n-1}\on^{M_{n-1}} p}$
    and $E^{M_0\times\dots\times M_{n-1} }\tup{a_0,\dots,a_{n-1}} = E^{M_0}a_0$.
\end{definition}

\begin{definition} If $M$ is such that for all $i<n$ we have $M=M_i$ we write $M^n$ for the
  $n$-fold \emph{presheaf power of} $M$ and note that
  \[|M^n| = \setof{\tup{a_0,\dots,a_{n-1}} }{\forall i < n \; a_i \in |M|
  \sss\sss\&\sss\sss \forall i,\ssss j < n\; E^{M}a_i = E^{M}a_j}.\]
\end{definition}

We now introduce morphisms between presheaves and related notions such as inclusions and monomorphisms.

\begin{definition}\label{defn_presheaf_morphism} Let
  $\Omega$ be a complete Heyting algebra and
  $M$, $N$ presheaves over $\Omega$. A \emph{presheaf morphism}
  $M \xlongrightarrow{f} N$ is a function ${f:|M|\longrightarrow |N|}$
  such that for all $a\in |M|$ and $p\in \Omega$ we have
  \[ E^N f(a) = E^M a \sss\sss\&\sss\sss f(a\on^M p) = f(a)\on^N p.
  \]
\end{definition}

\begin{definition}\label{defn_psh} If
  $\Omega$ is a complete Heyting algebra then
    $\psh{\Omega}$ is the category whose objects are presheaves over $\Omega$
    and whose morphisms are presheaf morphisms.
  \end{definition}

\begin{proposition}\label{Omega_terminal} If $\Omega$ is a complete Heyting algebra
  it is the terminal object in the category $\psh{\Omega}$.
  \end{proposition}
\begin{proof} Define $\on^\Omega$ and $E^\Omega$ by setting
  $q \on^\Omega p = q \wedge p$ and $E^\Omega q =q$, for $p$, $q\in\Omega$. Then it is immediate to
  check that $(\Omega,\on^\Omega,E^\Omega)\in\psh(\Omega)$ satisfies Definition (\ref{defn_presheaf}).
  Furthermore,  if $M \in \psh{\Omega}$ then, by Definition (\ref{defn_presheaf})
  and Definition (\ref{defn_presheaf_morphism}), $E^M$ is a presheaf morphism from
  $M$ to $(\Omega,\on^\Omega,E^\Omega)$.
\end{proof}

\vskip6pt

In view of Proposition (\ref{Omega_terminal}) we introduce a convention that otherwise perhaps would at first sight be surprising. 

\begin{notation}\label{defn_M^0} Let $M$ be a presheaf over $\Omega$. We set $M^0 = \Omega$.
\end{notation}
\vskip-6pt
We comment further on this notation when we discuss the interpretation of terms in $L$-structures prior to and in
Definition (\ref{defn_interp_terms_in_l_str_in_psho}) below.
\vskip12pt

\begin{lemma}\label{psh_morphism_push_up_equality}  Let
  $\Omega$ be a complete Heyting algebra, 
  $M$, $N$ presheaves over $\Omega$, 
  $M \xlongrightarrow{f} N$ a presheaf morphism, and $a$, $b\in |M|$. Then \[[a=b]_M \le [f(a) = f(b)]_N.\]
\end{lemma}

\begin{proof} If $p\in \Omega$ and $a\on^M p = b\on^M p$ then $f(a) \on^N p = f(a\on^M p) = f(b\on^M p)=f(b)\on^N p$ and
  $E^N f(a)\on^N p = E^N f(a\on^M p) = E^M a\on^M p$. So the set over which we take the disjunction in the definition of
  $[a=b]_M$ is a subset of that over which we take the disjunction to evaluate $[f(a)=f(b)]_N$.
\end{proof}
\vskip6pt

We give a converse to this last lemma.

\begin{lemma}\label{equal_push_up_ensures_morphism} Let
  $\Omega$ be a complete Heyting algebra and
  $M$, $N$ presheaves over $\Omega$.  If $f:M\longrightarrow N$ is a function such that for all  $a$, $b\in |M|$ we have
  $E^Ma = E^Nf(a)$ and $[a=b]_M \le [f(a)=f(b)]_N$ then $f$ is a presheaf morphism.
\end{lemma}

\begin{proof} Suppose $a \in |M|$ and $p\in \Omega$. We want to show $f(a\on p)= f(a)\on p$. By extensionality, it suffices to
  show \[ [f(a\on p)= f(a)\on p]_M = E^Nf(a\on p) = E^Nf(a)\on p.\]

  However, $E^Nf(a)\on p = p\land E^Nf(a) = p\land E^Ma = E^Ma\on p = E^Nf(a\on p)$. Also,
  $p\land E^Ma = [a\on p = a]_M \le [f(a\on p) = f(a)]_N \le E^Nf(a\on p) = p\land E^Ma$.
  \end{proof}

\begin{definition}\label{defn_subpresheaf} 
  Let $\Omega$ be a complete Heyting algebra,
  $M$, $N$ presheaves over $\Omega$ and $|M|\subseteq |N|$.
  Then $M$ is a \emph{subpresheaf} of $N$ if the identity immersion is a presheaf morphism.
\end{definition}

\begin{note-nono}
  Any subpresheaf of an extensional presheaf is extensional. (This is immediate from the definitions
  of extensional and subpresheaf.)
\end{note-nono}

  We next define the restriction of a presheaf to an element of $\Omega$.

  \begin{definition}\label{rest_presheaf_to_an_elt_of_Omega}
    Let $\Omega$ be a complete Heyting algebra,
    $M$ is a presheaf over $\Omega$ and $p\in \Omega$.
    Set $|M\on p| = \setof{a \in |M|}{Ea \le p}$
    and for all $a\in |M\on p|$ and $q\in\Omega$ set $a\on^p q = a\on q$
    and $E^p a = Ea$. Then define
    ${M\on p = (|M\on p|,\on^p, E^p)}$, 
    \end{definition} 

  \begin{note-nono} It is immediate from the definition of $M\on p$ and the definition of presheaf morphism that
    if $\Omega$ is a complete Heyting algebra,
    $M$ is a presheaf over $\Omega$ and $p\in \Omega$ then
    $M\on p$ is a subpresheaf of $M$.
  \end{note-nono}

\begin{definition}\label{defn_presheaf_mono}
  Let $\Omega$ be a complete Heyting algebra and
  $M$, $N$ presheaves over $\Omega$. A presheaf morphism
  $M \xlongrightarrow{f} N$ is a \emph{monomorphism} if
  $f:|M|\longrightarrow |N|$ is injective.
\end{definition}

\begin{proposition}\label{equiv_presheaf_mono} (\emph{cf}.~\cite{Mir06}, Lemma 25.21.)
    Let $\Omega$ be a complete Heyting algebra,
  $M$, $N$ presheaves over $\Omega$ and
  $M \xlongrightarrow{f} N$ be a presheaf morphism. 
  Suppose $A$, $B$ are subpresheaves of $M$, $N$, respectively, and
  that $f\on |A|:|A|\longrightarrow |B|$ is a presheaf morphism
  (\emph{i.e.}, $\rge(f\on |A|)\subseteq |B|$). Then
  $f\on |A|$ is a presheaf monomorphism if and only if
  for all $a$, $b\in |A|$ we have,
  in the notation of Definition (\ref{defn_evaluation_for_presheaves}),
  $ [a=b]_A = [f(a)=f(b)]_B$, or, more explicitly,
  \[ \bigvee \setof{p\wedge Ea \wedge Eb}{a\on p = b\on p} =
  \bigvee \setof{p\wedge Ea \wedge Eb}{f(a\on p) = f(b\on p)}.
  \]  
  \end{proposition}

  \begin{proof}
    If $f\on |A|$ is a presheaf monomorphism then for all $a$, $b\in
    A$ and $p\in\Omega$, $a\on p$, $b\on p\in |A|$ and we have
    $f(a\on p) = f(b\on p)$ if and only if $a\on p=b\on p$, so the
    disjunctions in the definitions of $[a=b]_A$ and $[f(a)=f(b)]_B$
    are taken over the same set.

    Conversely, suppose $a$, $b\in |A|$ and $f(a)=f(b)$. Then for all $p\in\Omega$ we have
    $f(a\on p)= f(a)\on p = f(b)\on p=f(b\on p)$, since $f$ is a presheaf morphism, and so
    $[f(a)=f(b)]_B = Ef(a)=Ef(b)$, by the definition of $[.=.]_B$. Moreover, $Ef(a)=Ef(b)= Ea=Eb$, since
    $f$ is a presheaf morphism. Thus, if $[a=b]_A=[f(a)=f(b)]_B$ we have $[a=b]_A = Ea=Eb$,
    and so $a=b$ by the extensionality of $M$.
  \end{proof}

\vskip6pt

We now discuss compatibility and glueing of sets of elements of a presheaf.
This will lead to the definition of a sheaf.

\begin{definition}\label{defn_compatible_set_of_sections}
  Let $\Omega$ be a complete Heyting algebra,
  $M$ a presheaf over $\Omega$ and $C\subseteq |M|$. $C$ is \emph{compatible}
  if for all $a$, $b\in C$ we have
  \[Ea\land Eb=\bigvee \setof{p\wedge Ea \wedge Eb}{a\on p = b\on p} \sss\sss\sss( = [a=b]_M).\]
  \end{definition}

  \begin{definition}\label{defn_glueing}
  Let $\Omega$ be a complete Heyting algebra,
  $M$ a presheaf over $\Omega$ and $C\subseteq |M|$ is compatible.
  A \emph{glueing of} $C$ is a element $b\in |M|$ such that
  \[ Eb = \bigvee_{a\in C} Ea \sss\sss\sss\sss\&\sss\sss\sss\sss
  \forall a\in C \sss\sss Ea = [a=b]_M \sss\sss\big( =
  \bigvee \setof{p\wedge Ea}{a\on p = b\on p}\big) .  \]
  \end{definition}

\begin{lemma}   Let $\Omega$ be a complete Heyting algebra,  $M$ a presheaf over $\Omega$ and $C\subseteq |M|$ a compatible set.
  Then $b\in |M|$ is a glueing of $C$ if and only if $Eb = \bigvee \setof{Ea}{a\in C}$ and for all $a\in C$ we have $b\on Ea = a$.
\end{lemma}

\begin{proof} For left-to-right, suppose $b$ is a glueing of $C$ and $a\in C$. We have
  \[ [b\on Ea = a]_M = Ea \land [b=a]_M  = Ea \land Ea = Ea  .\] However we also have $E(b\on Ea) = Eb\land Ea = Ea$. Thus, by extensionality, $b\on Ea = a$.

  For right-to-left, suppose $a\in C$ and $b\on Ea = a$. Then
  \[ Ea = E(b\on a) = Eb\land Ea = [b\on Ea = a]_M = Ea\land [b=a]_M,\]
  so $Ea\le [ a=b]_M$. By Lemma (\ref{[a=b]_below_Ea_Eb}), we also have $[a=b]_M \le Ea$, and thus have $[a=b]_M = Ea$. 
\end{proof}

     We show that no two elements of an extensional presheaf 
   are represented by the same system of pairs.

   \begin{lemma}\label{systems_rep_uniquely}   Let $\Omega$ be a complete Heyting algebra,
     $M$ an extensional presheaf over $\Omega$ and $D\subseteq |M|$.
     Let $\setof{(d_i,p_i)}{ i\in I}\subseteq D\times \Omega$. If $a$, $b\in |M|$,
     $Ea=\bigvee_{i\in I} p_i = Eb$ and for all $i\in I$ we have $a\on p_i = d_i \on p_i = b\on p_i$ then $a=b$
   \end{lemma}

   \begin{proof} By the definition of $[.=.]_M$ we have that
     $[a=b]_M = \bigvee X$ where $X=\setof{p\le Ea\wedge Eb}{a\on p = b\on p}$.
     However for each $i\in I$ we have that $p_i\le Ea\wedge Eb = Ea = Eb$ and
     $a\on {p_i}=d_i\on p_i = b\on p_i$, so $p_i\in X$. Thus
     \[Ea=Eb\ge [a=b]_M = \bigvee X \ge \bigvee_{i\in I}p_i = Ea = Eb.\]
     Hence $Ea=Eb=[a=b]_M$. By extensionality this gives that $a=b$
   \end{proof}
   \vskip12pt
   
  \begin{definition}\label{defn_sheaf} Let
    $\Omega$ be a complete Heyting algebra. A
    presheaf $M$ over $\Omega$ is a \emph{sheaf} (or \emph{is complete})
    if every compatible subset of $|M|$ has a (unique) glueing in $|M|$.
  \end{definition}

    \begin{definition}\label{defn_sh} Let
    $\Omega$ be a complete Heyting algebra.
    $\sh{\Omega}$ is the full subcategory of $\psh{\Omega}$ whose objects
    are sheafs over $\Omega$ (and, by fullness, whose morphisms are the
    presheaf morphisms between sheafs).
    \end{definition}

    \vskip12pt

    In the ensuing we will need the notions of a subpresheaf and a subsheaf
    generated by a subset of a presheaf or sheaf.
    Here we introduce the relevant definitions for presheaves and sheaves as defined in
    Definition (\ref{defn_presheaf}) and Definition (\ref{defn_sheaf}). We first give abstract definitions and
    then show these are equivalent to a concrete construction.
    
  \begin{lemma}\label{intersection_of_presh_lemma} (\cite{Mir06}, 24.10) Let
    $\Omega$ be a complete Heyting algebra, $M$ a presheaf over $\Omega$, $I$ a set and suppose for all $i\in I$ we have that $M_i$ a
  subpresheaf of $M$. Then $\bigcap\setof{M_i}{i\in I}$ is a subpresheaf of $M$. Further, if for each $i\in I$ we have
  $M_i$ is a subsheaf of a sheaf $M$ then $\bigcap\setof{M_i}{i\in I}$ is a subsheaf of $M$.
\end{lemma}

  \begin{definition}\label{defn_of_psh_gen_by} (\cite{Mir06}, 24.11, 24.13 )
    Let  $\Omega$ be a complete Heyting algebra,  $M$  a presheaf over $\Omega$
  and $A\subseteq |M|$. Write for $\gen{A}$, the \emph{subpresheaf of} $M$ \emph{generated by} $A$: 
   $\gen{A}^M = \bigcap\setof{Q\subseteq M}{Q\hbox{ is a presheaf and }A\subseteq |Q|}$.
  If $M$ is also a sheaf over $\Omega$,
  write $\gen{A}^M_s$ for $\bigcap\setof{Q\subseteq M}{Q\hbox{ is a sheaf and }A\subseteq |Q|}$,
  the \emph{subsheaf of} $M$ \emph{generated by} $A$. (Here we use the previous lemma to see that
  $\gen{A}$ really is a subpresheaf, resp.,~$\gen{A}_s$ a subsheaf.)
\end{definition}
We show this definition is unambiguous in the sense that going to larger ambient presheaves or sheaves does not
change the generated subpresheaf or subsheaf, and thus the superscripts in the definition are superfluous.

\begin{lemma}\label{uniqueness_of_gen_by} Let
    $\Omega$ be a complete Heyting algebra, $M$, $N$ presheaves over $\Omega$, with
  $M$ being a subpresheaf of $N$ and $A\subseteq |M|$. Then
  $\gen{A}^M = \gen{A}^N$.  
  If $M$, $N$ are sheaves over $\Omega$, $M$ is a subsheaf of $N$ and $A\subseteq |M|$ then
  $\gen{A}_s^M = \gen{A}_s^N$.  
\end{lemma}

\begin{proof} In each case $M$ is amongst the `$Q$'s whose intersections are taken
  in the definition of $\gen{A}^N$ and $\gen{A}_s^N$ respectively in Definition (\ref{defn_of_psh_gen_by}). 
  \end{proof}

\begin{notation} Let  $\Omega$ be a complete Heyting algebra,  $M$  a presheaf over $\Omega$
  and $A\subseteq |M|$. By Lemma (\ref{uniqueness_of_gen_by}), we simply write
  $\gen{A}$ for $\gen{A}^M$, and if $M$ is a sheaf over $\Omega$ we write
  $\gen{A}_S$ for $\gen{A}_s^M$
\end{notation}

\begin{proposition} Let  $\Omega$ be a complete Heyting algebra,  $M$  a presheaf over $\Omega$
  and $A\subseteq |M|$.  Then
  \[ |\gen{A}| = \setof{a\on p}{a\in A \sss\sss\&\sss\sss p \in \Omega} .\] 

  If $M$ is also a sheaf over $\Omega$, then
  \[ |\gen{A}_s| =  \setof{ b \in M}{b \hbox{ is the glueing of a set of compatible elements of }|\gen{A}|} .\]
\end{proposition}

\begin{proof} Clearly if $A\subseteq |Q|$, where $Q$ is a subsheaf of $M$, $a\in A$ and $p\in \Omega$, then $a\on p \in Q$ by the closure of $Q$ under restriction.
  However, if we equip $P = \setof{a\on p}{a\in A \sss\sss\&\sss\sss p \in \Omega}$ with the operators $E^P$ and $\on^P$ given by
  $E^Pa\on p = E^Ma \land p$ and $(a\on p) \on^P q= a \on p\land q$, then $(P,E^P,\on^P)$ is a subpresheaf of $M$ with $A\subseteq P$.
  Thus $P=\gen{A}$.
  
  If $M$ is a sheaf and $Q$ is a subsheaf of $M$ with $A\subseteq |Q|$ then since $Q$ is a subpresheaf of $M$ we have
  $|\gen{A}|\subseteq |Q|$ and the glueing of every compatible set of elements of $\gen{A}$ is an element of $Q$. 
  However, clearly the set of glueings of a set of compatible elements of $|\gen{A}|$ is closed under further glueings. 
\end{proof}

\vskip12pt
Whilst, as just noted, we do make use of the notion of a sub(pre)sheaf generated by a set of elements of a (pre)sheaf,
a separate notion of the closure of a set of elements of a (pre)sheaf will also be crucial. 

  \begin{definition} Let
    $\Omega$ be a complete Heyting algebra and $M$ a presheaf over $\Omega$.
  For  $A\subseteq |M|$ define $\hbox{cl}(A)$ to be the collection of all glueings in $M$
  of the sets of the form $\setof{a_i\on p_i}{a_i\in A\sss\sss\&\sss\sss p_i\in \Omega
    \sss\sss\&\sss\sss i\in I}$ \emph{which are compatible}.

  Explicitly,    \begin{eqn}
    \begin{split}
      \hbox{cl}(A) &=
      \{ \ssss b\in |M| \ssss : \ssss \exists I \sss\sss \tupof{(a_i,p_i)}{i\in I}\in {^I(A\times \Omega)}
      \\
 &\qquad   (\forall i,\sss j\in I \sss\sss ( Ea_i\on p_i = Ea_j\on p_j = [a_i\on p_i = a_j\on p_j]_M) ) \sss\sss\&\sss\sss \\
 &\qquad   Eb = \Vee Ea_i\on p_i \sss\sss\&\sss\sss \forall i\in I\sss\sss Ea_i\on p_i = [a_i\on p_i = b]_M \ssss \}.
   \end{split}
  \end{eqn}

    Note that if $M$ is a sheaf we can omit the words ``in $M$'' in
    the definition of $\hbox{cl}(A)$ since $M$ \emph{will} contain all
    the glueings of compatible systems.
  \end{definition}

  \begin{lemma}\label{reconstruction_from_dense_in_sheaves} Let
    $\Omega$ be a complete Heyting algebra, $M$, $N$ are sheaves over $\Omega$, with $M$ a
    subsheaf of $N$ and $A\subseteq |M|$ then $\hbox{cl}^M(A)=\hbox{cl}^N(A)$.
\end{lemma}

  \begin{proof} Immediate from the uniqueness of glueings in $M$.
  \end{proof}
  \vskip12pt

With this brief discussion of closure in hand we now turn to the
notion of a subset of a presheaf being dense in another subset. This consequent notion of the
density of a presheaf allows us to bound the cardinality of a presheaf in terms of its density.
That this can be done will be important later. (The precise bound is not so crucial.)

  \begin{definition}\label{dense_in} Let  $\Omega$ be a complete Heyting algebra,
    $M$ a presheaf over $\Omega$ and $D$, $A\subseteq |M|$. Then $D$ is \emph{dense in} $A$ if
for every $a\in A$ there is some $\setof{(d_i,p_i)}{d_i\in D\sss\sss\&\sss\sss p_i\in \Omega
  \sss\sss\&\sss\sss i\in I}$ such that $Ea=\bigvee_{i\in I} p_i$ and
   for all $i\in I$ we have $a\on p_i = d_i \on p_i$.
\end{definition}
    If $M$ is a presheaf over $\Omega$ and
    $D\subseteq A \subseteq |M|$ and $D$ is dense in $A$ then, the definition of
    density immediately gives that $A \subseteq \hbox{cl}(D)$.

    Of course, it may well be that
    $\hbox{cl}(D)\setminus A$ is non-empty, but that is very much to be expected for any
    notion of closure of a dense subset of a set (\emph{e.g}, topologically) unless the
    set itself has good closure properties.

\begin{lemma} Let  $\Omega$ be a complete Heyting algebra,
  $M$ a presheaf over $\Omega$ and $D$, $A\subseteq |M|$. Then $D$ is \emph{dense in} $A$ if and only if for all $a\in A$ we have
  $Ea = \bigvee \setof{[a=d]_M}{d\in D}$.
\end{lemma}

\begin{proof} For the left-to-right direction, suppose $a\in A$ and $\setof{(d_i,p_i)}{d_i\in D\sss\sss\&\sss\sss p_i\in \Omega
  \sss\sss\&\sss\sss i\in I}$ is such that $Ea=\bigvee_{i\in I} p_i$ and
  for all $i\in I$ we have $a\on p_i = d_i \on p_i$. Since for all $d\in D$ we have $[a=d]_M\le Ea$ (by Lemma (\ref{[a=b]_below_Ea_Eb}) we have
  $\bigvee \setof{[a=d]_M}{d \in D}\le Ea$.

  However, for $i\in I$, since $a\on pi = d_i\on p_i$ and $Ea = \bigvee \setof{p_i}{i\in I}$,
  we also have \[ p_i\land [a = d_i]_M = [a\on p_i=d_i\on p_i]_M =  p_i\land Ea = p_i\land Ed_i = p_i .\]
Thus for every $i\in I$ we have $p_i\le [a=d_i]_M$ and hence
\[ Ea = \bigvee\setof{p_i}{i\in I} \le \bigvee\setof{[a=d_i]_M}{i\in I} \le \bigvee \setof{[a=d_i]_M}{d\in D}.\]

Thus we must have $Ea = \bigvee \setof{[a=d]_M}{d\in D}$ as claimed.

For the right-to-left direction, if $a\in D$ and $Ea = \bigvee \setof{[a=d])M}{d\in D}$, then simply enumerate $D$ as $\setof{d_i}{i\in I}$ and for $i\in I$ set
$p_i$ to be $[a=d_i]_M$. Then the family $S= \setof{(d_i,p_i)}{i \in I}$ is as required.
\end{proof}
  
    \begin{corollary}\label{use_of_reconstruction_from_dense_in_sheaves} Let
    $\Omega$ be a complete Heyting algebra and $M$, $N$ and $P$ sheaves over $\Omega$ with
  $|M|$, $|N|\subset |P|$. Let $A\subseteq |M|$ and suppose $D\subseteq A$ is dense in $A$.
  If $D\subseteq |N|$, then $A\subseteq |N|$.
    \end{corollary}

    \begin{proof}  Immediate from Lemma (\ref{reconstruction_from_dense_in_sheaves}) and the definition of dense.
    \end{proof}
\vskip12pt
  
\begin{definition} Let
    $\Omega$ be a complete Heyting algebra, $M$ a presheaf over $\Omega$ or a  sheaf over $\sh{\Omega}$
  and $A\subseteq |M|$. We write $d(A)$ for the \emph{density} of $A$,
  the size of the smallest dense subset of $A$. We write $d(M)$ for $d(|M|)$.
    \end{definition}

   \begin{proposition}\label{bound_size_in_terms_of_density} Let
    $\Omega$ be a complete Heyting algebra. 
     If $M$ is a presheaf over $\Omega$ then $\card{|M|} \le \card{\Omega}^{d(M)}$.
\end{proposition}

\begin{proof} This is immediate from Lemma (\ref{systems_rep_uniquely}).
  For each $a\in M$ and $\setof{(d_i,p_i)}{d_i\in D\sss\sss\&\sss\sss p_i\in \Omega
  \sss\sss\&\sss\sss i\in I}$ such that $Ea=\bigvee_{i\in I} p_i$ and
  for all $i\in I$ we have $a\on p_i = d_i \on p_i$, let $D_a=\setof{d_i}{i\in I}$
  and define $p^a:D\longrightarrow \Omega$ by
  $p^a(d) = p_i$ if $d\in D_a$ and $d=d_i$ and $p^a(d)=\bot$ otherwise.
  Then $Ea=\bigvee_{d\in D} p^a(d)$ and for all $d\in D$ we have $a\on p^a(d) = d\on p^a(d)$.
  Lemma (\ref{systems_rep_uniquely}) shows that $p^a=p^b$ implies $a=b$, so
  the function $p^{.}:M\longrightarrow {{^D}\Omega}$ is injective.
\end{proof}

\vskip6pt

We now show how to amalgamate two presheaves or two sheaves when they are extensions of a common presheaf or sheaf.

  \begin{proposition}\label{amal_for_presheaves} (Amalgamation for presheaves)
    Let $\Omega$ be a complete Heyting algebra, let
    $A$, $M_0$ and $M_1$ be presheaves over $\Omega$
    and let $A \xlongrightarrow{f_0} M_0$ and
    $A \xlongrightarrow{f_1} M_1$ be presheaf morphisms.
    Then there is a presheaf $P$ over $\Omega$ and presheaf monomorphisms
    $M_0  \xlongrightarrow{g_0} P$ and
    $M_1 \xlongrightarrow{g_1} P$ such
    that $g_0\cdot f_0 = g_1\cdot f_1$. Furthermore, if $M_0$ and $M_1$ are sheaves
    over $\Omega$ there is a sheaf $S$ over $\Omega$ and sheaf monomorphisms
    $M_0  \xlongrightarrow{h_0} S$ and
    $M_1 \xlongrightarrow{h_1} S$ such that $h_0\cdot f_0 = h_1\cdot f_1$.
    then so are the $h_i$.
    \end{proposition}
  
  \begin{proof} We do the simplest thing possible: take the amalgamation
    of the underlying sets and impose the obvious restriction and extent
    functions by carrying over those from $M_0$ and $M_1$. We give the easy, but somewhat lengthy,
    details as we did not see them elsewhere in the literature.

    Let $|P|=\setof{(a,i)_{/{\sim}}}{i\in \set{0,1} \sss\sss\ \&\sss\sss a \in |M_i|}$
    where $(a,i)\sim (b,j)$ if and only if $(a,i)=(b,j)$ or
    $i\ne j$ and there is some $c\in A$ such that $a=f_i(c)$ and
    $b=f_j(c)$.

    For $(a,i)_{/{\sim}} \in |P|$ set $E^P (a,i)_{/{\sim}} =E^{M_i} a$. This
    is a good definition, since if $c\in A$ we have, since $f_0$ and $f_1$
    are morphisms, that $E^{M_0} f_0(c) = E^A a = E^{M_1} f_1(c)$.

    For $(a,i)_{/{\sim}} \in |P|$ and $p\in \Omega$ set
    $(a,i)_{/{\sim}}  \on^P p = (a \on^{M_i} p,i)_{/{\sim}}$.
    If $a\in |A|$ and $p\in \Omega$
    then since $A$ is a presheaf we have $a\on^A p\in |A|$.
    Since $f_0$ and $f_1$ are morphisms we have, for each $i\in \set{0,1}$, that
    $f_i(a)\on p =f_i(a\on p)$. Consequently,\begin{multline*}
    (f_0(a),0)_{/{\sim}} \on^{P} p =
    (f_0(a)\on^{M_0} p,0)_{/{\sim}} =
    (f_0(a\on^A p),0)_{/{\sim}} = \\
    (f_1(a\on^A p),1)_{/{\sim}} =
    (f_1(a)\on^{M_1} p,1)_{/{\sim}} =
    (f_1(a),1)_{/{\sim}} \on^{P} p ,
    \end{multline*}
    and we have shown that $\on^P$ is well defined.

    We have to check $E^P$ and $\on^P$ satisfy the
    three presheaf-defining properties.  Let $(a,i)_{/{\sim}} \in |P|$.

    First of all, if $p\in \Omega$ we have
    $E^P((a,i)\on p) = E^{M_i}(a\on p) = E^{M_i} a\wedge p
    = {E^P (a,i)_{/{\sim}}} \wedge p$, where the second equality holds because
    $M_i$ is a presheaf and the third by the definition of $E^P$.
    
    We also have
    $(a,i)_{/{\sim}} \on E^P (a,i)_{/{\sim}} =
    (a \on^{M_i} E^{M_i} a,i)_{/{\sim}}  =
    (a,i)_{/{\sim}} $, where the first equality holds
    by the definition of $E^P$ and the second holds
    since $M_i$ is a sheaf and
    hence $a \on^{M_i} E^{M_i} a = a$.

    Thirdly, if $p$, $q\in \Omega$ then
    \begin{multline*}((a,i)_{/{\sim}} \on^P p)\on^P q =
    ((a \on^{M_i} p,i)_{/{\sim}} )\on^P q =
    ((a \on^{M_i} p)\on^{M_i} q,i)_{/{\sim}} = \\
    ((a \on^{M_i} p\wedge q,i)_{/{\sim}} =
    ((a,i)_{/{\sim}}) \on^P p\wedge q .
    \end{multline*}

    For $i\in \set{0,1}$ and $a\in |M_i|$ set $g_i(a) = (a,i)_{/{\sim}}$. 

    We also have to check for $i\in \set{0,1}$ that
    $g_i$ is a monomorphism. However, if $a\in |M_i|$ and $p\in \Omega$
    then
    $E^P g_i(a) = E^P (a,i)_{/{\sim}} = E^{M_i} a$ by the definition of
    $E^P$, and $g_i(a \on^{M_i} p)= (a\on^{M_i} p,i)_{/{\sim}}
    = (a,i)_{/{\sim}} \on^P p = g_i(a) \on^P p$. Furthermore, if
    $a$, $b\in |M_i|$ and $(a,i)\sim (b,i)$ then $a=b$.

    The commuativity property holds since if $a\in |A|$
    then \[g_0\cdot f_0(a) =
    (f_0(a),0)_{/{\sim}} = (f_1(a),1)_{/{\sim}} = g_1\cdot f_1(a). \]

    Finally, if $M_0$ and $M_1$ are sheaves one can simply sheafify the resultant presheaf
    $P$ to obtain a sheaf $S$ with a sheafification presheaf morphism
    $P\xlongrightarrow{c} S$ and for $i\in\set{0,1}$ take $h_i = c\cdot g_i$.
  
    \end{proof}

  For the remainder of the paper we fix a complete Heyting algebra $\Omega$.

  \section{$L$-structures in  $\psho$}\label{L-structures}

  We begin with Scott's notion of  a characteristic function (\emph{cf.}~\cite{FS}). 
  
\begin{definition} Let $M$ be a presheaf over $\Omega$. 
  $h:|M|^n\longrightarrow \Omega$ is a\break {\emph{characteristic function}} if for $\bar{a}=\tup{a_0,\dots,a_{n-1}}$,
  {$\bar{b}=\tup{b_0,\dots,b_{n-1}} \in |M|^n$}, and writing
  $E\bar{a} $ for $\Wedge_{j<n} Ea_j$, we have
  \[h(\bar{a})\le E\bar{a}\hbox{ and }
  h(\bar{a}) \land\Wedge_{j<n} [a_j=b_j]_M \le h(\bar{b}),\]\vskip-18pt
  or equivalently,    \[ h(\bar{a})\le E\bar{a}\hbox{ and }
      h(\bar{a}) \land\Wedge_{j<n} [a_j=b_j]_M = h(\bar{b}) \land\Wedge_{j<n} [a_j=b_j]_M .
      \]
      Let ${\mathfrak k}_nM = \setof{ h:|M|^n\longrightarrow \Omega}{h\hbox{ is a characteristic function}}$.
\end{definition}

We note for fixed $n$ elements of ${\mathfrak k}_nM$ can be combined, and indeed form a complete Heyting algebra.

\begin{definition}\label{char_fns_form_cHa} (\emph{cf.}~\cite{Mir06}, \S{}37.) Let $M$ be a
  presheaf over $\Omega$ and $n<\omega$. Let $h$, $k\in {\mathfrak k}_nM$.

  Set $h\accentset{\star}\le k$ if for all $\bar{a}\in |M|^n$ we have $h(\bar{a}) \le k(\bar{a})$. 

  For $\diamondsuit\in \set{\land,\lor}$ define $h\accentset{\star}\diamondsuit k$ by for
  all $\bar{a}\in |M|^n$ setting $h\accentset{\star}\diamondsuit k(\bar{a}) = h(\bar{a})\diamondsuit k(\bar{a})$.

  Note, for all $\bar{a}\in |M|^n$ we have $h\accentset{\star}\diamondsuit k(\bar{a}) \le E\bar{a}$ since $h(\bar{a})$, $k(\bar{a})\le E\bar{a}$.

  Define $h\accentset{\star}\to k$ by for
  all $\bar{a}\in |M|^n$ setting $(h\to k)(\bar{a}) = E\bar{a}\land ( h(\bar{a})\to k(\bar{a})) $, and define
   $\accentset{\star}\neg h$ by for
  all $\bar{a}\in |M|^n$ setting $(\neg h)(\bar{a}) = E\bar{a}\land (\neg h(\bar{a})) $.

  The least element of ${\mathfrak k}_nM$ is the map $\bot_{{\mathfrak k}_nM}$, for which for all $\bar{a}\in |M|^n$ we have
  $\bot_{{\mathfrak k}_nM}(\bar{a})=\bot_\Omega$, and
  whilst the greatest element of  ${\mathfrak k}_nM$ is the map
  $\top_{{\mathfrak k}_nM}$, for which for all $\bar{a}\in |M|^n$ we have $\top_{{\mathfrak k}_nM}(\bar{a})=E\bar{a}$.
\end{definition}

\begin{lemma}\label{id_plays_nicely_with_char_fns} (\cite{M88}, Lemma 1.2), 
  Suppose $M$ is a presheaf over $\Omega$ and 
  ${h:|M|^n\longrightarrow \Omega}$ is a characteristic function.
  Then $h(\bar{a}\on E\bar{a}) = h(\bar{a})$.
\end{lemma}

\begin{proof} We have $h(\bar{a}\on E\bar{a}) \land \Wedge_{j<n} [a_j = a_j\on E\bar{a}]_M =
  h(\bar{a}) \land \Wedge_{j<n} [a_j = a_j\on E\bar{a}]_M$. Recall that for any $b\in |M|$ and $p\in \Omega$ we have
  $[b=(b\on p)]_M = p\wedge [b=b]_M = p\land Eb$. 
  Thus
  \begin{eqn}
    \begin{split}
  h(\bar{a}\on E\bar{a}) & = h(\bar{a}\on E\bar{a}) \land \Wedge_{j<n} (E (a_j) \wedge E\bar{a}) =
  h(\bar{a}\on E\bar{a}) \land \Wedge_{j<n} [a_j = a_j\on E\bar{a}]_M \\
  & = h(\bar{a}) \land \Wedge_{j<n} [a_j = a_j\on E\bar{a}]_M = 
  h(\bar{a}) \land \Wedge_{j<n} (E (a_j) \wedge E\bar{a}) \\
  & = h(\bar{a}) \land E\bar{a} = h(\bar{a}).
        \end{split}
  \end{eqn}
\end{proof}

\begin{lemma}
  Let $n<\omega$,  $M \in pSh(\Omega)$ and ${\bar{a} = \tup{a_0, \ldots, a_{n-1}} \in {^n|M|}}$. Let  $\bar{p} = \tup{p_0, \ldots, p_{n-1}} \in {^n\Omega}$  and
  write $\bar{a}\on\bar{p}$ for $\tup{a_0\on p_0,\dots,a_{n-1}\on p_{n-1}}$.
  Let $h : |M|^n \to \Omega$ be a characteristic function. Then,
\[ h(\bar{a} \on \bar{p}) = \bigwedge_{i <n} p_i \wedge h(\bar{a}) .\]
\end{lemma}

\begin{proof} It suffices to prove the lemma when $\tup{p_1,\dots,p_{n-1}} = \tup{Ea_1,\dots,Ea_{n-1}}$ and use induction to prove the lemma in the generality stated.
  Let $\bar{b} = a_0\on p_0 \concat \tup{a_1,\dots,a_{n-1}}$.  Then \[ h(\bar{a}) \land \bigwedge_{i<n} [a_i = b_i ]_M \le h(\bar{b}) \] by the second defining property of being a characteristic function.
  \begin{eqn}\begin{split} h(\bar{a}) \land p_0 = h(\bar{a}) \land p_0 \land E\bar{a} = & h(\bar{a}) \land \bigwedge_{i<n} [a_i = b_i ]_M = \\ & h(\bar{b}) \land \bigwedge_{i<n} [a_i = b_i ]_M
  = h(\bar{b}) \land E\bar{b}= h(\bar{b}). 
    \end{split}
  \end{eqn}
 \end{proof}

\begin{definition}\label{first_oreder_lang}
  A \emph{first order language with equality, $L$, over $\Omega$}
  consists of sets for each $n\in \omega\setminus\set{0}$ of (non-equality)
  $n$-ary relations $\rel(n,L)$ and $n$-ary functions $\fun(m,L)$, constant symbols $C$
which form a presheaf over $\Omega$, and variable symbols $\setof{v_i}{i \in \Nbb}$.
\end{definition}

We now build towards defining the notion of an \emph{$L$-structure in $\psh{\Omega}$},
originally given by \cite{FS}, for $L$ a first order language with equality over $\Omega$. (See also \cite{B16}, \cite{Borceux}, \cite{TvD}.)

We start by defining, by simultaneous induction, terms and their extents. 

\begin{definition}\label{defn_term_and_their_extents}   Let $L$ be a first order language with equality over $\Omega$. 
  \emph{Terms}, $\tau$, in $L$ are strings of symbols and their \emph{extents}, $E\tau$, are elements of $\Omega$. They are
  jointly generated, recursively, by the following three rules 

  \begin{itmz}
  \item    Every variable is a term of $L$ and has extent $\top$.
  \item Every constant of $L$ is a term of $L$ and as a term has extent equal to its extent as a constant.
     \item If $n > 0$, $f$ is an $n$-ary function symbol of $L$ and $\tau_0$,\dots, $\tau_{n-1}$ are terms
       so is  $f(\tau_0,\dots, \tau_{n-1})$, and it has extent
       $E\tau_0\wedge \dots \wedge E\tau_{n-1}$.
   \end{itmz}
\end{definition}
  
\begin{definition} Let $L$ be a first order language with equality over $\Omega$ and $\tau$ a term in $L$.
   \begin{itmz}
     \item        $\FV(\tau)$ is the set of variables occuring (recursively in) $\tau$.
       \item $v(\tau)$ is  $\card{\FV(\tau)}$ 
\item $k(\tau)$ is the   least $k$ such that $\FV(\tau)\subseteq \setof{v_i}{i<k}$. 
   \end{itmz}   
\end{definition}

\begin{definition} Let $L$ be a first order language with equality over $\Omega$ and
  let $\tau_0,\dots,\tau_{n-1}$ be a collection of terms in $L$.
 \begin{itmz}
 \item   $v(\tau_0,\dots,\tau_{n-1})$ is  $\card{\FV(\tau_0)\cup\dots\cup \FV(\tau_{n-1})}$
   \item $k(\tau_0,\dots,\tau_{n-1})$ is the least $k$ such that
     $\FV(\tau_0)\cup\dots\cup \FV(\tau_{n-1}) \subseteq \setof{v_i}{i<k}$.
     \end{itmz}
\end{definition}

\begin{definition} Let $L$ be a first order language with equality over $\Omega$ and $\tau$ a term in $L$. 
  If $v(\tau)=0$ (or, equivalently, $k(\tau)=0$)
  we say $\tau$ is a \emph{closed} term. Clearly, if  $\tau_0,\dots,\tau_{n-1}$ is a collection of terms in $L$
  then $v(\tau_0,\dots,\tau_{n-1})=0$ if and only if all of the $\tau_i$ are closed terms for $i<n$.
\end{definition}

\begin{definition}\label{set_up_for_formulas} (cf.~\cite{Marker})
  Let $L$ be a first order language with equality over $\Omega$.
  \emph{Formulae} are strings of symbols built using the relation, function, constant and
  variable symbols of $L$, the equality symbol $=$, the
  Boolean connectives $\land$, $\lor$, $\longrightarrow$ and $\neg$,
  the quantifiers $\exists$ and $\forall$, and parentheses $($ , $)$,
  as in the classical case. \emph{Sentences} are formulae with no free
  variables.
\end{definition}

\begin{definition}\label{defn_L-theory} Let $L$ be a first order language with equality over $\Omega$. A \emph{theory} (in $L$)
  or an \emph{$L$-theory} is simply a set of $L$-sentences.
\end{definition}
  
\begin{definition}\label{extent_of_formulas} (See \cite{B16}, Definition 1.22.)
  Let $L$ be a first order language with equality over $\Omega$. Let $\phi$ be
 a formula of $L$. Define
\[ E\phi = \bigwedge \setof{E^C c }{c \hbox{ occurs in }\phi}.\]
\end{definition}

\begin{definition}\label{restr_term_or_formula_to_elt_of_Omega}
  Let $L$ be a first order language with equality over $\Omega$.
For $p \in \Omega$ and $\gamma$ a term or a formula of $L$ define $\gamma\on p$, the \emph{restriction of $\gamma$ to $p$},
  as the term or formula obtained by substituting
  every occurrence of a symbol $c\in |C|$ by $c\on p$. 
\end{definition}

\begin{definition}\label{defn_L_str_in_psho} (\emph{cf}.~\cite{B16}, Definition 1.20) Let $L$ be a first order
  language with equality over $\Omega$.
  $M$ is an \emph{$L$-structure in $\psho$} or a \emph{presheaf of
  $L$-structures over $\Omega$} if it is a presheaf over $\Omega$ and there is
\begin{eqn}
    \begin{split}
&\hbox{ a characteristic function $[R(.)]_M:|M|^n\longrightarrow \Omega$ for every $R\in \rel(n,L)$,}\\
&\hbox{ a presheaf morphism $f^M:M^n \longrightarrow M$ for every $f\in \fun(n,L)$, and}\\
&\hbox{  a presheaf morphism $\cdot^M : C\longrightarrow |M|$, with $\cdot^M:c \longrightarrow {\mathfrak c}^M$.}
    \end{split}
\end{eqn}
Note that the equality symbol in $L$ is interpreted by $[.=.]_M$.\end{definition}
\vskip12pt

In the interests of conciseness we refer the reader to \cite{FS}, \cite{Hyland}, \cite{B16}, \cite{Borceux} and \cite{TvD}
for examples of $L$-structures and the other concepts introduced in this section.

\begin{proposition}\label{bound_on_number_of_presheaf_L-structures} Let $L$ be a first order language with equality over $\Omega$
  and let $\kappa$ be an infinite cardinal.
  The collection of isomorphism types of $L$-structures in $\psh{\Omega}$ with density $\kappa$
  is a set (and not a proper class). Hence the collection of isomorphism types
  of subsets of $L$-structures in $\psh{\Omega}$ with density $\kappa$
  is also a set (and not a proper class).
\end{proposition}

\begin{proof} For any $L$-structure $M$ in $\psh{\Omega}$ with $d(M)=\kappa$, by Proposition
  (\ref{bound_size_in_terms_of_density}), the size of
  $M$ is at most $\card{\Omega}^\kappa$. The possible isomorphism types now depends on
  the ways of assigning the appropriate structure to the set. This depends on $L$ and $\Omega$,
  but there is some cardinal bound on the number of ways of doing so.
    \end{proof}

A much more crude, but still useful,  result is the following

\begin{proposition} Let $L$ be a first order language with equality over $\Omega$ and suppose
  $\lambda$ is a strongly inaccessible cardinal.
  The collection of isomorphism types of $L$-structures in $\psh{\Omega}$ with density less than $\lambda$
  is a set (rather than a proper class). Hence the collection of isomorphism types
  of subsets of $L$-structures in $\psh{\Omega}$ with density less than $\lambda$
  is also a set (rather than a proper class).
\end{proposition}

\begin{proof} Immediate, since, by
  Proposition (\ref{bound_size_in_terms_of_density}), if $M\in
  \psh{\Omega}$ then $d(M)<\lambda$ if and only if $\card{M}<\lambda$.
\end{proof}

We now define the interpretation of terms in $L$-structures. Earlier treatments of the details
of interpretation were somewhat cavalier and below we need to work closely with the definitions.
Consequently, we give here (what we hope is) a rigorous account.

In classical model theory there are a number of ways of defining how terms are interpreted in structures.
At the heart of the matter is systematically replacing the free variables by elements of the structure.
However, one has to make provision for how a term could later be used in the process of generating another, more
complex, term or a formula. So, one has to set up a regime under which variables which are not mentioned in the term
can be later interpreted. There are various ways to do this.

One, attractive, way to set up how terms are interpreted in structures
is used in \cite{EFT}. They define the interpretation of a term under
an assignment of \emph{all} of the variables. In one way this is an
extremely clean procedure: if one's structure is $M$ one simply has one
function $\tau^M:M^{\omega}\longrightarrow M$ for each term $\tau$,
defined by induction on the complexity of $\tau$. However, there is a
small set theoretic price to be paid as one ends up with $\tau^M$
having domain of size $\card{M}^\omega$.

It is possibly because other authors prefer to end up with something with domain
of size $\card{M}$ that other approaches are more prevelant.  However,
these too have a cost, this time in terms of keeping track of how a
chain of increasing finite assignments of variables to elements of the
model cohere with each other. Again there are choices.

A finitary version of the approach in \cite{EFT} is to define functions for all
assignements of the first $k$ variables where $k$ is greater than or
equal to the largest index of a variable appearing in the term. This
is at least relatively simple conceptually. It is this approach we
generalize. However, initially one only has a function assigning exactly the
variables appearing in a term if these happen to be precisely the
variables with indices less than $k$ for some $k\in \omega$.
So an auxiliary definition is required to allow us to obtain functions
of the same arity as the number of variables appearing in the term.

Other approaches allow immediately for assignment of exactly the variables
occuring in the term. However, the plethora of all possible finite
supersets then has to be taken into account.

Clearly, either way
one ends up with countably many functions from a set of size the
cardinality of $M$ into $M$. Nevertheless, the questions of coherence appear to us to
be harder to visualize in the second approach. (Although, the proofs
of coherence are not, in fact, so difficult.)

The definitions we give here are similar to one of the classical
definitions just outlined, but with the twist that for closed terms we
need a function from $\Omega$ to $M$ rather than a function from
$\set{\emptyset}$ to $M$ -- or, equivalently, the choice of a single
element of $M$.

At this point we recall Notation (\ref{defn_M^0}), that
whenever $M$ be an $L$-structure in $\psho$ we set $M^0 = \Omega$.

The analogous convention in the classical model theoretic set-up is that $M^0$ is a
set of size $1$ (rather than being the empty set as one might imagine
if trying to mechanically apply the definition of $M^n$ for $n>0$ to
the case $n=0$). From a category theoretic point of view, while $\Omega$, as shown in Proposition (\ref{Omega_terminal}),
is the terminal object in the category $\psho$, a set of size $1$
is terminal in the category of sets.

\vskip12pt

\begin{definition}\label{defn_interp_terms_in_l_str_in_psho} Let $L$ be a first order language with equality over $\Omega$ and
  $M$ be an $L$-structure in $\psho$.
  For each term $\tau$ we define for each $n \ge k(\tau)$ a presheaf morphism
  $\tau^M: M^n\on E\tau \longrightarrow M\on E\tau$.
  $\tau^M$ is the \emph{interpretation of $\tau$ in $M$}.
  The definition is by induction on the complexity of terms.
  \begin{itmz}
  \item  If $\tau$ is a variable $v_i$, then $E\tau = \top$ and
    for all $a\in M^n$ we set $\tau^M(\bar{a})=v_i^M (\bar{a}) = a_i$.
  \item  If $\tau$ is some $c \in |C|$ then for every $\bar{a}\in M^n\on E\tau$ we set
    $c^M(\bar{a}) = {\mathfrak c}^M\on E\bar{a}$. (In fact we could make this definition for all
    $\bar{a} \in M^n$.) 
  \item If $\tau_0 $, \dots $\tau_{m-1}$ are terms in $L$, $f \in \fun(m, L)$ and
    $\tau = f(\tau_0, \dots , \tau_{m-1})$, then for all $n \ge k(\tau_0,\dots,\tau_{m-1})$
    since $E\bar{a} \le E\tau = E\tau_0\wedge \dots\wedge E\tau_{m-1}$ we have
    and thus we can set
    \[ \tau^M (\bar{a}) = f^M(\tau_0^M(\bar{a}), \dots , \tau_{m-1}^M(\bar{a})) .\]
  \end{itmz}

  Note that, consequent to Notation (\ref{defn_M^0}), in the second clause and in the third clause
  if $k(\tau_0,\dots,\tau_{m-1})=0$, when $n=0$ we have that ``$\bar{a}$'' is an element of
  $M^0 = \Omega$ and (so) $E\bar{a} = \bar{a}$.

  In order to see that the third clause of this definition is legitimate we must check, since $\dom(f^M)=M^n$, that 
  $E\tau_0^M(\bar{a})= \dots = E\tau_{m-1}^M(\bar{a})$.

  In order to show this we prove the following lemma follows by induction from the first two clauses of definition.

  \begin{lemma}\label{restrict_range_tau^M} Let $\tau$ be
    a term in $L$ as in the statement of the definition.
  For all $\bar{a}\in M^n \on E\tau$ we have
  $E\tau^M(\bar{a}) = E\bar{a}$.\end{lemma}

\begin{proof} If $\tau$ is a variable $v_i$ we have
  $E\tau^M(\bar{a}) = Ea_i = E\bar{a}$.

  If $\tau$ is a constant $c\in |C|$ we have 
  $E c^M(\bar{a}) = E({\mathfrak c}^M\on E\bar{a}) = E{\mathfrak c}^M\wedge E\bar{a}
  = E^Cc\wedge E\bar{a} = E\tau\wedge E\bar{a} = E\bar{a}$, since $\cdot:C\longrightarrow M$ is a
  presheaf morphism.

  Finally, if $\tau_0 $, \dots $\tau_{m-1}$ are terms in $L$, $f \in \fun(m, L)$ and
  $\tau = f(\tau_0, \dots , \tau_{m-1})$,
  and if, by induction, $E\tau_0^M(\bar{a})= \dots = E\tau_{m-1}^M(\bar{a}) = E\bar{a}$, we have
  \[ E\tau^M (\bar{a}) = Ef^M(\tau_0^M(\bar{a}), \dots , \tau_{m-1}^M(\bar{a})) = E\bar{a}
  ,\]
  since $f^M$ is a presheaf morphism.
  \end{proof}

Moreoever, by a similar inductive proof, in each case $\tau^M$ commutes with restrictions.

  \begin{lemma}\label{tau^M_commutes_with_restriction} Let $\tau$ be
    a term in $L$ as in the statement of the definition, $\bar{a}\in M^n \on E\tau$ and $p\in \Omega$. Then
  $\tau^M(\bar{a})\on p = \tau^M(\bar{a}\on p)$.\end{lemma}

\begin{proof} If $\tau$ is a variable $v_i$ we have
  \[ \tau^M(\bar{a})\on p = v_i^M(\bar{a})\on p = a_i\on p = v_i^M(\bar{a}\on p)=\tau^M(\bar{a}\on p).\]

  If $\tau$ is a constant $c\in |C|$ we have 
  \[ c^M(\bar{a}) \on p =( {\mathfrak c}^M\on E\bar{a}) \on p =  {\mathfrak c}^M\on E\bar{a}\land p
  =  ( {\mathfrak c}^M\on p) \on  E\bar{a} = c^M(\bar{a} \on p) .\]

 Finally, if $\tau_0 $, \dots $\tau_{m-1}$ are terms in $L$, $f \in \fun(m, L)$ and
  $\tau = f(\tau_0, \dots , \tau_{m-1})$,
  and if, by induction, for each $i<n$ we have $\tau_i^M(\bar{a})\on p = \tau_i^M(\bar{a}\on p)$
  \begin{eqn}\begin{split} \tau^M (\bar{a}) \on p & = f^M(\tau_0^M(\bar{a}), \dots , \tau_{m-1}^M(\bar{a})) \on p 
   = f^M(\tau_0^M(\bar{a})\on p, \dots , \tau_{m-1}^M(\bar{a}) \on p) \\
  & = f^M(\tau_0^M(\bar{a}\on p), \dots , \tau_{m-1}^M(\bar{a} \on p)) =  \tau^M (\bar{a} \on p)
  ,
  \end{split}
  \end{eqn}
  where the second equality holds since $f^M$ is a presheaf morphism.
\end{proof}

So we have shown that each $\tau^M$ is a presheaf morphism.

  Now let $v(\tau)\le n\le k(\tau)$ and let $\iota:n \longrightarrow
  k(\tau)$, $\iota:j\mapsto i_j$, be any increasing function such that
  $\FV(\tau) \subseteq \setof{v_{i_j}}{j<n} \subseteq
  \setof{v_l}{l<k(\tau)}$. Define a presheaf morphism $\tau_\iota^M:
  M^n\on E\tau \longrightarrow M\on E\tau$ by $\tau_\iota^M(\bar{a})=\tau^M(\bar{a}) =
  \tau^M(\bar{b})$ where $b\in M^{k(\tau)}\on E\tau$ is any sequence from $M\on E\tau$
  such that for all $e<n$ we have $b_{\iota(e)} = a_e$.
  \end{definition}
   
\begin{lemma} 
  The definition of the $\tau_{\iota}^M$ is a good one. \end{lemma}
\begin{proof} This is immediate from the three clause
  definition for the presheaf morphism $\tau^M:M^{k(\tau)}\on E\tau\longrightarrow M\on E\tau$.
  \end{proof}

\begin{note-nono} Note that without its penultimate paragraph
  this definition would only treat $n \ge k(\tau)$. For comparison we
  give a analogue of Marker's definition for classical model theory from \cite{Marker}.  In this
  definition an auxiliary clause (here given by the paragraph starting
  ``Furthermore'' below) is also essential, in this case in order for
  the definition of the interpretation of a composition of terms under
  a function to work correctly.

  \begin{definition}\label{defn_interp_terms_a_cote_Marker} Let $L$ be a first order language with equality over $\Omega$ and
    $M$ an $L$-structure in $\psho$. 
Let $\tau$ be a term with $\FV(\tau)=\set{v_{i_0},\dots,v_{i_{n-1}}}$.
We define a presheaf morphism
$\tau^M_*$, the \emph{interpretation of $\tau$ in $M$}, with
${\tau^M_*: M^{v(\tau)}\on E\tau \longrightarrow M\on E\tau}$.
  The definition is again by induction on the complexity of terms. Let $\bar{a} = \set{a_0,\dots,a_{n-1}}$.
  \begin{itmz}
  \item  If $\tau$ is a variable $v_{i_j}$, then $E\tau = \top$ and $\tau^M_*(\bar{a})=v_{i_j}^M (\bar{a}) = a_j$.
  \item  If $\tau$ is some $c \in |C|$, then for every $p \in M^0$, $ = \Omega$ with $p\le Ec$ we set
    $c^M_*(p) = {\mathfrak c}^M\on p$.
  \item If $\tau_0 $, \dots $\tau_{m-1}$ are terms in $L$, $f \in \fun(m, L)$ and
    $\tau = f(\tau_0, \dots , \tau_{m-1})$, then if $k(\tau)\ne 0$ we set 
    $\tau^M_* (\bar{a}) = f^M(\tau_{0*}^M(\bar{a}), \dots , \tau^M_{m-1*}(\bar{a}))$, and if the $\tau_i$ are all closed terms
        (\emph{i.e.}~$v(\tau_0,\dots,\tau_{m-1})=0$)
    and $p\in \Omega$ is such that $p\le E\tau$ then we set $\tau^M_* (p) = f^M(\tau_{0*}^M(p), \dots , \tau_{m-1*}^M(p))$.
  \end{itmz} 
  Furthermore, for every $k>v(\tau)$ and every order-preserving injection
  $\iota$ of $v(\tau)$ into $k$
  we define a presheaf morphism
  $\tau_{*\iota}^M = \tau^M_*: M^{k}\on E\tau \longrightarrow M\on E\tau$ as follows:
  \begin{itmz}
    \item If $v(\tau)=0$ then  $\bar{a}\in M^k\on E\tau$ we set 
    $c^M_*(\bar{a}) = {\mathfrak c}^M\on E\bar{a}$.
    \item If $v(\tau)\ne 0$ then $\tau_{*\iota}^M(\bar{a}) =
      \tau^M_*(\tupof{a_j}{j\in rge(\iota)})\on E\bar{a}$.
  \end{itmz} 
\end{definition}
\end{note-nono}

The upshot of this is that (under either approach) we have defined
enough presheaf morphisms to allow us to use $\tau^M$ (or $\tau^M_*$)
to refer both to an instantiation of $\tau$ in $M$ using a minimal set of elements of
$M$ and to instantiations using non-minimal sets of elements of $M$
giving us flexibility to compose terms.

Our official definition will be to use $\tau^M$, although everything would work fine using $\tau^M_*$ instead.

Now we define the interpretation of formulae in $L$-structures in
$\psho$. Each formula's interpretation will be a characteristic
function.

In order to give  this definition it is helpful to have some notation for
dealing with instantiations of variables appearing in subformulae of formulae.

\begin{notation} Let $L$ be a first order language with equality over $\Omega$ and $M$ an $L$-structure in $\psho$.
  Suppose $j:m\longrightarrow n$ is an increasing function, $\bar{a}\in M^n$ and
  $\psi(v_{i_{j(0)}},\dots,v_{i_{j(m-1)}})$, $=\psi(\bar{v}')$, is a subformula of $\phi(v_{i_0},\dots,v_{i_{n-1}})$, $=\phi(\bar{v})$.
  Write $\bar{a}_{\psi(\bar{v}')}$ for $\tupof{a_{j(l)}}{l<m}$, and where there is no confusion
  abbreviate $\bar{a}_{\psi(\bar{v}')}$ as $\bar{a}_\psi$.
\end{notation}

\begin{definition} Let $L$ be a first order language with equality over $\Omega$ and
  let $M$ be an $L$-structure in $\psho$. We
  associate to each formula $\phi(\bar{v})$ in $L$ with $n$-many free variables, where $n\ge 1$
  a characteristic function $[\phi(.)]_M : |M|^n \longrightarrow \Omega$, 
  the \emph{interpretation of $\phi$ in $M$}. We at the same time
  associate to each formula $\phi$ in $L$ with no free variables
  a valuation in $\Omega$ - in this $n=0$-case the reader should simply ignore
  all mentions of $\bar{v}$ and $\bar{a}$.  This definition is also by induction on complexity.

  \begin{itmz}
\item If $\tau_0 (\bar{v} )$, $\tau_{1}(\bar{v})$ are terms in $L$, and letting $p= E\tau_0\wedge E\tau_1$, then
  \[ [ (\tau_0 (\bar{v} ) = \tau_{1}(\bar{v}))(\bar{a})]_M = [\tau_0^M(\bar{a}\on p ) = \tau_{1}^M(\bar{a}\on p)]_M .\]
\item If $\tau_0 (\bar{v} )$, \dots, $\tau_{m-1}(\bar{v})$ are terms in $L$, 
  $q= E\tau_0\wedge\dots\wedge E\tau_{m-1}$ and  $R\in \rel(m,L)$ then
  \[ \qquad\qquad [R(\tau_0 (\bar{v}) ,\dots, \tau_{m-1}(\bar{v}))(\bar{a})]_M =
  [R(\tau_0^M(\bar{a}\on q) , \dots, \tau_{m-1}^M(\bar{a}\on q))]_M.\] 
\item
  If $\phi$ is a formula then
  \[ [\neg \phi(\bar{a})]_M = E\phi \wedge E\bar{a} \wedge \neg [\phi(\bar{a}\on E\phi )]_M \]
\item
  If $\phi$, $\psi$ are formulae with $q = E\phi \wedge E\psi$, then
  for  $\diamondsuit \in \set{\land, \lor,\longrightarrow} $
\[ [(\phi \sss\diamondsuit\ssss \psi)(\bar{a})]_M = q \wedge E\bar{a} \wedge
  ( [\phi(\bar{a}_\phi\on q )]_M \sss\sss\diamondsuit\sss\sss [\psi(\bar{a}_\psi\on q )]_M )
   . \]
\item $[\exists x\phi(x, \bar{a})]_M = \bigvee_{t\in |M|} [\phi(t, \bar{a})]_M$.
\item $[\forall x\phi(x, \bar{a})]_M =
  E\phi \wedge E\bar{a} \wedge \bigwedge_{t\in |M|}( Et \longrightarrow  [\phi(t, \bar{a})]_M)$.
  
\end{itmz}
\end{definition} 
\vskip12pt

The reader should note the asymmetry between clauses for existential and universal quantification in this (inductive)
definition. The definition of interpretation of the universal quantifier combined with the definition of `$\longrightarrow$' in $\Omega$ shows we take
\[ [\forall x\phi(x, \bar{a})]_M =
E\phi \wedge E\bar{a} \wedge \bigwedge_{t\in |M|}( \Vee \setof{ p \in \Omega }{ p \land Et \le  [\phi(t, \bar{a})]_M } ).\]

When we work with $M$ an $L$-structure in
$\psho$ we always have that there is a (unique) section with extent $\bot$, which can be obtained as
$t\on \bot$ for any $t\in |M|$. Consequently, for any $\phi(x,\bar{a})$ we have $\bigwedge_{t\in |M|}  [\phi(t, \bar{a})]_M = \bot$.
The `$Et\longrightarrow$' obviates this and comparable problems.

On the other hand, for any $\phi(x,\bar{a})$ and $t$ with $Et=\top$ we have
\[ Et \longrightarrow  [\phi(t, \bar{a})]_M = \top \longrightarrow  [\phi(t, \bar{a})]_M = [\phi(t, \bar{a})]_M .\]
Thus, if $\Omega$ is the two element complete Heyting algebra $\set{\bot,\top}$, we recover the usual interpretation of universal quantification in classical model theory, since then
\[ [\forall x\phi(x, \bar{a})]_M = \bigwedge_{t\in |M| \sss\sss\&\sss\sss Et\ne \bot} [\phi(t, \bar{a})]_M , \]
and so the definition is an appropriate generalization of the classical one.

Note, also, the contrast with the recent preprint of Chen, \cite{Chen},
where he takes $\forall x$ as an abbreviation for $\neg \exists x\neg$. 
This gives a strikingly different framework.\footnote{Chen works with \'etal\'e structures rather than working directly
with sheaves. However this is not a substantive difference: see, for example, Borceux, \cite{Borceux},
for the well known translations between \'etal\'e spaces and sheaves.}
\vskip12pt

\begin{lemma}\label{interp_of_fmla_bounded_by_extents_of_fmla_and_params}
  Let $L$ be a first order language with equality over $\Omega$,
  $M$ an $L$-structure in $\psho$ and $\phi$ a formula in $L$ with $n$-many free variables.
  If $n\ge 1$ and $\bar{a}\in |M|^n$ then
  $[\phi(\bar{a})]_M \le E\phi \wedge E\bar{a}$. If $n=0$ then $[\phi]_M \le E\phi$.
\end{lemma}

\begin{proof}
  The proof is by induction on the complexity of formulae.  We start
  with atomic formulae. If $\tau_0$ and $\tau_1$ are terms and $\phi$
  is the formula $\tau_0=\tau_1$ and we set $p=E\tau_0\wedge E\tau_1$ and note that
  $p = E(\tau_0 = \tau_1)$ by definition, we have
        \begin{eqn}
    \begin{split}
      &   [(\tau_0 =  \tau_1)(\bar{a})]_M =  [\tau_0^M(\bar{a}\on p) =  \tau_1^M(\bar{a}\on p)]_M =
      [\tau_0^M(\bar{a})\on p =  \tau_1^M(\bar{a})\on p]_M \\
      &   \sss\sss  = \bigvee \setof{E(\tau^M_0(\bar{a} \on p \wedge q) )}{q\in \Omega \sss\sss\&\sss\sss
        \tau^M_0(\bar{a}) \on p \wedge q =  \tau^M_1(\bar{a}) \on p \wedge q}  \\
  & \sss\sss   =  \bigvee \setof{E\bar{a}\wedge p \wedge q }{
        \tau^M_0(\bar{a}) \on p \wedge q = \tau^M_1(\bar{a}) \on p \wedge q} \le E\bar{a},\sss\sss p
    \end{split}
        \end{eqn}where the final equality holds because the $\tau_i^M$ are presheaf morphisms.
              
Similarly, if $R$ is a relation and $\tau_0$, \dots $\tau_{m-1}$ are
terms with $q = \bigwedge_{i<m} E\tau_i$, $ = E(R(\tau_0,\dots,\tau_{m-1}))$, 
since $[R(.)]_M$ is a characteristic function we have
        \begin{eqn}
          \begin{split}
            & [R(\tau_0^M(\bar{a}\on q),\dots,\tau_{m-1}^M(\bar{a}\on q))]_M
                         \le E\bar{a}\on q = Ea \wedge q.
    \end{split}
        \end{eqn}

        The inductive steps for the binary connectives, negation and
        universal quant\discretionary{-}{}{}ification are all immediate as the required
        result is explicit in the definition of interpretation of
        formulae.

        Finally, for the inductive step for existential quantification we have
        \[ [\exists x\phi(x, \bar{a})]_M =
        \bigvee_{t\in |M|} [\phi(t, \bar{a})]_M \le \bigvee_{t\in |M|} Et \wedge E\phi \wedge E\bar{a}
        \le E\phi \wedge E\bar{a}, \] where the first inequality is
        given by applying the inductive hypothesis to $\phi$ and
        $t\concat \bar{a}$.
        \end{proof}

We make a small, but useful, observation on substitution in the interpretation of equality of terms.

\begin{lemma}\label{char_funct}


  Let $\bar{v} = \tup{v_0, \ldots, v_{n-1}}$ be variables and let $\tau( \bar {v})$, $\sigma( \bar {v})$  be $L$-terms. Let $M \in pSh(\Omega, L)$, with
  $\bar{a} = \tup{a_0, \ldots, a_{n-1}}$,  $\bar{b}=\tup{ b_0, \ldots, b_{n-1}} \in {^n|M \on E\tau \wedge E\sigma| }$. Then
  \[ [\bar{a} = \bar{b}]_M \wedge [ \sigma^M (\bar{a}) = \tau^M (\bar{a})]_M = [\bar{a} = \bar{b}]_M \wedge [ \sigma^M (\bar{b}) = \tau^M (\bar{b})]_M .\]

Particularly, $ [\bar a = \bar b]_M \wedge [ \sigma^M (\bar{a}) = \tau^M (\bar{a})]_M \leq  [ \sigma^M (\bar{b}) = \tau^M (\bar{b})]_M$.

\end{lemma}

\begin{proof} Set $q= [\bar{a} = \bar{b}]_M$ 
  and $\bar{a} \on q$ for $\tup{a_0 \on q, \ldots, a_{n-1} \on q}$.
  For all $i <n$ we have $a_i \on [a_i = b_i]_M = b_i \on [a_i = b_i]_M$, so \emph{a fortiori}
  we also have  $a_i \on q = b_i \on q$.
  Using this fact in the third equality and Lemma (\ref{restrict_range_tau^M}) in the second equality in the displayed equations following we have 
\begin{eqn}\begin{split}
    [\bar{a} = \bar{b}]_M \wedge [ \sigma^M (\bar{a}) = & \tau^M (\bar{a})]_M = [ \sigma^M (\bar{a}) = \tau^M (\bar{a}) \on q]_M = \\ & [ \sigma^M (\bar{a}) = \tau^M (\bar{a} \on q) ]_M 
    = [ \sigma^M (\bar{b}) = \tau^M (\bar{b} \on q) ]_M  = \\
    & [\bar{a} = \bar{b}]_M \wedge [ \sigma^M (\bar{b}) = \tau^M (\bar{b})]_M .
  \end{split}
\end{eqn}
\end{proof}

We can now define when a model forces an instantiation of a formula or models a sentence.

\begin{definition}\label{defn_forces} Let $L$ be a first order language with equality over $\Omega$,
  $M$ an $L$-structure in $\psho$, $\phi(\bar{v})$ an $L$-formula and $\bar{a}\in |M|^n$.

    We say $M$ \emph{forces} $\phi$ at $\bar{a}$, written $M\forces
    \phi(\bar{a})$, if $[\phi(\bar{a})]_M = E\phi \wedge E\bar{a}$.

    If $\phi$ is an $L$-sentence we thus have
    $M\forces \phi$ if $[\phi]_M = E\phi$. In this case we 
    say $M$ \emph{models} $\phi$, or \emph{is a model of} $\phi$.

    If $\Lambda(\bar{v})$ is a set of $L$-formula with free variables in $v_0,\dots,v_{n-1}$ and $\bar{a}\in |M|^n$
    then $M\forces \Lambda(\bar{a})$ if $M\forces \phi(\bar{a})$ for each $\phi(\bar{v})\in \Lambda(\bar{v})$.

    If $\Lambda$ has no free variables we say $M$ \emph{is a model of} $\Lambda$. 

  \end{definition}

\begin{triteob} For any $L$ the set of all $L$-structures in $\psho$ is a model of the emtpy set of sentences
(\emph{i.e.}, where we take $\Lambda= \emptyset$).
\end{triteob}

This observation allow us to talk exclusively about collections of models of theories or
of sets of sentences, instead of $L$-structures more generally, without losing any generality.

That the forcing relation is closed under conjunction and disjunction follows from the definition of forcing.

\begin{proposition}\label{forced_closed_under_conj_and_disj} Let $L$ be a first order language with equality over $\Omega$,
  and $M$ an $L$-structure in $\psho$. Let $\phi$, $\psi$ be $L$-formulae and $\bar{a}\in |M|^n$.

  If $M\forces \phi(\bar{a}_\phi)$ and $M\forces \psi(\bar{a}_\psi)$ then $M\forces (\phi \lor \psi)(\bar{a})$ and
    $M\forces (\phi \land \psi)(\bar{a})$.
\end{proposition}

\begin{proof} Note, first of all, by Definition (\ref{extent_of_formulas}), we have
  $E(\phi\land\psi) = E(\phi\lor \psi) = E\phi \wedge E\psi$.
  By the definition of interpretation, for $\diamondsuit \in \set{\land,\lor}$, and letting
  $q = E\phi \wedge E\psi$, 
  we have \[ [(\phi \sss\diamondsuit\ssss \psi)(\bar{a})]_M = E\phi\wedge E\psi \wedge E\bar{a} \wedge
  ( [\phi(\bar{a}_\phi\on q )]_M \sss\sss\diamondsuit\sss\sss [\psi(\bar{a}_\psi\on q )]_M .\]
  By the definition of forcing we have $[\phi(\bar{a}_\phi)]_M = E\phi\wedge E\bar{a}_\phi$
  and $[\psi(\bar{a}_\psi)]_M = E\psi\wedge E\bar{a}_\psi$.

  For ``$\lor$'' we have \begin{eqn}
          \begin{split}
    [(\phi \sss\lor\ssss \psi)(\bar{a})]_M & = E\phi\wedge E\psi \wedge E\bar{a} \wedge
  (E\phi\wedge E\bar{a}_\phi \lor E\psi\wedge E\bar{a}_\psi) \\
   & = ( E\phi\wedge E\psi \wedge E\bar{a} \wedge
    (E\phi\wedge E\bar{a}_\phi)) \lor \\
   & \qquad\qquad
    ( E\phi\wedge E\psi \wedge E\bar{a} \wedge
    (E\psi\wedge E\bar{a}_\psi )) \\
    & = ( E\phi\wedge E\psi \wedge E\bar{a} ) \lor ( E\phi\wedge E\psi \wedge E\bar{a} )  \\
    & =  E\phi\wedge E\psi \wedge E\bar{a}  =  E(\phi\lor \psi) \wedge E\bar{a}  
          \end{split}
        \end{eqn}

  Similarly, but more simply, for ``$\land$'' we have
  \begin{eqn}
          \begin{split}
   [(\phi \sss\land\ssss \psi)(\bar{a})]_M & = E\phi\wedge E\psi \wedge E\bar{a} \wedge
   E\phi\wedge E\bar{a}_\phi \wedge E\psi\wedge E\bar{a}_\psi \\
   & = E\phi\wedge E\psi \wedge E\bar{a} =
  E(\phi\wedge \psi) \wedge E\bar{a} .
          \end{split}
  \end{eqn}
\end{proof}

However, while the forcing relation is closed more generally under intuition\discretionary{-}{}{}istic logic,
it not necessarily closed under classical logic.\footnote{See \cite{Fitting} and \cite{TvD}.}
Consequently, when discussing which collections of sentences are forced by a model,
it is reasonable to introduce classifications of formulae which do not rely on
formulae being reducible to prenex normal form.

There are a number of hierarchies of formulae to be found in the literature.
  For the discussion of preservation theorems in the following sections it is convenient to use
  the class forming operations defined by Fleischmann (\cite{Fle}) and we give a variant
  of one of Fleischmann's hierarchies. For other approaches to stratification
  of the collection of first order formulae see, inter alia, \cite{BNI}, \cite{Burr}, \cite{FuKu} and \cite{SUZ}. 

For the following definitions let $L$ be a first order language with equality over $\Omega$.
  
\begin{definition}\label{defn_cal_E} (\cite{Fle}, Definition 3.2)
  Let $\Gamma$ be a set of $L$-formulae. Define ${\mathcal E}\Gamma$ to be the closure of
  $\Gamma$ under $\land$, $\lor$ and $\exists$.
\end{definition}

\begin{definition}\label{defn_cal_U} (\cite{Fle}, Definition 3.3)
  Let $\Gamma$, $\Psi$ be a set of $L$-formulae. Define ${\mathcal U}(\Gamma,\Psi)$ to be the smallest set
  of formulae containing $\Gamma$, closed under $\land$, $\lor$ and $\forall$ and such that if
  $\psi\in \Psi$ and $\phi\in {\mathcal U}(\Gamma,\Psi)$ then $\psi\longrightarrow \phi \in {\mathcal U}(\Gamma,\Psi)$.
\end{definition}

\begin{definition}\label{defn_E-U_hier} Let $\Delta_0$ be the set of quantifier free $L$-formulae. Set
  ${\mathcal E}_0={\mathcal U}_0=\Delta_0$.
  For $n<\omega$ let ${\mathcal E}_{n+1} = {\mathcal E}{\mathcal U}_n$ and
  ${\mathcal U}_{n+1} = {\mathcal U}({\mathcal E}_n,{\mathcal E}_n)$. 
\end{definition}

\begin{definition}\label{defn_pos_exist} Let ${\mathcal E}^+_{1}$ be ${\mathcal E}\hbox{Atomic}$,
  the positive existential formulae.
\end{definition}

\begin{definition}\label{defn_univ_E} Let $\forall{\mathcal E}_{1}
  = \setof{\forall \bar{x}\sss \phi(\bar{x})}{\phi(\bar{x})\in {\mathcal E}_{1}}$.
\end{definition}

We now unpack the definition of forcing to show what is required for a model to force a $\forall{\mathcal E}_{1}$
sentence.
We will later use this example in the course of proving an analogue of
Tarski's classical theorem on the preservation of $\forall_2$-sentences under unions of chains
(see, \emph{e.g.}, \cite{hod}, Theorem (2.4.4)).

\begin{proposition}\label{characterization_of_forcing_forall_E1}
  Let $L$ be a first order language with equality over $\Omega$,
  $N$ an $L$-structure in $\psho$,
  $\phi(\bar{x},\bar{y})$ a formula in ${\mathcal E}_{1}$ with free variables $\bar{x}\concat\bar{y}$
  and $\bar{a}\in |N|^m$.  Then 
  \[ N\forces \forall \bar{x} \sss\sss \phi(\bar{x},\bar{a})\hbox{ if and only if }
  \sss\sss\forall \bar{b}\in |N|^n \sss\sss E\phi \wedge E\bar{a} \wedge E\bar{b} \sss\sss\le
  [\phi(\bar{b},\bar{a})]_N .\]
  \end{proposition}

\begin{proof} $N\forces \forall \bar{x}\sss \phi(\bar{x},\bar{a})$
  if and only if $[\forall \bar{x}\sss \phi(\bar{x},\bar{a})]_N = E\phi \land E\bar{a}$, by the definition of $\forces$,
  where,  
  $E\phi = \bigwedge\setof{E^{C}c}{c\in C \sss\sss\&\sss\sss c\hbox{ appears in }\phi}$,
      by Definition (\ref{extent_of_formulas}).

      Applying the inductive definition of interpretation of formulae we have
  \[ [\forall \bar{x}\sss \phi(\bar{x},\bar{a})]_N = 
    E\phi \wedge E\bar{a}\wedge \bigwedge_{\bar{b}\in |N|^n} (E\bar{b} \longrightarrow [\phi(\bar{b},\bar{a})]_N)
    .\]
    Thus $N\forces \forall \bar{x}\sss \phi(\bar{x},\bar{a})$ if and only if $\sss\forall\bar{b}\in |N|^n
    \sss (E\phi \wedge E\bar{a} \le E\bar{b} \longrightarrow [\phi(\bar{b},\bar{a})]_N )$.
    The latter holds if and only if
       $\sss\forall\bar{b}\in |N|^n
    \sss (E\phi \wedge E\bar{a} \wedge E\bar{b} \le [\phi(\bar{b},\bar{a})]_N )$,  by Lemma (\ref{adjunction}).
\end{proof}

\vskip24pt

\vskip12pt

We define six separate notions of morphism of presheafs (and of sheafs),
  recapitulating the first two from Definitions (\ref{defn_presheaf_morphism}) and (\ref{defn_presheaf_mono}).

\begin{definition}\label{morphism_notions} Let $L$ be a first-order language with equality  over $\Omega$. 
  Let $M$, $N$ be  $L$-structures in $\psho$. 
  Let $f:|M|\longrightarrow |N|$.
  
  $f$ is a \emph{presheaf morphism} if for all $a\in |M|$ and $p\in \Omega$ we have
  \[ E^N f(a) = E^M a \sss\sss\&\sss\sss f(a\on p) = f(a)\on p.
  \]

  If $f$ is a presheaf morphism it is a \emph{presheaf monomorphism} if additionally $f:|M|\longrightarrow |N|$ is injective.
  If $f$ is a presheaf morphism and $\bar{a}=\tup{a_0,\dots,a_{n-1}}\in |M|^n$ write $f``\bar{a}$ for
  $\tup{f(a_0),\dots,f(a_{n-1})}$.

    $f$ is a \emph{${L}$-morphism} if it is a presheaf morphism and for every $n \geq 1$ and
    \begin{itemize}[itemsep=6pt,topsep=-3pt]
    \item for all $n$-ary relations $R$ in $L$ and     $\bar{a} = \tup{a_1,\dots a_n} \in |M|^n$
        we have \[ [R(\bar{a})]_M  \leq [R(f``\bar{a})]_N ,\]
    \item for all $n$-ary functions $\omega$ in $L$  and  $\bar{a} = \tup{a_1,\dots a_n} \in M^n$ we have
      \[ f ( \omega^{M}(\bar{a}))  =  \omega^{N}(f``\bar{a}) , \]
    \item for all constants $c$ in $L$ we have $f({\mathfrak c}^M) = {\mathfrak c}^N$. 
            \end{itemize}

    $f$ is a \emph{weak $L$-monomorphism} if it is an
    $L$-morphism and for every $n \geq 1$, $n$-ary relation $R$
    in $L$ and $\bar{a} \in |M|^n$ we have $\neg \neg
    [R(\bar{a})]_M\ = \neg\neg [R(f``\bar{a})]_N$.

    $f$ is an \emph{${L}$-monomorphism} or \emph{$L$-embedding} if
      it is a presheaf monomorphism, it is an $L$-morphism and for every $n \geq 1$, $n$-ary
    relation $R$ in $L$ and     $\bar{a} \in |M|^n$ we have
     $ [R(\bar{a})]_M\ = [R(f(\bar{a}))]_N$.

    $f$ is an \emph{elementary ${L}$-monomorphism} or \emph{elementary ${L}$-embedding} if
     $f$ is a $L$-morphism and for every $L$-formula  $\phi(v_0,\dots,v_{n-1})$ and $\bar{a} \in |M|^n$ we have
     $[\phi(\bar{a})]_M = [\phi(f``\bar{a})]_N$.

    We refer to these morphisms as being of type $x$
      where $x$ is one of `psh', `psh-mono', `$L$',
      `wk-$L$-mono', `$L$-mono' and `elt-$L$-emb.' 
\end{definition}

This is also a convenient place for us to define a generalization of the notions of $L$-monomorphism and elementary $L$-monomorphism. We start with the more
widely applicable process of expanding a language $L$ by adding names for elements an $L$-structure in the context of presheaves of models.
\vskip12pt

\begin{definition}\label{defn_L_M} If $M$ is an $L$-structure in $\psho$
  let $\setof{\underaccent{\bar}{a}}{a\in |M|}$  be a presheaf of new constant symbols, with a distinct symbol for each element of $|M|$,
  such that for all $a\in |M|$ and $p\in\Omega$ we have $\underaccent{\bar}{a}\on p = \underline{a\on p}$. Equivalentely,
  the map $a\mapsto \underaccent{\bar}{a}$ is  a presheaf isomorphism.
  Let $L_M = L \cup\setof{\underaccent{\bar}{a}}{a\in |M|}$ and define
  $\tilde{M} = (M,a)_{a\in |M|}$, the natural expansion of $M$ to an $L_M$-structure, where for each $a\in |M|$ we have
  $\underaccent{\bar}{a}^{\tilde{M}} =a$.

  If $N$ is another $L$-structure in $\psho$ and $f:M\longrightarrow N$ is a presheaf morphism
  let $\tilde{N}= (N,f(a))_{a\in |M|})$, where for each $a\in |M|$ we have
  $\underaccent{\bar}{a}^{\tilde{N}} =f(a)$.
\end{definition}

The following lemma gives a converse to the last part of Definition (\ref{defn_L_M}).

\begin{lemma}\label{lemma_expansions_give_psh_morphisms} If $M$, $N$ are $L$-structures in $\psho$ and
  $\tilde{N} = (N,\underaccent{\bar}{a}^{\tilde{N}})_{a\in |M|}$ is an $L_M$-structure,
  then the map $f:M\longrightarrow N$ given by $f(\bar{a}) \mapsto \underaccent{\bar}{a}^{\tilde{N}}$ is a presheaf morphism.
\end{lemma}

\begin{proof} The map given by $\underaccent{\bar}{a} \mapsto \underaccent{\bar}{a}^{\tilde{N}}$ is a presheaf morphism by the definition of
  $\tilde{N}$ being an $L_M$-structure, and the map given by $\underaccent{\bar}{a}^{\tilde{M}} \mapsto \underaccent{\bar}{a}$ is a presheaf isomorphism
  by Definition (\ref{defn_L_M}). By the same definition, for all $a\in |M|$ we have $a = \underaccent{\bar}{a}^{\tilde{M}}$. Thus composing the
    morphisms gives a presheaf morphism from $M$ to $N$.
\end{proof}

  \vskip12pt
\begin{definition}\label{Gamma_L_monomorphism}  Let $M\in\psho$ and let $\Gamma$ be a class of $L_M$ formulae.
  We say a presheaf morphism $f:M\longrightarrow N$ is a $\Gamma$-$L$-monomorphism if for all
  $\phi(\underaccent{\bar}{a}_0,\dots,\underaccent{\bar}{a}_{n-1})\in\Gamma$
    we have $[\phi(\bar{a})]_M = [\phi(f``\bar{a})]_N$.
  \end{definition}

  Observe, $\Gamma$ being a set of $L_M$-sentences in this definition is simply a concise way of
  simultaneously specifying a class of $L$-formulae and for each element of the class a set of instantiations which we would like to consider.



\vskip12pt

    \begin{remark}
      In Definition (3.19) of \cite{BM14} (and Definition (4.1) of \cite{B16}),
      the equivalent of $f``\bar{a}$  is given the meaning
      $\tup{f(a_0)\on E\bar{a},\dots,f(a_{n-1})\on E\bar{a}}$.
      We believe that usage to be less desirable than
      the one here, which is the one typical in mathematics more generally.
      However, in virtue of Lemma (\ref{id_plays_nicely_with_char_fns}), 
      none of the proofs of \cite{BM14} are effected by the choice of either meaning of
       $f``\bar{a}$.
    \end{remark}

    \begin{definition} Let $L$ be a first order language with equality over $\Omega$, and let $M$, $N$
    be $L$-structures in $\psho$.
  If the inclusion map is an $L$-embedding we say $M$ is an \emph{${L}$-substructure} of $N$.
 If it is  an elementary $L$-monomorphism we say $M$ is an \emph{elementary ${L}$-substructure} of $N$.

     In all cases, if $L$ is clear from the context we may drop it from these names. So, for example, if
     $L$ is clear from the context and the inclusion map is an $L$-embedding
     (resp.~an elementary $L$-monomorphism) we say $M$ is
     a \emph{substructure} of $N$ (resp.~\emph{elementary substructure} of $N$).
\end{definition}

    We have the analogue of the classical model theoretic notion of a substructure generated by a set of elements of a structure.
    We do use these notions here, but will do so in \cite{BMPIII}.

    \begin{definition}\label{defn_of_substr_gen_by} Let $L$ be a first order language with equality over $\Omega$, let $M$ an $L$-structure in $\psho$ and suppose
      $A\subseteq |M|$. Write for $\gen{A}$, the \emph{presheaf substructure of} $M$ \emph{generated by} $A$: 
   $\gen{A}^M = \bigcap\setof{Q\subseteq M}{Q\hbox{ is an $L$-structure in $\psho$ and }A\subseteq |Q|}$.
  If $M$ is also a sheaf over $\Omega$,
  write $\gen{A}^M_s$ for $\bigcap\setof{Q\subseteq M}{Q\hbox{ is a an $L$-structure in $\sho$  and }A\subseteq |Q|}$,
  the \emph{sheaf substructure of} $M$ \emph{generated by} $A$. (Here we use the Lemma (\ref{intersection_of_presh_lemma}) to see that
  $\gen{A}$ really is a subpresheaf, resp.,~$\gen{A}_s$ a subsheaf.)
    \end{definition}

    We show this definition is unambiguous in the sense that going to larger ambient structures does not
change the generated presheaf substructure or sheaf substructure, and thus the superscripts in the definition are superfluous.

\begin{lemma}\label{uniqueness_of_gen_by_struct} Let
    $\Omega$ be a complete Heyting algebra, $M$, $N$ $L$-structures in $\psho$  over $\Omega$, with
  $M$ being a substructure of $N$ and $A\subseteq |M|$. Then
  $\gen{A}^M = \gen{A}^N$.  
  If $M$, $N$ are $L$-structures in $\psho$ over $\Omega$, $M$ is a subsheaf of $N$ and $A\subseteq |M|$ then
  $\gen{A}_s^M = \gen{A}_s^N$.  
\end{lemma}

\begin{proof} In each case $M$ is amongst the `$Q$'s whose intersections are taken
  in the definition of $\gen{A}^N$ and $\gen{A}_s^N$ respectively in Definition (\ref{defn_of_substr_gen_by}). 
  \end{proof}

\begin{notation} Let  $\Omega$ be a complete Heyting algebra,  $M$  an $L$-structure in $\psho$ over $\Omega$
  and $A\subseteq |M|$. By Lemma (\ref{uniqueness_of_gen_by_struct}), we simply write
  $\gen{A}$ for $\gen{A}^M$, and if $M$ is an $L$-structures in $\sho$ over $\Omega$ we write
  $\gen{A}_S$ for $\gen{A}_s^M$
\end{notation}

\begin{proposition} Let  $\Omega$ be a complete Heyting algebra,  $M$  an $L$-structure in $\psho$  over $\Omega$
  and $A\subseteq |M|$.   Then
  $|\gen{A}| = \setof{\tau^M(\bar{a}\on p)}{a\in {^nA} \sss\sss\&\sss\sss p \in \Omega \sss\sss\&\sss\sss \tau\hbox{ is a term}}$. 

  If $M$ is also a sheaf over $\Omega$, then
  \[ |\gen{A}_s| =  \setof{ b \in M}{b \hbox{ is the glueing of a set of compatible elements of }|\gen{A}|} .\]
\end{proposition}

\begin{proof} Clearly if $A\subseteq |Q|$, where $Q$ is an $L$-substructure in $\psho$ of $M$,
  then $P = \setof{\tau^M(\bar{a}\on p)}{a\in {^nA} \sss\sss\&\sss\sss p \in \Omega \sss\sss\&\sss\sss \tau\hbox{ is a term}} \subseteq |Q|$. 
  Since the $\tau^M$ are presheaf morphisms
  for each $\bar{a}\in {^nA}$, $p$, $q\in \Omega$ and term $\tau$ we have $\tau^M(\bar{a}\on p)\on q = \tau^M(\bar{a}\on p\land q)\in P$.
  Moreover, $P$ is the underlying set of an $L$-structure in $\psho$.  Thus $P=\gen{A}$.
  
  If $M$ is a sheaf and $Q$ is a subsheaf of $M$ with $A\subseteq |Q|$ then since $Q$ is a subpresheaf of $M$ we have
  $|\gen{A}|\subseteq |Q|$ and the glueing of every compatible set of elements of $\gen{A}$ is an element of $Q$. 
  However, clearly the set of glueings of a set of compatible elements of $|\gen{A}|$ is closed under further glueings.

  We check that this set is closed under instantiation of terms.

  Let $\tau$ be a term. Suppose for each $i<n$ that ${b_i}$ is the glueing of a compatible set $C_i$ of elements of $|\gen{A}|$ and for all
  $i<j<n$ we have $Eb_i=Eb_j$, $=p$, say, with $p\le E\tau$. Since $ \Vee_{a\in C_i} Ea = Eb_i = Eb_j = \Vee_{a\in C_i} Eb_j\land Ea$ we may assume that for all
  $i<j< n$ we have $C_i=C_j$, $=C$, say. Write $\bar{b}$ for $\tup{b_0,\dots,b_{n-1}}$. Now for $a$ $a'\in C$ we have
  $E\tau^M(\bar{b})\on Ea \land E\tau^M(\bar{b})\on Ea'  = p\land Ea \land p \land Ea' = Ea\land Ea'$. Moreover
  \begin{eqn}
    \begin{split} [\tau^M(\bar{b})\on Ea =  E\tau^M(\bar{b})\on Ea'] = & Ea\land Ea'\land [\tau^M(\bar{b})  = E\tau^M(\bar{b})] \\
      = & Ea\land Ea'\land p = Ea\land Ea'.
    \end{split}
  \end{eqn}
  Thus $\tau^M(\bar{b})\in \gen{A}_s$.
  \end{proof}

In Definition (\ref{defn_of_psh_gen_by}) we gave a definition of the presheaf generated by a set of element of a presheaf
and of the sheaf genereated by a set of elements of a sheaf.
These definitions can be thought of as the instance of Definition (\ref{defn_of_substr_gen_by}) in the case of $L$ being the empty language with equality.

 \vskip12pt   
We now state and prove a result that will be important for our goals in \S{}\ref{accessibility}.
 
\begin{proposition}\label{coherence} 
Let $L$ be a first order language with equality over $\Omega$.
  Suppose that
  $M$, $N$ and $P\in \psh{\Omega,L}$ (resp.~$\sh{\Omega,L}$), that
  $|M| \subseteq |N| \subseteq |P|$. 
  Then the following hold:
  \begin{enumerate}[align=parleft,labelsep=1cm, label=(\arabic*),itemsep=6pt,topsep=0pt] 
  \item If the inclusions $N\inclusion P$ and $M\inclusion P$ are presheaf morphisms then
    so is the inclusion $M\inclusion N$.
  \item If the inclusion $N\inclusion P$ is an $L$-monomorphism and the inclusion $M\inclusion P$ is an
    $L$-morphism, a weak $L$-monomorphism or an $L$-monomorphism then the inclusion $M\inclusion N$ is also,
    respectively, an $L$-morphism, a weak $L$-monomorphism or a $L$-monomorphism.
  \item If $\Gamma$ is a class of $L_M$-formulae and the inclusions $M\inclusion P$  and $N\inclusion P$ are $\Gamma$-$L$-monomorpisms
     then the inclusion $M\inclusion N$ is also a $\Gamma$-$L$-monomorpism. In particular, if the inclusions $M\inclusion P$ and
     $N\inclusion P$ are elementary $L$-embeddings then the inclusion $M\inclusion N$ is also an elementary $L$-embedding.
\end{enumerate} 
\end{proposition}

\begin{proof}  First of all, suppose the inclusions $M\inclusion P$ and $N\inclusion P$ are a presheaf morphisms. Let
  $a\in |M|$ and $p\in \Omega$. By the former hypothesis we have $E^Ma=E^Pa$ and $a\on^M p = a\on^P p$;
  by the latter, since $|N|\subseteq |P|$, $E^Pa=E^Na$ and $a\on^P p = a\on^N p$. Hence
  $E^Ma=E^Na$ and $a\on^M p = a\on^N p$, and thus the inclusion $M\inclusion N$ is also a presheaf morphism.

  Now suppose $N\inclusion P$ is an $L$-monomorphism and $M\inclusion P$ is an 
  \mbox{$L$-morphism}, a weak $L$-monomorphism or an $L$-monomorphism.

  For the interpretation of constants and function symbols there is little to be proved.
  If $c$ is a constant symbol in $L$ then, since both inclusions into $P$ are $L$-morphisms, we have
  ${\mathfrak c}^M = {\mathfrak c}^P = {\mathfrak c}^N$. If $f$ is an $n$-ary function symbol in $L$ and $\bar{a}\in |M|^n$ then,
  again since both inclusions into $P$ are $L$-morphisms,
  $f^M(\bar{a}\on E\bar{a}) = f^P(\bar{a}\on E\bar{a}) = f^N(\bar{a}\on E\bar{a})$.

  Now suppose $R$ is an $n$-ary relation symbol in $L$ and $\bar{a}\in |M|^n$. Since
  $|M|^n\subseteq |N|^n$, we have $\bar{a}\in |N|^n$ and so, since the inclusion of $N$ into $P$ is an $L$-monomorphism,
  $[R(\bar{a})]_N = [R(\bar{a}\on E\bar{a})]_P$.
  Moreover, $[R(\bar{a})]_N = {[R(\bar{a}\on E\bar{a})]_N}$, and $[R(\bar{a})]_P = [R(\bar{a}\on E\bar{a})]_P$,
  by Lemma (\ref{id_plays_nicely_with_char_fns}).

  However, we also have $[R(\bar{a})]_M \le [R(\bar{a}\on E\bar{a})]_P$,
  hence $[R(\bar{a})]_M \le [R(\bar{a})]_N = [R(\bar{a}\on E\bar{a})]_N$. So $M\inclusion N$ is an $L$-morphism.

  If $\bar{a}\in |M|^n$, we may also have either
  $\neg\neg [R(\bar{a})]_M = \neg\neg [R(\bar{a}\on E\bar{a})]_P$
  or $[R(\bar{a})]_M = [R(\bar{a}\on E\bar{a})]_P$, depending on the case. In the first case we have
  $\neg\neg [R(\bar{a})]_M =
  \neg\neg [R(\bar{a})]_N = \neg\neg [R(\bar{a}\on E\bar{a})]_N$, and in the second 
  $[R(\bar{a})]_M = [R(\bar{a})]_N = [R(\bar{a}\on E\bar{a})]_N$. 
  Hence, $M\inclusion N$ is a weak $L$-monomorphism or
  an $L$-monomorphism if $M\inclusion P$, respectively, is.

  Lastly, suppose $\Gamma$ is a class of $L_M$-formulae and $M\inclusion P$,
  $N\inclusion P$ are $\Gamma$-$L$-embeddings. Let $\phi(\underaccent{\bar}{a}_0,\dots,\underaccent{\bar}{a}_{n-1})\in \Gamma$.

  Since $|M|\subseteq |N|$ and $N\inclusion P$ is an elementary, we have $[\phi(\bar{a})]_N = [\phi(\bar{a})]_P$.

  By the hypotheses on $M\inclusion P$ we have 
  $[\phi(\bar{a})]_M = [\phi(\bar{a})]_P$. Hence
  $[\phi(\bar{a})]_M = [\phi(\bar{a})]_N$, and we have $M\inclusion N$ is a $\Gamma$-$L$-monomorphism.
  \end{proof}

For the remainder of the paper we fix $L$ a first order language with equality over $\Omega$,

\section{Preservation phenomenon}\label{pres_phenom_section}

We will be concerned below with the various types of morphisms between models of theories defined in the previous section.
In particular, we will be concerned with
    situations in which we have a morphism $f:M\longrightarrow N$ of presheaves and know $M\forces T$ for some
    set of sentences $T$ and want to know whether the further properties of $f$, for example, being an $L$-morphism,
    an $L$-monomorphism and so on, ensure that $N\forces T$ as well. 
    We refer to a positive answer, generically, as a \emph{preservation} phenomenon. In general a positive answer will
    depend on the properties of $T$ as well as of the morphisms.
 
    As a step to understanding preservation phenomena, we summarize
    some information from \cite{BM14} about the preservation of
    sentences being forced under various types of morphisms.  We
    restrict attention here initially to the case of sentences in the
    language $L$ for which $M$, $N\in \psho$ are
    $L$-structures. However, analogously to Theorems (\ref{thm_diagram_equivs_L-morph}),
    (\ref{thm_diagram_equivs_weak_L-monomorph}) and  (\ref{thm_m-o-diagram_equivs}) and
    Proposition (\ref{graded_m-o-d_thm}) below, the full strength results
    of \cite{BM14} treat the extended language
    $L\cup\setof{\underaccent{\bar}{a}}{a\in |M|}$ with additional constant
    symbols for each element of $|M|$ and moreover give biconditional
    equivalents of various types of morphism.
    
    \begin{proposition} (See \cite{BM14},\S{}\S{}4,5.) Suppose $f:M\longrightarrow N$ is
      a presheaf morphism. Let $\phi$ be a sentence in $L$ and suppose $M\forces \phi$. Let $D$ be a dense set of
      elements of $\Omega$. The following conditions
      are sufficient for us to conclude that $N\forces \phi$ also.
      
      \begin{tabular}{l l}        
condition on morphism $f$ & form of formula $\phi$ \\
\midrule\addlinespace[8pt]
$L$-morphism (see Lemma (\ref{pres_atomic})) & atomic \\\addlinespace[4pt]
weak $L$-monomorphism & atomic or negation of atomic\\\addlinespace[4pt]
$L$-monomorphism & atomic or negation of atomic\\\addlinespace[2pt]
elementary $L$-embedding & arbitrary $L$-sentence\\\addlinespace[2pt]
\end{tabular}

\end{proposition}   
We next prove some lemmas concerning the preservation of interpretations of formulae under
$L$-morphisms. The lemmas lead us ultimately to two
results on $L$-monomorphisms which extend those of \cite{BM14},
Corollary (\ref{mono_pres_qf_forcing_both_ways}) and Corollary (\ref{mono_pres_exist_forcing_up}). 
However, they will also be useful below in their own right.

\begin{lemma}\label{pres_interp_of_terms}
Suppose $f:M\longrightarrow N$ is an $L$-morphism, $\tau$ is a term
  and $\tau^M$, $\tau^N$ are interpretations of $\tau$ in $M$ and $N$, respectively. If $\bar{a}\in \dom( \tau^M)$ then
  $f(\tau^M(\bar{a})) = \tau^N(f``\bar{a})$. If $p\in \Omega$ is such that $p\in \dom(\tau^M)$ then $f(\tau^M(p)) = \tau^N(p)$.
\end{lemma}

\begin{proof} We work by induction on the structure of $\tau$ and so there are five cases to consider.

  First of all, suppose $\tau = v_i$ for some $i$. Then we have $\tau^N(f``\bar{a}) =f(a_i) = f(\tau^M(\bar{a}))$.

  Secondly, suppose $c\in |C|$. We have $c^M(\bar{a}) = {\mathfrak c}^M\on E\bar{a}$,
  $c^N(\bar{a}) = {\mathfrak c}^N\on E\bar{a}$ and by definition, since
  $f$ is an $L$-morphism, ${\mathfrak c}^N = f({\mathfrak c}^M)$ and $Ef``\bar{a}=E\bar{a}$.
  Hence ${\mathfrak c}^N\on Ef``\bar{a} = f({\mathfrak c}^M)\on Ef``\bar{a} = f({\mathfrak c}^M\on E\bar{a})$.
  
  Furthermore, if $p\le Ec$ we have
    $c^M(p) = {\mathfrak c}^M \on p$, $c^N(p) = {\mathfrak c}^N\on p$ and ${\mathfrak c}^N\on p = f({\mathfrak c}^M)\on p
    = f({\mathfrak c}^M\on p)$.

  Fourthly, if $\tau = g(\tau_0,\dots, \tau_{m-1})$ then, since $f$ is an $L$-morphism, we have 
  \[ f(\tau^M(\bar{a})) = f( g^M(\tau_0^M(\bar{a}), \dots , \tau_{m-1}^M(\bar{a}) )) =
  g^N(f(\tau_0^M(\bar{a})), \dots, f(\tau_{m-1}^M(\bar{a}))).\] However, inductively,
  we have for each $j<m$ that $f(\tau_j^M(\bar{a})) = \tau_j^N(f``\bar{a})$, and hence we have that
  $f(\tau^M(\bar{a})) = g^N(\tau_0^N(f``\bar{a}), \dots, \tau_{m-1}^M(f``\bar{a})) = \tau^N(f``\bar{a})$.

  Finally,
        if the $\tau_i$ are all closed terms 
        then for $p\le E\tau$ we have
        $\tau^N (p) = g^N(\tau_0^N(p), \dots , \tau_{m-1}^N(p))$.
        Yet, for each $j<m$ we have
        $\tau^N_j(p) =   f(\tau_j^M(p))$. So we have
  \begin{eqn}
    \begin{split} \tau^N (p)& = g^N(\tau_0^N(p), \dots , \tau_{m-1}^N(p)) = g^N(f(\tau_0^M(p),\dots,f(\tau_{m-1}^M(p))\\
&        = f(g^M(\tau_0^M(p),\dots,\tau_{m-1}^M(p))) = f(\tau^M(p).
    \end{split}
  \end{eqn}
\end{proof}

\begin{lemma}\label{pres_atomic} Suppose
  $f:M\longrightarrow N$ is an $L$-morphism, $\phi(\bar{v})$ is atomic and $\bar{a}\in |M|^n$. Then
    $[\phi(\bar{a})]_M \le [\phi(f``\bar{a})]_N$.
  \end{lemma}
  
  \begin{proof}
    First of all we deal with the atomic formulae of the form $\tau_0 = \tau_1$ where
    $\tau_0$, $\tau_1$ are terms over $L$. There are six separate cases. In each case the property of $f$ used in
    the proof is that it is a presheaf morphism and hence Lemma (\ref{psh_morphism_push_up_equality}) can be applied.

    Case 1. $\tau_0 = v_i$, $\tau_1 = v_j$ for some $i$, $j$, with $i$ and $j$ not necessarily different.

    We have
  \begin{eqn}
    \begin{split}
        [\phi(\bar{a})]_M =   [(v_i = v_j)(\bar{a})]_M  & =  [v_i^M(\bar{a}\on Ev_i \wedge Ev_j) = 
          v_j^M(\bar{a}\on Ev_i \wedge Ev_j) ]_M \\
        & = [v_i^M(\bar{a}) = v_j^M(\bar{a})]_M = [a_i = a_j]_M.
    \end{split}
  \end{eqn}Similarly, $ [\phi(f``\bar{a})]_N = [f(a_i) = f(a_j)]_N $.
  Finally,
  \[ [a_i = a_j]_M \le [f(a_i) = f(a_j)]_N  \hbox{, 
    since $f$ is an presheaf morphism.}\]

  Case 2. $\tau_0 = v_i$, $\tau_1 = c$ for some $i$ and some $c\in |C|$.

    We have
  \begin{eqn}
    \begin{split}
        [\phi(\bar{a})]_M & =   [(v_i = c)(\bar{a})]_M   = [v_i^M(\bar{a}\on Ev_i \wedge Ec) = 
          c^M(\bar{a}\on Ev_i \wedge Ec) ]_M  \\
        & = [v_i^M(\bar{a}\on Ec) = {\mathfrak c}^M\on E\bar{a}]_M = [a_i\on Ec = {\mathfrak c}^M\on E\bar{a}]_M.
    \end{split}
  \end{eqn}Similarly, since $Ef``\bar{a} = E\bar{a}$,
  we have
  \[
    [\phi(f``\bar{a})]_N = [f(a_i)\on Ec = {\mathfrak c}^N\on E\bar{a}]_N 
   = [f(a_i\on Ec) = f({\mathfrak c}^M\on E\bar{a})]_N, \]
  since  $f$ is a presheaf morphism and an $L$-morphism. Finally, since $f$ is an presheaf morphism,
  by Lemma (\ref{psh_morphism_push_up_equality}) we have
    \[ [a_i\on Ec = {\mathfrak c}^M\on E\bar{a}]_M \le [f(a_i\on Ec) = f({\mathfrak c}^M\on E\bar{a})]_N
  \]

  Case 3. $\tau_0 = d$, $\tau_1 = c$ for some $c$, $d \in |C|$. Let $p = Ec\wedge Ed $.

  We have
  \begin{eqn}
    \begin{split}
      [\phi]_M  & = [c=d]_M = [{\mathfrak c}^M \on p = {\mathfrak d}^M \on p]_M  \\
      &  \le [f({\mathfrak c}^M \on p) = f({\mathfrak d}^M \on p)]_M ,
      \hbox{ by Lemma (\ref{psh_morphism_push_up_equality}), } \\
      &  = [{\mathfrak c}^N\on p = {\mathfrak d}^N\on p]_N = [c=d]_N = [\phi]_N.
    \end{split}
  \end{eqn}

  Case 4. $\tau_0 = v_i$ and there are $\sigma_0$, \dots, $\sigma_{m-1}$ and a function $g$ such that
  $\tau_1 = g(\sigma_0,\dots,\sigma_{m-1})$. Note $E\tau_1 = E\sigma_0\wedge\dots\wedge E\sigma_{n-1}$.

  We have 
  \begin{eqn}
    \begin{split}
      [\phi(\bar{a})]_M & = [(v_i = g(\sigma_0,\dots,\sigma_{m-1}))(\bar{a})]_M \\
     & = [a_i\on E\tau_1 = 
       g^M(\sigma_0^M(\bar{a}\on E\tau_1),\dots,\sigma_{m-1}^M(\bar{a}\on E\tau_1))]_M  \\
     & \le [f(a_i \on E\tau_1)  = 
         f( g^M(\sigma_0^M(\bar{a}\on E\tau_1),\dots,\sigma_{m-1}^M(\bar{a}\on E\tau_1)) )]_N \\
     & = [f(a_i) \on E\tau_1 = 
         g^N ( f(\sigma_0^M(\bar{a}\on E\tau_1)),\dots,f(\sigma_{m-1}^M(\bar{a}\on E\tau_1))) ]_N \\
     & = [f(a_i) \on E\tau_1 = 
         g^N ( \sigma_0^N(f``\bar{a}\on E\tau_1),\dots,\sigma_{m-1}^N(f``\bar{a})\on E\tau_1) ]_N \\
      & \qquad \qquad \qquad\qquad \qquad\qquad \qquad \qquad\qquad \qquad
      \hbox{(by Lemma (\ref{pres_interp_of_terms}))} \\
     & = [v_i =  g(\sigma_0,\dots,\sigma_{m-1}))(f``\bar{a})]_N = [\phi(f``\bar{a})]_N 
    \end{split}
  \end{eqn}

  Case 5. $\tau_0 = c\in |C|$ and there are $\sigma_0$, \dots, $\sigma_{m-1}$ and a function $g$ such that
  $\tau_1 = g(\sigma_0,\dots,\sigma_{m-1})$. Let $p = E\tau_0\wedge E\tau_1$, 
  $= Ec\wedge E\sigma_0\wedge\dots\wedge E\sigma_{n-1}$.

  We have
  \begin{eqn}
    \begin{split}
      [\phi(\bar{a})]_M  & =   [(c = g(\sigma_0,\dots,\sigma_{m-1}))(\bar{a})]_M \\
      & = [{\mathfrak c}^M \on p = 
       g^M(\sigma_0^M(\bar{a}\on p),\dots,\sigma_{m-1}^M(\bar{a}\on p))]_M  \\
      & \le [f({\mathfrak c}^M \on p) = 
         f( g^M(\sigma_0^M(\bar{a}\on p),\dots,\sigma_{m-1}^M(\bar{a}\on p))))]_N  \\
     & = [f({\mathfrak c}^M ) \on p = 
         g^N ( f(\sigma_0^M(\bar{a}\on p)),\dots,f(\sigma_{m-1}^M(\bar{a}\on p))) ]_N  \\
     & = [{\mathfrak c}^N \on p)  =
          g^N ( \sigma_0^N(f``\bar{a}\on p),\dots,\sigma_{m-1}^N(f``\bar{a}\on p))]_N \\
      &  \qquad \qquad \qquad\qquad \qquad\qquad \qquad 
      \hbox{(by Lemma (\ref{pres_interp_of_terms}))}\\
     & = [(c = g(\sigma_0,\dots,\sigma_{m-1}))(f``\bar{a})]_N = [\phi(f``\bar{a})]_N 
    \end{split}
  \end{eqn}

  Similarly, if all the $\sigma_j$ for $j<m$ are closed terms and then we have  
  \begin{eqn}
    \begin{split}
     [\phi]_M & =   [c = g(\sigma_0,\dots,\sigma_{m-1})]_M \\
     &  = [{\mathfrak c}^M \on p = 
       g^M(\sigma_0^M ,\dots,\sigma_{m-1}^M) \on p ]_M  \\
     & \le [f({\mathfrak c}^M \on p) = 
         f( g^M(\sigma_0^M(p),\dots,\sigma_{m-1}^M(p)))]_N  \\
     & = [f({\mathfrak c}^M ) \on p = 
         g^N ( f(\sigma_0^M(p)),\dots,f(\sigma_{m-1}^M(p))) ]_N  \\
     & = [{\mathfrak c}^N \on p  =  g^N ( \sigma_0^N(p),\dots,\sigma_{m-1}^N(p))]_N \\
     & = [(c = g(\sigma_0,\dots,\sigma_{m-1}))]_N = [\phi]_N 
    \end{split}
  \end{eqn}

  Case 6. $\tau_0 = \tau_1$ and there are $\sigma_0$, \dots, $\sigma_{m-1}$, $\rho_0$, \dots $\rho_{k-1}$
  and functions $g$ and $h$
  such that $\tau_0 = h(\rho_0,\dots,\rho_{k-1})$ and $\tau_1 = g(\sigma_0,\dots,\sigma_{m-1})$. Let
  $p = E\tau_0\wedge E\tau_1$, $= E\rho_0\wedge\dots\wedge E\rho_{k-1}\wedge E\sigma_0\wedge\dots\wedge E\sigma_{m-1}$.
   
  We have
  \begin{eqn}
    \begin{split}
      [\phi(\bar{a})]_M & =[(h(\rho_0,\dots,\rho_{k-1}) =
        g(\sigma_0,\dots,\sigma_{m-1}))(\bar{a})]_M  \\
    & = [h^M(\rho_0^M(\bar{a}),\dots,\rho_{k-1}^M(\bar{a})) \on  p = 
       g^M(\sigma_0^M(\bar{a}),\dots,\sigma_{m-1}^M(\bar{a})) \on p ]_M  \\
     & \le [f( h^M(\rho_0^M(\bar{a}),\dots,\rho_{k-1}^M(\bar{a})) \on  p ) = \\
     & \qquad    f( g^M(\sigma_0^M(\bar{a}),\dots,\sigma_{m-1}^M(\bar{a}))\on p )]_N  \\
     & = [h^N ( f(\rho_0^M(\bar{a})),\dots,f(\rho_{k-1}^M(\bar{a}))) \on p = \\
     & \qquad    g^N ( f(\sigma_0^M(\bar{a})),\dots,f(\sigma_{m-1}^M(\bar{a}))) \on p)]_N  \\
     &  = [g^N ( \rho_0^N(f``\bar{a}),\dots,\rho_{k-1}^N(f``\bar{a})) \on p = \\
          & \qquad g^N ( \sigma_0^N(f``\bar{a}),\dots,\sigma_{m-1}^N(f``\bar{a})) \on p]_N
      \hbox{, by Lemma (\ref{pres_interp_of_terms}),} \\
     &  =
      [(g(\rho_0,\dots,\rho_{k-1}) = g(\sigma_0,\dots,\sigma_{m-1}))(f``\bar{a})]_N = [\phi(f``\bar{a})]_N .
    \end{split}
  \end{eqn}

  Similarly, if all the $\sigma_j$ for $j<m$ and all the $\rho_e$ for $e<k$ are closed terms,
  \begin{eqn}
    \begin{split}
      [\phi]_M  & =
      [h(\rho_0,\dots,\rho_{k-1}) =
        g(\sigma_0,\dots,\sigma_{m-1})]_M  \\
     & = [h^M(\rho_0^M(p),\dots,\rho_{k-1}^M(p))  = \\
     & \qquad  g^M(\sigma_0^M (p),\dots,\sigma_{m-1}^M(p))  ]_M  \\
     & \le [f( h^M(\rho_0^M(p),\dots,\rho_{k-1}^M(p)) ) = \\
     & \qquad    f( g^M(\sigma_0^M(p),\dots,\sigma_{m-1}^M(p)) )]_N  \\
     & = [h^N ( f(\rho_0^M(p)),\dots,f(\rho_{k-1}^M(p))) = \\
     & \qquad      g^N ( f(\sigma_0^M(p)),\dots,f(\sigma_{m-1}^M(p))) )]_N  \\
     & = [h^N ( \rho_0^N(p),\dots,\rho_{k-1}^N(p))  = \\
          & \qquad     g^N ( \sigma_0^N(p),\dots,\sigma_{m-1}^N(p)) ]_N,
\hbox{by Lemma (\ref{pres_interp_of_terms}),} 
      \\
      & =
     [(h(\rho_0,\dots,\rho_{k-1}) = g(\sigma_0,\dots,\sigma_{m-1}))]_N = [\phi]_N .
    \end{split}
  \end{eqn}
  
  This completes the proof for atomic formulae of the form $\tau_0=\tau_1$.

  Now suppose $R$ is a relation symbol in $L$. 
  We need to consider atomic formulae of the form $R(\tau_0,\dots,\tau_{n-1})$ where for $i<n$ the $\tau_i$ are terms.
  Here we also need that $f$ is an $L$-morphism.
  
    By definition we have \[ [R(\tau_0,\dots,\tau_{n-1})(\bar{a})]_M = 
  [R(\tau_0^M(\bar{a}\on q) , \dots, \tau_{m-1}^M(\bar{a}\on q))]_M,\] where
  $q= E\tau_0\wedge\dots\wedge E\tau_{m-1}$.

  Since $f$ is an $L$-morphism we have that 
  \begin{eqn}
    \begin{split}
      & \; [R(\tau_0^M(\bar{a}\on q) , \dots, \tau_{m-1}^M(\bar{a}\on q))]_M  \\
     \le & \; [R(f(\tau_0^M(\bar{a}\on q)) , \dots, f(\tau_{m-1}^M(\bar{a}\on q)) )]_N \\
     = & \; [R(f(\tau_0^M(\bar{a})\on q) , \dots, f(\tau_{m-1}^M(\bar{a})\on q) )]_N  \\
     = & \; [R(\tau_0^N(f``\bar{a})\on q) , \dots, \tau_{m-1}^N(f``\bar{a})\on q) )]_N \hbox{, using Lemma (\ref{pres_interp_of_terms}), }
              \\
      = & \; [R(\tau_0,\dots,\tau_{n-1})(f``\bar{a})]_N . 
    \end{split}
  \end{eqn}
        \end{proof}

  \begin{lemma}\label{pres_atomic_mono}
    Suppose $f:M\longrightarrow N$ is an $L$-monomorphism, $\phi(\bar{v})$ is atomic and $\bar{a}\in |M|^n$. Then
    $[\phi(\bar{a})]_M = [\phi(f``\bar{a})]_N$.
  \end{lemma}

  \begin{proof} The proof is almost identical to the proof of the previous lemma (Lemma (\ref{pres_atomic}), except that
    in the six cases concerning equality between terms we use the fact
    that $f$ is a presheaf monomorphism and so Proposition
    (\ref{equiv_presheaf_mono}) is applicable, where previously we
    applied Lemma (\ref{psh_morphism_push_up_equality}), and in
    dealing with relations we use the defining property regarding relations of $L$-monomorphism rather than
    that of $L$-morphisms. The result is that we preserve equality throughout.
  \end{proof}

  \begin{lemma}\label{pres_interp_by_neg_ineq}
    Suppose $f:M\longrightarrow N$ is an presheaf morphism. Let $\phi(\bar{v})$ be a formula
    in $L$ and $\bar{a}\in |M|^n$. 
    If $[\phi(\bar{a})]_M \le [\phi(f``\bar{a})]_N$, then
       $[\neg\phi(\bar{a})]_M \ge [\neg \phi(f``\bar{a})]_N$.
  \end{lemma}
    \begin{proof}
  By hypothesis $[\phi(\bar{a})]_M \le [\phi(f``\bar{a})]_N$, so
    $\neg [\phi(\bar{a})]_M \ge \neg [\phi(f``\bar{a})]_N$. \mbox{Using} the two clauses of the definition of
      presheaf morphism, we also have \mbox{$E(f``\bar{a}\on E\phi) = E(\bar{a}\on E\phi)$}. Hence
      \begin{eqn}
        \begin{split}
          [\neg \phi(f``\bar{a})]_N & = E(f``\bar{a}\on E\phi) \land \neg [\phi(f``\bar{a}\on E\phi)]_N \\
          & \le E(\bar{a}\on E\phi) \land \neg [\phi(\bar{a})]_M = [\neg \phi(\bar{a})]_M .
        \end{split}
      \end{eqn}    
  \end{proof}

  \begin{lemma}\label{pres_interp_by_neg}
    Suppose $f:M \longrightarrow N$ is a presheaf morphism. Let $\phi(\bar{v})$ be a formula
    in $L$ and $\bar{a}\in |M|^n$. 
    If $[\phi(\bar{a})]_M = [\phi(f``\bar{a})]_N$, then
       $[\neg\phi(\bar{a})]_M = [\neg \phi(f``\bar{a})]_N$.
  \end{lemma}
    \begin{proof} The proof is similar to that of the previous lemma, but with inequality replaced by equality.
  By hypothesis $[\phi(\bar{a})]_M = [\phi(f``\bar{a})]_N$, so
    $\neg [\phi(\bar{a})]_M = \neg [\phi(f``\bar{a})]_N$. \mbox{Using} the two clauses of the definition of
      presheaf morphism, we also have \mbox{$E(f``\bar{a}\on E\phi) = E(\bar{a}\on E\phi)$}. Hence
      \begin{eqn}
        \begin{split}
          [\neg \phi(f``\bar{a})]_N & = E(f``\bar{a}\on E\phi) \land \neg [\phi(f``\bar{a}\on E\phi)]_N \\
          & = E(\bar{a}\on E\phi) \land \neg [\phi(\bar{a})]_M = [\neg \phi(\bar{a})]_M .
        \end{split}
      \end{eqn}    
  \end{proof}

  \begin{lemma}\label{up_pres_interp_by_conj_disj}
    Suppose $f:M \longrightarrow N$ is an presheaf morphism. Let $\diamond \in\set{\lor, \land}$.
    Let $\phi(\bar{v})$, $\psi(\bar{v})$ be formulae in $L$ and $\bar{a}\in |M|^n$.
    Suppose $[\phi(\bar{a})]_M \le [\phi(f``\bar{a})]_N$ and $[\psi(\bar{a})]_M \le [\psi(f``\bar{a})]_N$.
    Then $[(\phi \diamond \psi)(\bar{a})]_M \le [(\phi \diamond \psi)(\bar{a})]_N$.
  \end{lemma}

  \begin{proof} Let $q=E\phi \land E\psi$. Then, again using the definition of a presheaf morphism,
    $E(f``\bar{a} \on q) = E(\bar{a} \on q)$, $=r$, say. Then we have
      \begin{eqn}
        \begin{split}
          [\phi \diamond \psi(f``\bar{a})]_N &
          = r \land ([\phi(f``\bar{a}\on q)]_N \diamond [\psi(f``\bar{a}\on q)]_N )  \\
          & = (r \land [\phi(f``\bar{a}\on q)]_N) \diamond (r\land [\psi(f``\bar{a}\on q)]_N)  \\
          & \ge (r \land  [\phi(\bar{a}\on q)]_M) \diamond (r \land [\psi(\bar{a}\on q)]_M))  \\
          & = r \land  ([\phi(\bar{a}\on q)]_M \diamond [\psi(\bar{a}\on q)]_M) =
          [(\phi \diamond \psi)(\bar{a})]_M .
        \end{split}
      \end{eqn}          
    \end{proof}

    \begin{lemma}\label{pres_interp_by_bin_conns}
    Suppose $f:M \longrightarrow N$ is an presheaf morphism. Let $\diamond \in\set{\lor, \land,\to}$.
    Let $\phi(\bar{v})$, $\psi(\bar{v})$ be formulae in $L$ and $\bar{a}\in |M|^n$.
    Suppose $[\phi(\bar{a})]_M = [\phi(f``\bar{a})]_N$ and $[\psi(\bar{a})]_M = [\psi(f``\bar{a})]_N$.
    Then $[(\phi \diamond \psi)(\bar{a})]_M = [(\phi \diamond \psi)(\bar{a})]_N$.
  \end{lemma}

    \begin{proof} Let $q=E\phi \land E\psi$. Then, once more using the definition of a presheaf morphism,
      \begin{eqn}
        \begin{split}
          [\phi \diamond \psi(f``\bar{a})]_N &
          = E(f``\bar{a} \on q) \land [\phi(f``\bar{a}\on q)]_N \diamond [\psi(f``\bar{a}\on q)]_N = \\
          & E(\bar{a} \on q) \land  [\phi(\bar{a}\on q)]_M \diamond [\psi(\bar{a}\on q)]_M =
          [(\phi \diamond \psi)(\bar{a})]_M .
        \end{split}
      \end{eqn}    
      \end{proof}
    Note that `$\to$' is one of the connectives treated in this last lemma, whereas in the previous
    previous lemma it is not.

    \vskip12pt
    
  \begin{proposition}\label{mono_conserves_qf}
    Suppose $f:M\longrightarrow N$ is an $L$-monomorphism,
    $\phi(\bar{v})$ is quantifier free and $\bar{a}\in |M|^n$. Then
    $[\phi(\bar{a})]_M = [\phi(f``\bar{a})]_N$.
  \end{proposition}
    
  \begin{proof}
   The proof is by induction on the structure of $\phi$. Lemma (\ref{pres_atomic_mono}) deals with atomic formulae.
   Since $f$ is an $L$-morphism inductive steps involving the binary connectives (${\lor, \land,\to}$) can
   be handled by Lemma (\ref{pres_interp_by_bin_conns}) and those involving negation can be handled by
   Lemma (\ref{pres_interp_by_neg}). 
  \end{proof}

  \begin{corollary}\label{mono_pres_qf_forcing_both_ways}
    Suppose $f:M\longrightarrow N$ is an $L$-monomorphism,
    $\phi(\bar{v})$ is quantifier free and $\bar{a}\in |M|^n$. Then
    $M\forces \phi(\bar{a})$ if and only if $N\forces
    \phi(f``\bar{a})$.
  \end{corollary}

  \begin{proof} By the proposition, $[\phi(\bar{a})]_M = [\phi(f``\bar{a})]_N$, whilst
    $Ef``\bar{a}=E\bar{a}$ since $f$ is a presheaf morphism. So
    $[\phi(\bar{a})]_M = E\phi \wedge
    E\bar{a}$ if and only if  $[\phi(f``\bar{a})]_N = E\phi \wedge Ef``\bar{a}$, as required.
  \end{proof}

  We next give preservation results for formulae involving existential quantifiers.

\begin{lemma}\label{pres_up_exist_quant} Suppose $f:|M| \longrightarrow |N|$ is a function. 
  Let $\phi$ be a formula in $L$ with $n+1$ free variables and $\bar{a}\in |M|^n$.
  \begin{itmz}
    \item If for all $b\in |M|$ we have
      $[\phi(b,\bar{a})]_M \le [\phi(f(b),f``\bar{a})]_N $, then
      $[\exists t\sss\phi(t,\bar{a})]_M \le [\exists t\sss \phi(t,f``\bar{a})]_N $.
  \item  If $f$ is a surjection and  for all $b\in |M|$ we have
      $[\phi(b,\bar{a})]_M = [\phi(f(b),f``\bar{a})]_N $, then
      $[\exists t\sss\phi(t,\bar{a})]_M = [\exists t\sss \phi(t,f``\bar{a})]_N $.
    \end{itmz}

    \end{lemma}

    \begin{proof}
      For the first statement, using the hypothesis for the first inequality and that $f``|M| \subseteq |N|$ for the second, we have
      \begin{eqn}
        \begin{split}
      [\exists t \sss \phi(t,\bar{a})]_M & = \Vee_{b\in |M|} [\phi(b,\bar{a})]_M \le 
      \Vee_{b\in |M|} [\phi(f(b),f``\bar{a})]_N \\
      & \le \Vee_{d\in |N|} [\phi(d,f``\bar{a})]_N =  [\exists t\sss \phi(t,f``\bar{a})]_N .    
        \end{split}
      \end{eqn}
     For the second statement, using the hypothesis for the second equality and that $f``|M| = |N|$ for the third, we have
      \begin{eqn}
        \begin{split}
      [\exists t \sss \phi(t,\bar{a})]_M & = \Vee_{b\in |M|} [\phi(b,\bar{a})]_M =
      \Vee_{b\in |M|} [\phi(f(b),f``\bar{a})]_N \\
      & = \Vee_{d\in |N|} [\phi(d,f``\bar{a})]_N =  [\exists t\sss \phi(t,f``\bar{a})]_N .    
        \end{split}
      \end{eqn}
     \end{proof}

  \begin{corollary}\label{pres_positive_existential_up}
    Suppose $f:M\longrightarrow N$ is an $L$-morphism, $\phi(\bar{v})\in {\mathcal E}^+_1$ is positive existential
    and $\bar{a}\in |M|^n$. Then
    $[\phi(\bar{a})]_M \le [\phi(f``\bar{a})]_N$.
  \end{corollary}

  \begin{proof} The proof is by induction on the structure of $\phi$. Lemma (\ref{pres_atomic})) deals with
    atomic formulae and the previous lemma (Lemma (\ref{pres_up_exist_quant})) deals with existential quantification.

    Finally, suppose $\diamond \in\set{\lor, \land}$.
    Let $\chi(\bar{v})$, $\psi(\bar{v})$ be formulae in ${\mathcal E}^+_1$ and $\bar{a}\in |M|^n$.
    Suppose, by induction,
    $[\chi(\bar{a})]_M \le [\chi(f``\bar{a})]_N$ and $[\psi(\bar{a})]_M \le [\psi(f``\bar{a})]_N$.
    Then by Lemma (\ref{up_pres_interp_by_conj_disj}) we have
    $[\chi \diamond \psi(f``\bar{a})]_N  \le  [(\chi \diamond \psi)(\bar{a})]_M$.
  \end{proof}
  \vskip12pt
  
  \begin{proposition}\label{pres_existential_up_mono}
    Suppose $f:M\longrightarrow N$ is an $L$-monomorphism,
    $\phi(\bar{v})\in {\mathcal E}_1$ and $\bar{a}\in |M|^n$. Then
    $[\phi(\bar{a})]_M \le [\phi(f``\bar{a})]_N$.
  \end{proposition}

  \begin{proof} This is immediate from Proposition (\ref{mono_conserves_qf}), since again
    existential  quantification can again be handled by Lemma
    (\ref{pres_up_exist_quant}) and conjunction and disjunction can be handled by
    Lemma (\ref{up_pres_interp_by_conj_disj}).
  \end{proof}

  \begin{corollary}\label{mono_pres_exist_forcing_up}
    Let $f:M\longrightarrow N$ be an $L$-monomorphism, $\phi(\bar{v}) \in {\mathcal E}_1$ and $\bar{a}\in |M|^n$. Then
    $M \forces \phi(\bar{a})$ implies $N\forces \phi(f``\bar{a})$.
  \end{corollary}

  \begin{proof}     By Proposition (\ref{pres_existential_up_mono}),
    $[\phi(\bar{a})]_M \le [\phi(f``\bar{a})]_N$. However,
    by hypothesis $M\forces \phi(\bar{a})$, so we have $[\phi(\bar{a})]_M = E\phi \wedge E\bar{a}$, and hence 
$E\phi \wedge E\bar{a} \le [\phi(f``\bar{a})]_N $. However, by
Lemma (\ref{interp_of_fmla_bounded_by_extents_of_fmla_and_params}),
 $[\phi(f``\bar{a})]_N \le E\phi \wedge Ef``\bar{a} =E\phi \wedge E\bar{a} $, so we also have
$[\phi(f``\bar{a})]_N = E\phi \wedge Ef``\bar{a}$, and thus $N\forces \phi(f``\bar{a})$. 
  \end{proof}

  With the results proved so far in hand we can improve on some of the
  results of \cite{BM14}, cited previously, concerning the \emph{method of diagrams},
  introduced in classical model theory by Robinson. 

  \begin{definition}\label{defn-Gamma-diagram}
    Let $M$ be an $L$-structure in $\psho$. For $\Gamma$ a set of ${\mathcal L}_M$ sentences, the \emph{$\Gamma$-diagram of ${\mathcal M}$}
  is the function $[.]_{\tilde{M}}:\Gamma\longrightarrow \Omega$. When $\Gamma$ is the set of all ${\mathcal L}_M$ sentences which are either
  atomic or negations of atomic we call this function the \emph{(Robinson) diagram of $\mathcal M$}, and when
  $\Gamma$ is the set of all ${\mathcal L}_M$ sentences we call it the \emph{elementary (Robinson) diagram of $\mathcal M$}.
  When $\Omega$ is the algebra $\set{\top,\bot}$ these are
  characteristic functions of the usual $\Gamma$-diagram, Robinson diagram and elementary Robinson diagram, respectively.
  \end{definition}

  \begin{proposition}\label{graded_m-o-d_thm} Let ${\mathcal M}$ and ${\mathcal N}$ be $L$-structures in $\psho$.
    The following are equivalent.
  \begin{eqn}\begin{split} 
      (1)&\;\; \hbox{There is a $\Gamma$-$L$-monomorphism } f:M\longrightarrow N. \\
      (2)&\;\; \hbox{There is an expansion of }N \hbox{ to an }L_M\hbox{-structure }\tilde{N}\hbox{ such that }\\
      &\;\;\;\; [.]_{\tilde{N} } \on \Gamma = [.]_{\tilde{M}} \on \Gamma .
    \end{split}
  \end{eqn}    
  \end{proposition}

  \begin{proof} (1) implies (2). Let $f:M \longrightarrow N$ be a $\Gamma$-$L$-monomorphism.
   Thus, if $\phi(\underaccent{\bar}{a}_0,\dots,\underaccent{\bar}{a}_{n-1})\in\Gamma$
   we have $[\phi(\bar{a})]_M = [\phi(f``\bar{a})]_N$. Let $\tilde{M} = (M,a)_{a\in |M|}$ and $\tilde{N} = (N,f(a))_{a\in |M|}$.
   We also have $[\phi(\bar{a})]_M = [\phi(\underaccent{\bar}{\bar{a}}^{\tilde{M}})]_M$  and
   $[\phi(f``\bar{a})]_N = [\phi(\underaccent{\bar}{\bar{a}}^{\tilde{N}})]_N$.
   Thus $[\phi(\underaccent{\bar}{\bar{a}}^{\tilde{M}})]_M = [\phi(\underaccent{\bar}{\bar{a}}^{\tilde{N}})]_N$. Hence,
   $\tilde{N}$ is an expansion of $N$ to and $L_M$-structure with $[.]_{\tilde{N} } \on \Gamma = [.]_{\tilde{M}} \on \Gamma$.
 
   (2) implies (1). Define $f:M \longrightarrow N$ by $f(a) =  \underaccent{\bar}{a}^{\tilde{N}}$.
   By Lemma (\ref{lemma_expansions_give_psh_morphisms}),
   $f$ is a presheaf morphism. If
    $\phi(\underaccent{\bar}{a}_0,\dots,\underaccent{\bar}{a}_{n-1}) \in \Gamma$ then 
    $ [\phi(\underaccent{\bar}{a}_0^{\tilde{M}},\dots,\underaccent{\bar}{a}_{n-1}^{\tilde{M}})]_{\tilde{M} }=
    [\phi(\underaccent{\bar}{a}_0^{\tilde{N}},\dots,\underaccent{\bar}{a}_{n-1}^{\tilde{N}})]_{\tilde{N} }$,
    and so $ [\phi({a}_0,\dots,{a}_{n-1})]_{{M} } =
    [\phi(f({a}_0),\dots,f({a}_{n-1}))]_{{N} }$. Thus $f$ is a $\Gamma$-$L$-monomorphism.
  \end{proof}
  
  \begin{definition} Define
    $\Sent^{M}_{0} = \setof{\phi}{\phi\hbox{ is a quantifier-free } L_M\hbox{-sentence}}$ and
    $\Sent^{M} = \setof{\phi}{\phi\hbox{ is a } L_M\hbox{-sentence}}$. Define
    $\Rel^{M}  = \setof{ R(\bar{\underaccent{\bar}{a}})}{R\hbox{ is a }\allowbreak\hbox{relation in }L,\sss\bar{a}\in |M|^n}$.
 \end{definition}
  
  \begin{theorem}\label{thm_m-o-diagram_equivs} Let $M$ and $N$ be $L$-structures in $\psho$. The following are equivalent.
  \begin{eqn}\begin{split}
     1(a)&\;\;\hbox{There is an }L\hbox{-monomorphism } f:M\longrightarrow N .\\
      1(b)&\;\;\hbox{There is an expansion of }N\hbox{ to an }L_M\hbox{-structure } \tilde{N}\hbox{ with } \\
      &\;\;\;\;\;\;      [.]_{\tilde{N}} \on \hbox{$\Sent_0^M$} = [.]_{\tilde{M}} \on \hbox{$\Sent_0^M$}   .  \\  
1(c)&\;\;\hbox{There is an expansion of }N\hbox{ to an }L_M\hbox{-structure }\tilde{N}\hbox{ with } \\
&\;\;\;\;\;\;      [.]_{\tilde{N}} \on \hbox{$\Rel^M$} = [.]_{\tilde{M}} \on \hbox{$\Rel^M$}   .  \\[6pt]
      &\;\;\;\; \hbox{and} \\[6pt]
     2(a)&\;\;\hbox{There is an elementary }L\hbox{-monomorphism } f:M\longrightarrow N .\\
      2(b)&\;\;\hbox{There is an expansion of }N\hbox{ to an }L_M\hbox{-structure } \tilde{N}\hbox{ with } \\
      &\;\;\;\;      [.]_{\tilde{N}} \on \hbox{$\Sent^M$} = [.]_{\tilde{M}} \on \hbox{$\Sent^M$}   .  
    \end{split}
  \end{eqn}
  \end{theorem}

  \begin{proof} Immediate from Proposition (\ref{graded_m-o-d_thm}), using Proposition (\ref{pres_existential_up_mono}) to show
    1(c) implies 1(b).
  \end{proof}

  \begin{remark} Theorem (\ref{thm_m-o-diagram_equivs}) generalizes Robinson's original results on the method of diagrams as these
    are exactly what the theorem gives in the case of $\Omega$ being the 2 element (Boolean) algebra $\set{\top,\bot}$. (\emph{Cf.}, for example,
    \cite{CK}, Chapter 2, or \cite{Kirby}, Chapter 13.1.)
  \end{remark}
  
  \begin{lemma}\label{forcing_preservation_lemma} Let $M$ and $N$ be $L$-structures in $\psho$,
    let $\tilde{M}$ be the natural expansion of $M$ to an $L_M$-structure
    (as before) and let $\tilde{N}$ be an expansion of $N$ to an
    $L_M$-structure. Let $\phi(\bar{\underaccent{\bar}{a}})$ be an
    $L_M$-sentence. Let $q=[\phi(\bar{\underaccent{\bar}{a}}^{\tilde{M}})]_{\tilde{M}}$, $p=[\phi(\bar{\underaccent{\bar}{a}}^{\tilde{N}})]_{\tilde{N}}$ and $d=p\to q$.
    \begin{eqn}\begin{split}
        (1)&\;\; \tilde{M} \forces \phi(\bar{\underaccent{\bar}{a}}\on q).\\
        (2)& \;\; \tilde{N} \forces \phi(\bar{\underaccent{\bar}{a}}\on q)\hbox{ if and only if } q\le p \hbox{, \emph{i.e.} }
        [\phi(\bar{\underaccent{\bar}{a}}^{\tilde{M}})]_{\tilde{M}} \le [\phi(\bar{\underaccent{\bar}{a}}^{\tilde{N}})]_{\tilde{N}}.\\
        (3)& \;\;\hbox{If } \tilde{N} \forces  \phi(\bar{\underaccent{\bar}{a}}\on q),\sss \neg\phi(\bar{\underaccent{\bar}{a}}\on q)
        \hbox{ then }   \neg\neg [\phi(\bar{\underaccent{\bar}{a}}^{\tilde{M}})]_{\tilde{M}} = \neg\neg [\phi(\bar{\underaccent{\bar}{a}}^{\tilde{N}})]_{\tilde{N}}.\\
        (4)& \;\; \hbox{If  }q\le p \hbox{ then }  (E\bar{a}\on p) \land d = q = [\phi(\bar{a}\on p)]_M.
    \end{split}
  \end{eqn}
  \end{lemma}

  \begin{proof} (See also the proof of \cite{BM14}, Theorem (5.8).)

    (1). We have $ [\phi(\bar{\underaccent{\bar}{a}}^{\tilde{M}}\on q)]_{\tilde{M}} =
    [\phi(\bar{\underaccent{\bar}{a}}^{\tilde{M}})]_{\tilde{M}} \land [\phi(\bar{\underaccent{\bar}{a}}^{\tilde{M}})]_{\tilde{M}}  =
    [\phi(\bar{\underaccent{\bar}{a}}^{\tilde{M}})]_{\tilde{M}}$. However, we also have $[\phi(\bar{\underaccent{\bar}{a}}^{\tilde{M}})]_{\tilde{M}} \le
    E\phi(\bar{\underaccent{\bar}{a}}^{\tilde{M}})$, by Lemma (\ref{interp_of_fmla_bounded_by_extents_of_fmla_and_params}), while
    $E\phi(\bar{\underaccent{\bar}{a}}^{\tilde{M}}\on q) = E\phi(\bar{\underaccent{\bar}{a}}^{\tilde{M}}) \land q\le q =
    [\phi(\bar{\underaccent{\bar}{a}}^{\tilde{M}})]_{\tilde{M}}$.

    (2). If $\tilde{N} \forces \phi(\bar{\underaccent{\bar}{a}}\on q)$ then $q = E\bar{\underaccent{\bar}{a}}\on q \land E\phi =
    [\phi(\bar{\underaccent{\bar}{a}}^{\tilde{N}}\on q)]_{\tilde{N}} = q \land [\phi(\bar{\underaccent{\bar}{a}}^{\tilde{N}})]_{\tilde{N}}$.
    If $[\phi(\bar{\underaccent{\bar}{a}}^{\tilde{M}})]_{\tilde{M}} \le [\phi(\bar{\underaccent{\bar}{a}}^{\tilde{N}})]_{\tilde{N}}$ then
    $[\phi(\bar{\underaccent{\bar}{a}}^{\tilde{N}})\on q)] = [\phi(\bar{\underaccent{\bar}{a}}^{\tilde{N}})]_{\tilde{N}} \land q = q = q\land q =
    [\phi(\bar{\underaccent{\bar}{a}}^{\tilde{M}}\on q)]_{\tilde{M}} = E\bar{\underaccent{\bar}{a}}^{\tilde{M}}\on q \land E\phi =
     E\bar{\underaccent{\bar}{a}}^{\tilde{N}}\on q \land E\phi$, where the penultimate equality holds by (1).
    
        (3). By (2) one immediately has
        $\neg\neg [\phi(\bar{\underaccent{\bar}{a}}^{\tilde{M}})]_{\tilde{M}} \le  \neg\neg [\phi(\bar{\underaccent{\bar}{a}}^{\tilde{N}})]_{\tilde{N}}$.

        Applying (2) for $\neg\phi$ in place of $\phi$, one also has
        $[\neg\phi(\bar{\underaccent{\bar}{a}}^{\tilde{M}})]_{\tilde{M}} \le [\neg\phi(\bar{\underaccent{\bar}{a}}^{\tilde{N}})]_{\tilde{N}}$, that is
        $E\bar{a} \land \neg[\phi(\bar{\underaccent{\bar}{a}}^{\tilde{M}})]_{\tilde{M}} \le E\bar{a}\land \neg[\phi(\bar{\underaccent{\bar}{a}}^{\tilde{N}})]_{\tilde{N}}$.
        Taking the $\land$ with $\neg\neg[\phi(\bar{\underaccent{\bar}{a}}^{\tilde{N}})]_{\tilde{N}}$ on both sides give
        $E\bar{a} \land \neg[\phi(\bar{\underaccent{\bar}{a}}^{\tilde{M}})]_{\tilde{M}} \land \neg\neg[\phi(\bar{\underaccent{\bar}{a}}^{\tilde{N}})]_{\tilde{N}} =\bot$,
        and hence one has
        $E\bar{a} \land \neg\neg[\phi(\bar{\underaccent{\bar}{a}}^{\tilde{N}})]_{\tilde{N}} \le  \neg\neg[\phi(\bar{\underaccent{\bar}{a}}^{\tilde{M}})]_{\tilde{M}} $.
        However, $[\phi(\bar{\underaccent{\bar}{a}}^{\tilde{N}})]_{\tilde{N}} \le E\bar{a}$, so 
        $\neg\neg[\phi(\bar{\underaccent{\bar}{a}}^{\tilde{N}})]_{\tilde{N}} \le \neg \neg[\phi(\bar{\underaccent{\bar}{a}}^{\tilde{M}})]_{\tilde{M}} $.

        (4.) Using modus ponens (Lemma \ref{modus_ponens}) for the second equality, the hypothesis that $q\le p$ for the third and penultimate equalities,
        and (1) for the fourth equality, we have
        \begin{eqn}\begin{split} (E\bar{a}\on p ) \land d & = E\bar{a}\on (p\land d) = E\bar{a}\on (p\land q) = E\bar{a} \on q  = [\phi(\bar{a}\on q)]_M\\
            & = [\phi(\bar{a})]_M \land q = q \land q = q = q\land p = [\phi(\bar{a}\on p)]_M.
          \end{split}
        \end{eqn}

  \end{proof}
\vskip12pt

  \begin{definition} Let $M$ be an $L$-structure in $\psho$. 
    Let \begin{eqn}
      \begin{split} \Delta^+_M     = & \sss\sss \setof{\phi}{\phi \hbox{ is an atomic }L_M\hbox{-sentence and }\tilde{M}\forces \phi }\hbox{ and }\\
        \Delta_M= & \sss\sss \setof{\phi}{\phi \hbox{ is an atomic }L_M\hbox{-sentence or a negation of }\\
          &\qquad \qquad\hbox{an atomic } L_M\hbox{-sentence}\hbox{ and }\tilde{M}\forces \phi}.
      \end{split}
    \end{eqn} 
\end{definition}
  
  \begin{definition}Let
    ${\mathcal E}^{M}_1$ be the collection of $L_M$-sentences obtained by closing
    the quantifier-free $L_M$ formulae under the operator $\mathcal E$ (see Definition (\ref{defn_cal_E})).
  \end{definition}

\begin{theorem}\label{thm_diagram_equivs_L-morph} (See \cite{BM14}, Lemma (4.6), showing (1) and (3) are equivalent.)
  Let $M$ and $N$ be $L$-structures in $\psho$.  The following are equivalent.
  \begin{eqn}\begin{split}
      (1)&\hbox{There is an }L\hbox{-morphism } f:M\longrightarrow N .\\
      (2)&\hbox{There is an expansion of }N\hbox{ to an }L_M\hbox{-structure } \tilde{N}\hbox{ with }
\tilde{N}\forces {\mathcal E}^+_1.  \\  
(3)&\hbox{There is an expansion of }N\hbox{ to an }L_M\hbox{-structure }\tilde{N}\hbox{ with }
\tilde{N}\forces \Delta^+_M
    \end{split}
  \end{eqn}
  \end{theorem}

\begin{proof} The equivalence of (2) and (3) is follows from Corollary (\ref{pres_positive_existential_up}).
  The proof that (3) implies (1) follows immediately from Lemma (\ref{forcing_preservation_lemma}(2)).
  To see that (1) implies (3) note that if $\phi(\bar{\underaccent{\bar}{a}}) \in  \Delta^+_M$ then
  $[\phi(\bar{a})] = E\bar{a}\land E\phi$, as $f$ is an $L$-morphism we have $[\phi(\bar{a})]\le [\phi(f``\bar{a})]$, and
  $E\bar{a}\land E\phi = Ef``\bar{a}\land E\phi$.
  \end{proof} 

\begin{theorem}\label{thm_diagram_equivs_weak_L-monomorph} Let $M$ and $N$ be $L$-structures in $\psho$. The following are equivalent.
  \begin{eqn}\begin{split}
      (1)&\hbox{There is a weak }L\hbox{-monomorphism } f:M\longrightarrow N .\\
      (2)&\hbox{There is an expansion of }N\hbox{ to an }L_M\hbox{-structure } \tilde{N}\hbox{ with }
\tilde{N}\forces {\mathcal E}\Delta_M.  \\  
(3)&\hbox{There is an expansion of }N\hbox{ to an }L_M\hbox{-structure }\tilde{N}\hbox{ with }
\tilde{N}\forces \Delta_M.
    \end{split}
  \end{eqn}
  \end{theorem}

\begin{proof} This theorem is similarly immediate from the proof of (1) implies (3) in
  the corresponding result \cite{BM14}, Proposition (4.7), which follows directly from Lemma (\ref{forcing_preservation_lemma}(3)),
  using Lemma (\ref{up_pres_interp_by_conj_disj}) and Lemma (\ref{pres_positive_existential_up}).
\end{proof}

\begin{remark}
These preservation results could be extended to the infinitary logic $L_{\infty,\omega}$
    where arbitrary (set-sized) conjunctions and disjunctions are allowed (and a class sized collection of variables), and
    where the inductive definitions of $[\bigwedge_{i\in I} \phi_i]$ and
    $[\bigvee_{i\in I} \phi_i]$ are simply the infinite versions of the definitions for $\land$ and $\lor$.
\end{remark}

\begin{remark} In \cite{BMPIII}, \S{}\ref{additional_connectives} we discuss extending the languages $L$ by adding
  a unary connectives $\squaresub{p}$ for each $p\in\Omega$ in the style of \cite{BM14},\S{}5.
  At the close of \cite{BMPIII}, \S{}\ref{additional_connectives} we show how one can re-write results on the method of diagrams
  from this section in the form of assertions between the equivalence of the existence of a suitable $L$-$\Gamma$-monomorphism
from $M$ to $N$ and the existence of an extension $\tilde{N}$ of $N$ which forces the relevant $L^{\squ}_M$-sentences
which $M$ forces.
\end{remark}

\section{Categories of interest and directed colimits}\label{cats_and_dir_colims}

\vskip12pt

Now we define the categories of our principal interest.

\begin{definition}\label{categories_of_interest} (Categories of interest) Let $L$ be a first order language
  with equality over $\Omega$ and $T$ an $L$-theory.
Obejcts in our categories of interest are pairs $(A,M)$, taking one of four forms. 
  \begin{eqn}
    \begin{split}
&M\in \psh{\Omega,L},\, M\forces T,\, A\in \psh{\Omega}\hbox{ and }A\hbox{ is a subpresheaf of }M,\hbox{ or}\\
      &M\in \sh{\Omega,L}, \, M\forces T,\, A\in \psh{\Omega}\hbox{ and }A\hbox{ is a subpresheaf of }M,\hbox{ or}\\
      &M\in \sh{\Omega,L},\ssss M\forces T,\ssss A\in \sh{\Omega},\hbox{ and } A\hbox{ is a subsheaf of }M,\hbox{ or}\\
      &M\in \sh{\Omega,L},\ssss M\forces T,\ssss A\in \sh{\Omega},\ssss A\hbox{ is a subsheaf of }M\hbox{ and }d(A)\le d(M).
    \end{split}
  \end{eqn}

We equip each of these collections of objects with the six collections of morphisms given
in Definition (\ref{morphism_notions}).

That is, $(A,M) \xlongrightarrow{f} (B,N)$ is a morphism  if
$f$ is
  \begin{itmz}
  \item a presheaf morphism from $M$ to $N$,
  \item presheaf monomorphism, 
  \item an $L$-morphism, 
  \item a weak $L$-monomorphism
  \item an $L$-monomorphism,
          \item an elementary $L$-monomorphism.
  \end{itmz}
 respectively.

Note, we really do mean to take $A\in \psh{\Omega}$ or $A\in \sh{\Omega}$ rather than $\psh{\Omega,L}$ or $\sh{\Omega,L}$.
We could prove versions of the results in \S{}\S{}5-7 for the latter variants, however we would need stronger notions of morphisms,
which would also preserve logical properties of $A$, for example, to take the strongest case, an elementary $L$-monomorphism whose restriction
to $A$ is an elementary $L$-monomorphism from $A$ to $B$. 
 
 For fixed $\Omega$ and $T$, we give these categories (somewhat ugly) names:
 \begin{eqn}
   \begin{split}
& (\hbox{\normalfont \bfseries SubpSh}(\Omega),\hbox{\normalfont \bfseries pShMod}(T))_x,\\
& (\hbox{\normalfont \bfseries SubpSh}(\Omega),\hbox{\normalfont \bfseries ShMod}(T))_x ,\\
& (\hbox{\normalfont \bfseries SubSh}(\Omega),\hbox{\normalfont \bfseries ShMod}(T))_x\hbox{ and} \\
     & (\hbox{\normalfont \bfseries bdSubSh}(\Omega),\hbox{\normalfont \bfseries ShMod}(T))_x\hbox{, respectively,}
   \end{split}
 \end{eqn}depending on the notion of object, where
the index $x$ indicates the type of morphisms: one of psh, psh-mono,
$L$, wk-$L$-mono, $L$-mono, or elt-$L$-emb.  For example, $\kern-2pt
(\hbox{\normalfont \bfseries bdSubSh}(\Omega),\hbox{\normalfont
  \bfseries ShMod}(T))_{\hbox{\scriptsize elt-$L$-emb}}$ is the
category whose objects are pairs $(A,M)$ such that $M\in
\sh{\Omega,L}$, $M\forces T$, $A\in \sh{\Omega}$ and $A$ is a subsheaf
of $M$ with $d(A)\le d(M)$, and whose morphisms are elementary
$L$-embeddings.

When there is no confusion we will shorten these names by omitting
$\Omega$ and/or $T$, \emph{e.g.}, the same example will be called
$\kern-2pt (\hbox{\normalfont \bfseries bdSubSh},\hbox{\normalfont
  \bfseries ShMod})_{\hbox{\scriptsize elt-$L$-emb}}$.
\end{definition}

\vskip12pt

The next proposition discusses directed colimits in the categories of interest in the case that $T=\emptyset$,
that is that the ``models'' are simply $L$-structures. 
It shows that directed colimits exist under all of the types of morphism we consider.

\begin{proposition}\label{directed_colimits_exist} (\cite{M88}). Let ${\mathfrak C}$ be any
  of the categories defined in Definition
  (\ref{categories_of_interest}) \emph{in the case of $T=\emptyset$}.
  (So the ``models'' are simply $L$-structures in $\psho$.)  Let
  $\tup{\tupof{(A_i, M_i)}{i\in I},\tupof{g_{ij}}{i\pret j}}$ be a
  directed system in ${\mathfrak C}$, where if $i\pret j$ then
  $(A_i,M_i) \xlongrightarrow{g_{ij}} (A_j,M_j)$ is a morphism in the
  category.  Then the colimit of the system exists in
  ${\mathfrak C}$.
\end{proposition}

We will give a proof since Miraglia's proof leaves out many details and
has some unfortunate typos, making it difficult for the reader to flesh out what
was skimmed.
\vskip6pt

\begin{proof}   Throughout we write $[. = .]_i$ for $[.=.]_{M_i}$, where the latter
  is as defined in Definition (\ref{defn_evaluation_for_presheaves}).

  We start by defining the underlying set of the colimit and showing this is a good definition. We then show how
  to equip the set with the structure of a presheaf and show the resulting presheaf is extensional. As a final step
  in the presheaf part of the construction define morphisms from the $(A_i,M_i)$ into the colimit and show
  these are presheaf morphisms or presheaf monomorpisms, respectively if the morphisms
  in directed system have the corresponding property.
  
  Let $|M^*| = \bigcup (|M_i| \times \set{i})$ be the disjoint union of the underlying sets of the $M_i$s.
  Also, let $|A^*| = \bigcup (|A_i| \times \set{i})$.

  We define an equivalence relation on $|M^*|$ by
  $(x,i) \sim (y,j)$ if there is some $l\in I$ with $i$, $j\pret l$ and $g_{il}(x)=g_{jl}(y)$.

  It is immediate that $\sim$ is reflexive and symmetric. The transitivity of $\sim$ follows
  rapidly from directedness. If
  $(x,i) \sim (y,j) \sim (z,k)$ with $l$, $m$ such that $i$, $j\pret l$ and $j$, $k\pret m$,
  while $g_{il}(x)=g_{jl}(y)$ and $g_{jm}(y)=g_{km}(z)$ we can take $n\in I$ with $l$, $m\pret n$.
  We then have $i$, $j$, $k\pret l$ and
  \[
  g_{in}(x) = g_{ln}\cdot g_{il}(x) = g_{ln}\cdot g_{jl}(y) =
     g_{mn}\cdot g_{jm}(y) = g_{mn}\cdot g_{km}(z) = g_{kn}(z),
  \]
  and hence $(x,i)\sim (z,k)$.

  For $(x,i)\in |M^*|$ write $(x,i)_{\sim}$ for its equivalence class under $\sim$.
  We set $|M| = \setof{ (x,i)_{\sim}}{(x,i) \in |M^*|}$, the quotient of $|M^*|$ under $\sim$. 
  We set $|A| = \setof{ (x,i)_{\sim}}{(x,i) \in |M^*| \sss\sss\&\sss\sss x\in |A_i|}$.

  We will define $\on^M : |M| \times\Omega \longrightarrow |M|$ and
  $E^M : |M| \longrightarrow \omega$. For $(x,i)_{\sim} \in |M|$ and $p\in\Omega$ define
  $((x,i)_{\sim} ) \on^M p$ to be $(x\on^{M_i} p,i)_{\sim}$ and
  $E^M ((x,i)_{\sim}) = E^{M_i}x$.  We set $M = (|M|,\on^M,E^M)$.

  We define $\on^A$ and $E^A$ similarly: for $(x,i)_{\sim} \in |A|$ and $p\in\Omega$ define
  $((x,i)_{\sim} ) \on^A p$ to be $(x\on^{A_i} p,i)_{\sim}$ and
  $E^A ((x,i)_{\sim}) = E^{A_i}x$. We set $M = (|A|,\on^A,E^A)$.

  We need to check these are good definitions. So suppose $(x,i)\sim (y,j)$ and let $k\in I$ be such that
  $i$, $j\pret k$ and $g_{ik}(x)=g_{jk}(y)$. Let $p\in\Omega$.

  By the first clause in the definition of presheaf morphisms we have that
  $E^{M_i} x = E^{M_k} g_{ik}(x)= E^{M_k} g_{jk}(y) = E^{M_j} y$. Hence, $E^M$ is well defined.

  By the second clause in the definition of presheaf morphism we have that
  $g_{ik}(x\on p) = g_{ik}(x) \on p = g_{jk}(y) \on p = g_{jk}(y\on p)$. Consequently we have that
  $(x\on p,i)\sim (y\on p,i)$. This shows that $\on^M$ is well defined.

  The same proofs verbatim with $M_i$, $M_j$, $M_k$ and $M$ replaced by $A_i$, $A_j$, $A_k$ and $A$ respectively show that
  $E^A$ and $\on^A$ are also well defined.

  We next check that $M=(|M|,\on^M,E^M)$ is indeed a presheaf over $\Omega$. Let $a = (x,i)_{\sim} \in |M|$
  and $p$, $q\in \Omega$.

  We have 
      $a\on^M (p\wedge q) = (x\on^{M_i} (p\wedge q),i)_{\sim} = ( (x\on^{M_i} p)\on^{M_i} q,i)_{\sim} 
       = {((x\on^{M_i} p),i)_{\sim} \on^M q} = ((x,i)_{\sim} \on^M p) \on^M q = (a \on^M p)\on^M q, $
            as required. 
  
  We also have that $a \on^M E^Ma = ((x,i)_{\sim}) \on^M (E^{M_i} x) =
  (x\on (E^{M_i}x,i)_{\sim} = (x,i)_{\sim}=a$, as required.

  We, finally, have $E^M (a\on^M p) = E^M((x\on^{M_i} p,i)_{\sim}) = E^{M_i}( x\on^{M_i} p)
  = (E^{M_i}x)\wedge p = (E^M ((x,i)_{\sim})) \wedge p = E^Ma \wedge p$, again, as required.
  
  Again, the same proofs verbatim with $M_i$ and $M$ replaced by $A_i$ and $A$ show that $A$ is a presheaf
  and hence is a subpresheaf of $M$.

  We also show that $M$ is extensional. So suppose $a$, $b\in |M|$ with
  $a= (x,i)_{\sim}$ and $b= (y,i)_{\sim}$ (where we may take the second coordinate of the representative of
  each equivalence class to be the same by directedness), and suppose that $D\subseteq \Omega$ is such that
  for all $p\in D$ we have $a\on^M p=b\on^M p$ and $E^Ma=E^Mb=\Vee D$. (The latter is long-hand for
  $E^Ma=E^Mb = [a=b]_M$.) 

  We have $E^{M_i}x = E^{M_i}y  =\Vee D$ and $x\on^{M_i} p =y\on^{M_i} p$. By the extensionality of $M_i$ this
  gives that $x=y$. Hence $a=b$.

  Once more, the same proofs with  $M_i$ and $M$ replaced by $A_i$ and $A$ show that $A$ is extensional.

  \vskip6pt

  Finally, as far as the pure presheaf side of things goes, we define presheaf morphisms from the
  $M_i$ to $M$ and show they are monomorphisms if the maps witnessing the directed system are.

  For each $i\in I$ we define $(A_i,M_i)   \xlongrightarrow{g_{i}} (A,M)$
  by for each $x\in M_i$ setting $g_i(x) = (x,i)_{\sim}$. Note that if $x\in |A_i|$ then $g_i(x)\in |A|$.

  If $x\in M_i$ and $p\in\Omega$ then $E^{M} g_i(x) = E^{M_i}x$, by the definition of $E^M$ and, further,
  $g_i(x \on^{M_i} p) = (x \on^{M_i} p,i)_{\sim} = (x,i)_{\sim}\on^M p$, by the definition of $\on^M$. Hence for
  $i\in I$ we have that $g_i$ is a presheaf morphism.

  Now, suppose $g_{ij}$ is an presheaf monomorphism for all $i$, $j\in I$ with $i\pret j$.
  
  Let $i\in I$. Suppose $x$, $y\in M_i$ and $g_i(x) = g_i(y)$. Thus $(x,i)_{\sim} = (y,i)_{\sim}$.
  However, by the definition of $(x,i)_{\sim}$, the latter implies there is some $k\in I$ with $i\pret k$ and
  $g_{ik}(x)=g_{ik}(y)$. Since $g_{ik}\on M_i$ is injective this gives $x=y$.
  Hence, we've shown $g_i\on M_i$ is also injective and that $g_i$ is a presheaf monomorphism.
  \vskip6pt

  We suppose from here on that for each $i\in I$ we have $M_i\in\psh{\Omega,L}$ and for $i$, $j \in I$ such that $i\pret j$
    we have $g_{ij}$ is an $L$-morphisms. This assumption on the morphisms is necessary in order to equip the whole of $M$ as an
    $L$-structure.

  We will next equip $M$ with structure so that $M\in \psh{\Omega,L}$. We
  have to define the interpretation of constants, relations and
  functions and show these definitions are good. (Note ``$=$'' is also one of the relations we need to
  check is interpreted correctly.)  We will then check the morphisms
  are of the appropriate type if the ones in the directed system are
  of that type.
  
  We have a presheaf of constants $C$ and need to define a presheaf morphism
  $\cdot^M:C\longrightarrow M$ given by $c\mapsto {\mathfrak c}^{M}$ for $c\in |C|$.

  We already have the interpretations in the $M_i$, given by the morphisms $\cdot^{M_i}$. We also
  have for all $i$, $j\in $ with $i\pret j$ and $c\in C$ that $g_{ij}({\mathfrak c}^{M_i})={\mathfrak c}^{M_j}$.

  So we define $\cdot^M$ by choosing $i\in I$ and setting ${\mathfrak c}^{M} = ({\mathfrak c}^{M_i},i)_{\sim}$. The usual argument using
  directedness shows this definition is independent of the choice of $i$. For any $j\in I$ then there is
  $k\in I$ such that $i$, $j\pret k$ and we have $g_{ik}({\mathfrak c}^{M_i})={\mathfrak c}^{M_k} = g_{jk}({\mathfrak c}^{M_j})$, thus
  $({\mathfrak c}^{M_i},i) \sim ({\mathfrak c}^{M_j},j)$.

  Recall from Definition (\ref{defn_presheaf_prod}) that
    \[|M^n| = \setof{\tup{a_0,\dots,a_{n-1}} }{\forall i < n \; a_i \in |M|
      \sss\sss\&\sss\sss \forall i,\ssss j < n\; E^{M}a_i = E^{M}a_j}.\]   

    For every function symbol $\omega$, we define 
  $\omega^M : |M^n| \longrightarrow |M|$, a presheaf morphism, given by $\bar{a} \mapsto \omega^M(\bar{a})$ as follows.
  
  Suppose $\bar{a} = \tup{a_0,\dots,a_{n-1}} =
  \tup{(x_0,i_0)_\sim,\dots,(x_{n-1},i_{n-1})_\sim}$. Exploiting directedness once more, let
  $k\in I$ be such that for all $j<n$ we have $i_j \pret k$ and for $j<n$ write $y_j$ for $g_{i_j k}(x_j)$.
  Then $\bar{a} = \tup{(y_0,k)_\sim,\dots,(y_{n-1},k)_\sim)}$. We set
  \[\omega^M(\bar{a}) = (\omega^{M_k}(y_0,\dots,y_{n-1}),k)_\sim.\]

  We must check this definition is independent of the choice of representatives of
  $a_0$, \dots, $a_{n-1}$ and of the choice of $k$. So, for $j<n$ let $(x'_j,i'_j)$ be such that
  $(x_j,i_j) \sim (x'_j,i'_j)$ and let
  $k'\in I$ be such that for all $j<n$ we have $i'_j \pret k'$;
  for $j<n$ write $y'_j$ for $g_{i_j k'}(x'_j)$.

  Let $l\in I$ be such that $k$, $k'\pret l$. For all $j<n$ we have
  $g_{kl}(y_j)=g_{k'l}(y'_j)$, $=z_j$, say, and,
  since $g_{kl}$, $g_{k'l}$ are $L$-morphisms,
  \[\omega^{M_k}(y_0,\dots,y_{n-1}) = \omega^{M_l}(z_0,\dots,z_{n-1}) = \omega^{M_{k'}}(y'_0,\dots,y'_{n-1}),\]
  as required.

  It is clear from the definition of $g_i$ that if $\bar{a}=\tup{a_0,\dots,a_{n-1}}\in M_i^n$ then
  \[\omega^M(\tup{g_i(a_0),\dots,g_i(a_{n-1})})= g_i(\omega^{M_i}(\bar{a})).\]
  
  Now we define the interpretation of relation symbols (including equality).
  Let $R$ be an $n$-ary relation symbol in $L$ and $\bar{a} = \tup{a_0,\dots,a_{n-1}} \in |M|^n$.
  (Note, here, in contrast to the function symbol case just treated, $\bar{a}\in |M|^n$, and not $|M^n|$.)

  Set
     \begin{equation*}
           \begin{split}
             [R^M(\bar{a})]_M = & \bigvee  \{\ssss [R^{M_k}(g_{i_0k}(x_0),\dots,g_{i_{n-1}k}(x_{n-1}) ) ]_{M_k} \ssss : \\
             & \ssss \;\;\;  (x_0,i_0)\in a_0 , \dots , (x_{n-1},i_{n-1}) \in a_{n-1} \sss\sss\hbox{ and }\sss\sss
               i_0,\dots,i_{n-1}\pret k  \ssss \} .        
           \end{split}
     \end{equation*}\vskip-6pt
     If $(x_0,\dots,x_{n-1})\in |M_i|^n$ it is immediate from the definition of $R^M$ that
     $[ R^{M_i}(x_0,\dots,x_{n-1}) ]_{M_i} \le [ R^M(g_i(x_0),\dots,g_i(x_{n-1})) ]_M$.

     So once we have shown that each $[R(.)]_M$ is a characteristic function
     we will have shown that for $i\in I$ we have that $g_i$ is an $L$-morphism.
     \vskip6pt
     
     We next give an alternative characterization, in the general situation where the
     $g_{ij}$ for $i\pret j$ are simply $L$-morphisms, of the interpretation of
     relation symbols in $M$. This simplifies some of the proofs below.

     \begin{lemma}\label{char_interp_relns}
       Let $R$ be a $n$-ary relation in $L$.
       Let $\bar{a}=\tup{a_0,\dots,a_{n-1}}\in |M|^n$ and suppose
       for all $j<n$ that $(x_j,i_j)\in a_j$. Then
    \[ [R(\bar{a}) ]_{M} = \Vee \setof{ [R(g_{i_0k}(x_0),\dots,g_{i_{n-1}k}(x_{n-1}) )]_{M_k}}{i_0,\dots, i_k \pret k}.\]
      \end{lemma}

       That is, there is no need to take the supremum over all representatives of
  $\bar{a}$ in the definition of $[R(\bar{a})]_M$, we can simply fix some representatives and
       take the appropriate supremum with respect to their upwards images
       as in the right hand side of the equation of this paragraph.

     \begin{proof} Let $\bar{a}$ and the $(x_j,i_j)$ be as in the statement of the lemma.
     Let $P = \setof{ p\in I }{ i_0,\dots i_{n-1}, i'_0,\dots i'_{n-1} \pret l
    \hbox{ and } \forall j<n\;  g_{i_jl}(x_j)=g_{i'_jl}(x'_j)   } $,
     $K = \setof{ k \in I }{ \forall j< n \; i_j \pret k }$ and
     $K' = \setof{k' \in I }{\forall j< n\; i'_j \pret k' }$.

  On the one hand, we clearly have $P \subseteq K \cap K'$.
  Thus,
  \begin{eqn}
    \begin{split}
  \Vee_{p \in P}  [ R&(g_{i_0p}(x_0),\dots,g_{i_{n-1}p}(x_{n-1}) ) ]_{M_p}
 \le \\
 & \Vee_{k \in K} [R(g_{i_0k}(x_0),\dots,g_{i_{n-1}k}(x_{n-1}) ) ]_{M_k}
  \hbox{, }\\
 & \Vee_{k' \in K'} [R(g_{i'_0k'}(x'_0),\dots,g_{i'_{n-1}k'}(x'_{n-1}) ) ]_{M_{k'}}.
    \end{split}
  \end{eqn}

  However, on the other hand, for any $k \in K$ we have some
  $p \in P$ such that
  \begin{eqn}
    \begin{split}   [R&(g_{i_0k}(x_0),\dots,g_{i_{n-1}k}(x_{n-1}) ) ]_{M_k} \le \\
      &[ R(g_{i_0p}(x_0),\dots,g_{i_{n-1}p}(x_{n-1}) ) ]_{M_p}
          \end{split}
  \end{eqn}
  and, likewise, for any  $k' \in K'$ we have some $p'\in P$ such that
 \begin{eqn}
    \begin{split}
       [R&(g_{i'_0k'}(x'_0),\dots,g_{i'_{n-1}k'}(x'_{n-1}) ) ]_{M_k} \le \\
      & [ R(g_{i'_0p'}(x'_0),\dots,g_{i'_{n-1}p'}(x'_{n-1}) ) ]_{M_{p'}} . 
    \end{split}
  \end{eqn}

  Hence,    \begin{eqn}
    \begin{split}
 \Vee_{k \in K} & [R(g_{i_0k}(x_0),\dots,g_{i_{n-1}k}(x_{n-1}) ) ]_{M_k} \hbox{, }\\
 & \Vee_{k' \in K'} [R(g_{i'_0k'}(x'_0),\dots,g_{i'_{n-1}k'}(x'_{n-1}) ) ]_{M_{k'}} \le \\
 &  \qquad \Vee_{p \in P}  [ R(g_{i_0p}(x_0),\dots,g_{i_{n-1}p}(x_{n-1}) ) ]_{M_p} .
    \end{split}
  \end{eqn}

  Hence, all three disjunctions are equal.

  Consequently, we have shown that
  \[ [R(\bar{a}) ]_{M} = \Vee \setof{ [R(g_{i_0k}(x_0),\dots,g_{i_{n-1}k}(x_{n-1}) )]_{M_k}}{i_0,\dots, i_k \pret k} .\]

  \end{proof}
  
  \vskip6pt

  We next show that each $[R(.)]_M$ is a characteristic function. So,
  suppose $R$ is an $n$-ary relation symbol in $L$. For $\bar{a}\in |M|^n$ let
  \begin{eqn}
    \begin{split}
      X = \{\ssss (\tup{(x_0,i_0)_\sim,\dots,(x_{n-1},i_{n-1})_\sim},k) \ssss : \ssss\forall & j < n \; 
       ( x_j\in |M_{i_j}|\sss\sss\&\sss\sss \\
  &  (x_j,i_j)\in a_j \sss\sss\&\sss\sss i_0\pret k ) \} 
    \end{split}
  \end{eqn}

  Let $\tilde{x}=(\tup{(x_0,i_0)_\sim,\dots,(x_{n-1},i_{n-1})_\sim},k) \in X$. For $j<n$ set $z_j=g_{i_jk}(x_j)$.
  Set $\bar{z}=\bar{z}_{\tilde{x}} = \tup{z_0,\dots,z_{n-1}}$.
  Since $[R(.)]_{M_k}$ is a characteristic function, we have $[R(\bar{z})]_{M_k} \le E^{M_k}\bar{z}$.
  Further, since $g_k$ is a presheaf morphism we have $E^M\bar{a} = E^{M_k}\bar{z}$.

  Thus, $[R(\bar{a})]_M  = \Vee \setof{[R(\bar{z}_{\tilde{x}})]_{M_k}}{\tilde{x} \in X} \le E^Ma$.

  Now suppose $\bar{a}$, $\bar{b}\in |M|^n$. By directedness choose $k\in I$ and $\bar{x}$, $\bar{y}\in |M_k|^n$
  such that for all $j<n$ we have $g_k(x_j)=a_j$ and $g_k(y_j)=b_j$.

  By Lemma (\ref{char_interp_relns}), for the first equality, the
  distributivity of $\wedge$ over $\Vee$ twice, for the second
  equality, the fact that the $[R(.)]_{M_l}$ are characteristic functions for all $l\in I$ for the inequality and
  Lemma (\ref{char_interp_relns}), again, for the final equality, we have
  \begin{eqn}
    \begin{split}
  [\bar{a}=\bar{b}]_M  & \wedge [R(\bar{b})]_M  \\
    & = \Vee \setof{ [g_{kl}``\bar{x} = g_{kl}``\bar{y}]_{M_l}}{k \pret l} \wedge
    \Vee \setof{ [R(g_{kl}``\bar{y})]_{M_l}}{k \pret l}  \\
    & = \Vee \setof{ [g_{kl}``\bar{x} = g_{kl}``\bar{y}]_{M_l} \wedge [R(g_{kl}``\bar{y})]_{M_l}}{k \pret l} \\
    & \le \Vee\setof{[R(g_{kl}``\bar{x})]_{M_l}}{k \pret l} =  [R(\bar{a})]_M .
          \end{split}
  \end{eqn}

  Thus we have shown that $[R(.)]_M$ is a characteristic function.
  
  \vskip6pt

  We now need to address the cases where the $g_{ij}$ have further preservation properties.

  First of all, let us suppose that for all $i$, $k\in I$ with $i\pret k$ we have that
  $g_{ij}$ is a weak $L$-monomorpism. We shall show for each $i\in I$ that
  $g_i$ is also a weak $L$-monomorphism.

  Suppose $\bar{x}\in |M_i|^n$ and $\bar{a}=(x_0,i),\dots,x_{n-1},i)_\sim$.
  
  For each $k$ such that $i\pret k$ we have
  \[ [R(g_{ik}``\bar{x})]_{M_k} \le \neg\neg [R(g_{ik}``\bar{x})]_{M_k} = \neg\neg [R(\bar{x})]_{M_i} .\]
  Thus, by Lemma (\ref{char_interp_relns}) above, we have
  \begin{eqn}
    \begin{split}
  [R(\bar{a})]_{M} = \Vee & \setof{[R(g_{ik}``\bar{x})]_{M_k} }{i\pret k}  \le \\
  & \Vee \setof{\neg\neg [R(g_{ik}``\bar{x})]_{M_k} }{i\pret k} = \neg\neg [R(\bar{x})]_{M_i}       
    \end{split}
  \end{eqn}

  However, since, as just shown above
  $g_i$ is an $L$-morphism, we also have
  $[R(\bar{x})]_{M_i} \le [R(\bar{a})]_M$, and hence
  $\neg\neg [R(\bar{x})]_{M_i} \le \neg \neg [R(\bar{a})]_M$.

  Thus we have shown $\neg\neg [R(\bar{x})]_{M_i} = \neg \neg [R(\bar{a})]_M$, as required.
  \vskip6pt

  Now, suppose let us suppose that for all $i$, $k\in I$ with $i\pret k$ we have that
  $g_{ij}$ is an $L$-monomorpism. We shall show for each $i\in I$ that
  $g_i$ is also an $L$-monomorphism.

  Let $i\in I$ and $\bar{x}\in |M_i|^n$ and $\bar{a}=g_i``\bar{x}$.
  By Lemma (\ref{char_interp_relns}) we have that
  $[R(g_i``\bar{x})]_M = \Vee \setof{ [R(g_{i_0k}(x_0),\dots,g_{i_{n-1}k}(x_{n-1}) )]_{M_k}}{i_0,\dots, i_k \pret k}$.

  However, by the hypothesis that for each $k\in I$ with $i\pret k$ we have $g_{ik}$ is an restricted $L$-monomorphism,
  we have for each such $k$ that $[R(\bar{x})]_{M_i} =  [R(g_{ik} ``\bar{x})]_{M_{k}}$.

  Hence $[R(g_i``\bar{x})]_M = [R(\bar{x})]_{M_i}$.  

  We remark that, by Proposition (\ref{equiv_presheaf_mono}), for equality,
  ``$=$'', we simply need the $g_{ik}$ to be presheaf monomorphisms.

  \vskip6pt

  Penultimately, suppose for all $i$, $k\in I$ with $i\pret k$ we have $g_{ij}$ is an elementary $L$-embedding.
  We will show that for all $i\in I$ we have $g_i$ is an elementary $L$-embedding. We prove this by
  induction on the complexity of formulae. 

  By Lemma (\ref{mono_conserves_qf}) we only need to 
  consider the quantifiers.

  Suppose $\phi(t,v_0,\dots,v_n )$ is an $L$-formula with free variables shown.
  Let $i\in I$ and $\bar{x}\in |M_i|^n$, and for $j<n$ let $a_j = (x_j,i)_{\sim}$.
  Then
  \begin{eqn}
    \begin{split}
      [ \exists t \;\phi(t,\bar{a})]_M & = \Vee \setof{[\phi(b,\bar{a})]_M }{b\in |M|}\\
                & = \Vee \setof{ [\phi(y,g_{ik}``\bar{x})]_{M_k} }{i\pret k  \ssss\ssss\&\ssss\ssss y\in |M_k|}\\
      & =\Vee \setof{[\exists t \; \phi(t,g_{ik}``\bar{x})]_{M_k}}{i\pret k} \\
        & =  [\exists t \;\phi(t,\bar{x})]_{M_i}.
    \end{split}
  \end{eqn}Here the first and third equality are given by the definition of interpretation of existential quantification in
  $M$ and the $M_k$ respectively, while the second equality is given by the inductive hypothesis and the final
  equality is given by the elementarity of the $g_{ik}$.

  Sentences with a leading universal quantifier are handled similarly.

  \begin{eqn}
    \begin{split}
      [ \forall t \;\phi(t,\bar{a})]_M
          & =  E\bar{a} \land E\phi \land \bigwedge \setof{Eb \longrightarrow [\phi(b,\bar{a})]_M }{b\in |M|}\\
          & = E\bar{x} \land E\phi \land  \bigwedge
               \setof{Ey \longrightarrow [\phi(y,g_{ik}``\bar{x})]_M }{i\pret k \ssss\ssss\& \ssss\ssss y\in |M_k|}\\
          & = \bigwedge \setof{  [ \forall t \;\phi(t,g_{ik}``\bar{x})]_{M_k}}{i\pret k} \\
          & =  [\forall t \;\phi(t,\bar{x})]_{M_i},
    \end{split}
  \end{eqn}

  \vskip6pt
  
  We now address sheafification. 

  By \cite{BM14}, Theorem 2.15, which summarizes \cite{Mir06}, Theorem 27.9, Theorem 37.8 and Corollary 37.9,
  if we apply the sheafification process to both $A$ and $M$ (separately), obtaining $cA$ and $cM$
  there is a unique morphism $d$ from $cA$ to $cM$ such that we get a commutative square  

\[  
\begin{tikzcd}
A \arrow{r}{c^A} \arrow{d}{id} & cA \arrow{d}{d} \\
M \arrow{r}{c^M} & cM
\end{tikzcd}
\]

By the construction (see \cite{Mir06}, Example 27.12) we have that $g$
is a presheaf monomorphism.  Specifically, for $a\in |A|$ we have that
$c^A(a) \in ^{|A|}\Omega$ is the function taking each $y\in |A|$ to
$[a=y]_A$ and $c^M(a) \in ^{|M|}\Omega$ is the function taking each
$y\in |M|$ to $[a=y]_M$, and, moreover, for all such $a$ and $y$ we
have $[a=y]_A = [a=y]_M$. Thus $c^A(a) = c^M(a) \on |A|$, where here
`$\on$' is used in the sense of the restriction of a function between
two sets to a subset of its domain. Thus for all $a\in |A|$ we have
$d(c^A(a)) = c^M(a)$.

It is immediate from what we have already shown  that for all $a\in |M|$ that we have
$Ea = \bigvee_{i\in I} \bigvee_{b\in |M|_i} [g_i(b) = a]_M$. 

We thus have that $(d``cA,cM)$ is the directed colimit.
  
\end{proof}

We now consider preservation phenomena of a different type to those of the previous section.
Rather than concerning individual morphisms, 
these phenomenona concern situations where $T$ is a theory, and there is a directed colimit for which
all of the objects in the diagram whose colimit we are taking
force $T$ and all of the morphisms in the diagram satisfy certain properties. The preservation property
is that the colimit also forces $T$.

Proposition (\ref{directed_colimits_exist}) allows us to prove the following theorem and corollary which
should be considered as showing various instances such phenomena.

\begin{theorem}\label{gen_of_tarski's_thm} Let $\tup{\tupof{M_i}{i\in I},\tupof{g_{ij}}{i\pret j}}$ be a
      directed system in the category whose objects are $L$-structures in $\psho$ and whose morphisms
      are $L$-monomorphisms. Let $\tup{M,\tupof{g_{i}}{i\in I}}$ be the directed colimit of the
      system as given by Proposition (\ref{directed_colimits_exist}).   
      Let $\bar{a} \in |M|^m$. For each $i\in I$ with $\bar{a} \subseteq \rge(g_i)$ let
      $\bar{a}_i$ be such that $\bar{a} = g_i``\bar{a}_i$.

      Let $ \phi(\bar{x},\bar{y})$ be a ${\mathcal E}_{1}$-formula with free variables shown, where
      $\lh(\bar{x})=n$ and $\lh(\bar{y})=m$.
      If for all $i\in I$ with $\bar{a}\subseteq \rge(g_i)$ we have
      $M_i\forces \forall \bar{x} \sss\sss \phi(\bar{x},\bar{a}_i)$, then we also have
      $M\forces \forall \bar{x} \sss\sss \phi(\bar{x},\bar{a})$. 
\end{theorem}

\begin{proof} By Proposition (\ref{characterization_of_forcing_forall_E1}) it suffices to show
  that for all $\bar{b}\in |M|^n$ we have $E\phi \wedge E\bar{a} \wedge E\bar{b} \sss\sss\le
  [\phi(\bar{t},\bar{a})]_M$. So let $\bar{b}\in |M|^n$ and let $i\in I$ be such that
  $\bar{b}$, $\bar{a} \subseteq \rge(g_i)$. Let $\bar{b}_i$ be such that $g_i``\bar{b}_i = \bar{b}$.

  By the hypothesis on $M_i$
  we have $E\phi \wedge E\bar{a}_i \wedge E\bar{b}_i \sss\sss\le  [\phi(\bar{b}_i,\bar{a}_i)]_{M_i}$
  (again using Proposition (\ref{characterization_of_forcing_forall_E1})).

  As $g_i$ is a presheaf morphism $E\bar{a} = E\bar{a}_i$ and $E\bar{b}=E\bar{b}_i$.
  As $g_i$ is an \mbox{$L$-monomorphism}, by Lemma (\ref{pres_existential_up_mono}) we have
      $[\phi(\bar{b}_i,\bar{a}_i)]_{M_i} \le [\phi(\bar{b},\bar{a})]_M$.
  Thus
  $ E\phi \wedge E\bar{a} \wedge E\bar{b} = 
  E\phi \wedge E\bar{a}_i \wedge E\bar{b}_i \sss\sss\le  [\phi(\bar{b},\bar{a}_i)]_{M_i}
  \le [\phi(\bar{b},\bar{a})]_M$, as required.
\end{proof}

\vskip12pt

\begin{corollary}\label{directed_colimits__for_models_of_theories_exist}
  Let ${\mathfrak C}$ and $T$ be any of the following pairs:
  \begin{itmz}
  \item \sloppy{${\mathfrak C}$ is any
    of the four categories whose morphisms are elementary
    \mbox{$L$-embeddings} defined in Definition (\ref{categories_of_interest})
  and $T$ is any set of \mbox{$L$-sentences}}
  \item \sloppy{${\mathfrak C}$ is any
  of the four categories whose morphisms are \mbox{$L$-monomorphisms} defined in Definition (\ref{categories_of_interest})
  and $T$ is any set of $\forall{\mathcal E}_{1}$-sentences}
  \item\sloppy{${\mathfrak C}$ is any
  of the categories defined in Definition (\ref{categories_of_interest})
  and $T =\emptyset$}
\end{itmz}
  Let $\tup{\tupof{(A_i, M_i)}{i\in I},\tupof{g_{ij}}{i\pret j}}$ be a
  directed system in ${\mathfrak C}$, where if $i\pret j$ then $(A_i,M_i)
  \xlongrightarrow{g_{ij}} (A_j,M_j)$ is a morphism in the category.
  Then the colimit of the directed system exists
  in ${\mathfrak C}$.
\end{corollary}

\begin{proof} Immediate from Proposition (\ref{directed_colimits_exist}),
  the definition of elementary $L$-embedding
and Theorem (\ref{gen_of_tarski's_thm}).
  \end{proof}

\section{Presheaves and sheaves of models and accessibility}\label{accessibility}

We start by introducing the category theoretic notions of presentability and accessiblity.
The reader can consult, for example, \cite{Adamek-Rosicky} and \cite{Makkai-Pare}
for further information on this area of category theory.

\begin{definition} 
  Let $\lambda$ be a regular cardinal and ${\mathcal C}$ a category.

  A ${\mathcal C}$-object $M$ is $\lambda$-\emph{presentable} if given $\tup{N,\tupof{\phi_i}{i\in I}}$, the
  colimit in ${\mathcal C}$ of a $\lambda$-directed system
  $\tup{\tupof{N_i}{i\in I},\tupof{\phi_{ij}:N_i\longrightarrow N_j}{i\le_I j,\sss\sss i,\sss j\in I}}$,
  and  $f : M \longrightarrow N$, a ${\mathcal C}$-morphism, 
  there is some $i\in I$ and ${\mathcal C}$-morphism $g:M\longrightarrow N_i$ such
  that $f = \phi_i\cdot g$.

  The \emph{presentability rank} of an object $M$ of ${\mathcal C}$ is the smallest regular cardinal
  $\kappa$ such that $M$ is $\kappa$-presentable.

  The \emph{inaccessible presentability rank} of an object $M$ of ${\mathcal C}$ is the
  smallest strongly inaccessible cardinal
  $\kappa$ such that $M$ is $\kappa$-presentable.
  \end{definition}

\begin{definition}
  A category ${\mathcal C}$ is $\lambda$-\emph{accessible}
  if it is closed under $\lambda$-directed colimits
  (\emph{i.e.}, colimits indexed by a $\lambda$-directed partially ordered set)
  and contains, up to isomorphism, a set $R$ of $\lambda$-presentable objects
  such that each object of $\mathcal C$ is a $\lambda$-directed colimit of objects from $R$.
\end{definition}
\vskip12pt

\begin{definition}
  A category ${\mathcal C}$ is \emph{accessible} if there is some $\lambda$ for which it is $\lambda$-accessible.
\end{definition}

We now show how certain of the categories discussed in the previous
section are accessible.  Note for our second theorem that we require
the mild set theoretic hypothesis that there is an strongly
inacessible cardinal larger than the sizes of $L$ and
$\Omega$. Clearly to prove the theorem for all possible $L$ or all
possible $\Omega$ we would need a proper class of strongly
inaccessible cardinals.

\begin{theorem}\label{accessibility_thm_1}  Let $\lambda>\card{L}$, $\chi(\Omega)$ be a regular cardinal.
  Let ${\mathfrak C}$ be any of the categories
  \mbox{$(\hbox{\normalfont \bfseries bdSubShMod}(T),\hbox{\normalfont \bfseries ShMod}(T))_x$}, where
  either $x=L$-mono and the theory $T$ is a set of $\forall{\mathcal E}_{1}$ sentences, or
  $x=$elt-$L$-emb (and $T$ is arbitrary).
       Then ${\mathfrak C}$ is
      $\lambda$-accessible.

  Moreover, if $\kappa$ is a regular cardinal then
  $(A,M)\in {\mathfrak C}$ has presentability rank $\kappa^+$ if and only if $d(M)=\kappa$.
  \end{theorem}
\vskip6pt

\begin{theorem}\label{accessibility_thm_2}   Let $\lambda>\card{L}$, $\card{\Omega}$ be an strongly inaccessible cardinal.
  Let ${\mathfrak C}$ be any of the categories of the form
  $(\hbox{\normalfont \bfseries SubpShMod}(T),\hbox{\normalfont \bfseries ShMod}(T))_x$,
  $(\hbox{\normalfont \bfseries SubpShMod}(T),\hbox{\normalfont \bfseries pShMod}(T))_y$ or
  \mbox{$(\hbox{\normalfont \bfseries SubShMod}(T),\hbox{\normalfont \bfseries ShMod}(T))_z$}, for each $v\in \set{x,y,z}$
  where either $v=L$-mono and the theory $T$ is a set of $\forall{\mathcal E}_{1}$ sentences, or $v=$elt-$L$-emb (and $T$ is arbitrary).
  Then ${\mathfrak C}$ is $\lambda$-accessible.
  \end{theorem}

We prove the two theorems together, with the argument bifurcating in a couple of places in ways corresponding to the
differing hypotheses. The idea of the proof to use the collection of isomorphism types of models of density
less than $\lambda$ as the presentable objects. This is a candidate collection of representatives
as we have see in Proposition (\ref{bound_size_in_terms_of_density})
that there are only a set many of these isomorphism types, and not a proper class
as one might have imagined \emph{a priori}. It remains to show
that every element of ${\mathfrak C}$ can be represented as a $\lambda$-directed colimit
of such models and that every such model is $\lambda$-presentable; it is only in showing the
latter that the proof of the statements of the two theorem bifurcates.

\vskip12pt

In the course of the proof we  will need (the general version of)
Miraglia's downward L\"owenheim-Skolem theorem for sheaves of structures, and thus
we introduce this here.

\begin{definition} For any complete Heyting algebra $\Omega$ the \emph{cardinal character} of $\Omega$,
  $\chi(\Omega)$, is the least cardinal $\kappa$ such that for all $X\subseteq \Omega$ there
  are $Y$, $Z\in [\Omega]^{\le\kappa}$ such that $\bigvee X=\bigvee Y$ and $\bigwedge X = \bigwedge Z$.

Note that the cardinal character of a complete Heyting algebra
is always defined since, at worst, we always have
$\chi(\Omega) \le \card{\Omega}$.\end{definition}

\begin{definition} Let $\mu(\Omega,L) =  \max(\set{\chi(\Omega),\card{L}})$
\end{definition}

\begin{proposition}\label{downward-LS} (Downward L\"owenheim-Skolem theorem, \cite{M88}) Let $\Omega$ be a
  complete Heyting algebra and $L$ a first order language with
  equality.  Let $M\in \psh{\Omega,L}$ (resp.~$\sh{\Omega,L}$). For
  all $\kappa\ge \mu(\Omega,L)$ and all $X\in [|M|]^\kappa$ there is
  an elementary substructure $N$ of $M$ such that $X\subseteq |N|$ and
  $d(N)\le \kappa$.

  Moreover, if $A$ is a subpresheaf of $M$ then ${A\cap N}$ is a
  subpresheaf of $N$; if $M$ is a sheaf and $A$ is a subsheaf of $M$ then ${A\cap N}_s$ is a subsheaf of $N$; and 
  if $M$ is a sheaf and $Y\subseteq X$ then $\gen{Y}_s$ is a subsheaf of $N$.
  \end{proposition}

\begin{proof} See \cite{M88}, Theorem 3.1 and Corollary 3.1.
  The proof given there is for complete Heyting algebras with
  countable character, countable languages and countable subsets of the model, but the alterations
  required for the general case (already envisaged in \cite{M88} and stated here) are trivial.

  The clauses in the second paragraph of the statement are all immediate from Lemma (\ref{intersection_of_presh_lemma}) and Lemma (\ref{uniqueness_of_gen_by}).
  \end{proof}

We now start the proofs of Theorems (\ref{accessibility_thm_1}) and (\ref{accessibility_thm_2}).

\begin{proof}

  \begin{lemma}\label{lemma_1_of_access_thm_1}
    Let ${\mathfrak C}$ be a category as in Theorem (\ref{accessibility_thm_1}) and
    $(A,M) \in \ob({\mathfrak C})$. Then $(A,M)$ is
    a $\lambda$-directed colimit of the set of its elementary substructures $(A',M')$
    with $d(A') \le d(M') < \lambda$.
\end{lemma}

  \begin{proof}
  Let ${\mathcal M}$ be the collection of all substructures $(A',M')$ 
  of $(A,M)$ with $d(A')\le d(M') < \lambda$ and such that the inclusion map is a morphism.
  Let ${\mathcal N} \in
  [\mathcal M]^{<\lambda}$ and for $(B,N)\in {\mathcal N}$ let $D_B$ be dense in $B$ and
  let $D_N$ be a superset of $D_B$ dense in $N$  and such that $\card{D_N} <\lambda$.

  Let $S= \bigcup_{(B,N)\in {\mathcal N}} D_N$ and let $T=\bigcup_{(B,N)\in {\mathcal N}} D_B$.
  Note, here  $\card{S}<\lambda$, as it is the union of
  fewer than $\lambda$ many sets of size less than $\lambda$, and
  $d(S) \le \card{S}$.

  By the downward L\"owenheim-Skolem theorem
  (Proposition \ref{downward-LS}) there is some elementary
  substructure $P$ of $M$ with $d(P)<\lambda$ and $S\subseteq |P|$.
  Thus $P\in {\mathcal M}$.

  Further, by the second paragraph of the statement of the theorem,
  $C=\gen{T}_s$ is a subsheaf of $P$
  and $d(C) \le d(P)$. 
  
  Since for every $(B,N)\in {\mathcal N}$ we have $D_N \subseteq S\subseteq P$,
  by Lemma (\ref{reconstruction_from_dense_in_sheaves}) we have that $N\subseteq P$.
  
  Hence, by Proposition
  (\ref{coherence}), we have for every $(B,N)\in {\mathcal N}$ that $(B,N)$
  is a substructure of $(C,P)$ such that the inclusion map is a morphism.
  
  Thus we have shown that ${\mathcal M}$ is $\lambda$-directed.
  \end{proof}
\vskip6pt
  
  \begin{lemma}\label{lemma_1_of_access_thm_2} Let ${\mathfrak C}$ be a category as in Theorem (\ref{accessibility_thm_2}) and
    $(A,M) \in \ob({\mathfrak C})$. Then $(A,M)$ is
    a $\lambda$-directed colimit of the collection of its elementary substructures $(A',M')$
    with $d(M') < \lambda$.
\end{lemma}

  \begin{proof}
  Let ${\mathcal M}$ be the collection of all substructures
  of $(A,M)$ with density less than $\lambda$, equivalently of size less than $\lambda$ and such that the inclusion map is a morphism.
  
  By Proposition (\ref{bound_size_in_terms_of_density}) and elementary cardinal
  arithmetic $\card{\bigcup {\mathcal N}} < \lambda$.
  Set $S= \bigcup {\mathcal N}$ and note that
  $d(S) \le \card{S}$.

  By the downward L\"owenheim-Skolem theorem
  (Proposition \ref{downward-LS}) there is some elementary
  substructure $P$ of $M$ with $d(P)<\lambda$ and $S\subseteq |P|$.
  Thus $P\in {\mathcal M}$.

  It is immediate from the definition of $S$ that
  for every $(B,N)\in {\mathcal N}$ we have that $N\subseteq P$.
  
  Hence, by Proposition
  (\ref{coherence}), we have for every $N\in {\mathcal N}$ that $N$
  is an substructure of $P$ such that the inclusion map is a morphism.

  Now set $C = \tup{\bigcup \setof{B\subseteq M}{(B,N)\in {\mathcal N}}}$. Then $C$ is a subpresheaf of $P$ and 
  clearly for $(B,N)\in {\mathcal N}$ we have $B\subseteq C $. Consequently $(C,P)$ witnesses
  $\lambda$-directed in the instance of $\mathcal N$.

  Thus we have shown that ${\mathcal M}$ is $\lambda$-directed.
  \end{proof}

  The next lemma is where in the proof of Theorem (\ref{accessibility_thm_1})
  we have need of the constraint on ${\mathfrak C}$-objects which are pairs consisting of a sheaf and a subsheaf
  that the subsheaves are of density no larger than the that of the larger sheaves.

  Here we need also that the maps in the diagram, and hence the morphisms in the category, are in fact
  of the form required for Proposition (\ref{coherence}).
  This is what restricts our results to categories whose morphisms are 
  $L$-monomorphisms or elementary $L$-embeddings.

  \begin{lemma}\label{lemma_2_of_access_thm_1} Let ${\mathfrak C}$ be a category as in
    Theorem (\ref{accessibility_thm_1}). Every $(A,M) \in \ob({\mathfrak C})$ with
    $d(A) \le d(M)<\lambda$ is $\lambda$-presentable.
\end{lemma}

\begin{proof} Suppose $d(A)\le d(M)<\lambda$,
  $(B,N)\in {\mathfrak C}$ is a colimit of the $\lambda$-directed
  system $\tupof{\phi_{ij}:(B_i,N_i)\longrightarrow (B_j,N_j)}{i\le_I j,\sss\sss i,\sss j\in I}$ and
  $f:(A,M)\longrightarrow (B,N)$ is an $x$-morphism.
  By Proposition (\ref{directed_colimits_exist}), $|N| = \bigcup_{i\in I} \phi_i ``N_i$ and
  $|B|=\bigcup_{i\in I} \phi_i ``B_i$.

  { \begin{sublemma}  There is some $i\in I$ such that $f``M \subseteq \phi``N_i$.
  \end{sublemma}

  \begin{proof} As $f$ is a sheaf morphism, we have $d(f ``M) = d(M)<\lambda$.
  Let $D\in [f ``M]^{<\lambda}$ be dense in $f `` M$.
  By the $\lambda$-directedness of the system there is some $i\in I$ such that
  $D\subseteq \phi_i ``N_i$. By Corollary (\ref{reconstruction_from_dense_in_sheaves})
  (with $f ``M$ for $A$, $D$ for $D$, $\phi_i `` N_i$ for $N$ and $N$ for $P$), we have
  $f `` M \subseteq \phi_i `` N_i$.
  \end{proof} }
  
  Fix $N_i$ as in the sublemma. By Proposition (\ref{coherence}),
  we have that $f``M$ is an submodel of $\phi_i``N_i$ and the inclusion is a morphism.

  Thus
  $\phi_i^{-1}\cdot f$ is the desired factorization of $f$ through $\phi_i$.
  Moreover, the factorization is unique: if
  $g : M\longrightarrow N_i$ with $\phi_i \cdot g = f = \phi_i \cdot (\phi_i^{-1}\cdot f)$
  then since $\phi_i$ is a monomorphism we have $g = \phi_i^{-1}\cdot f$.
\end{proof}

\begin{lemma}\label{lemma_2_of_access_thm_2} Let $\lambda$ be strongly inaccessible and
  ${\mathfrak C}$ a category as in
    Theorem (\ref{accessibility_thm_2}). Every $(A,M) \in \ob({\mathfrak C})$ with
    $d(M)<\lambda$ is $\lambda$-presentable.
\end{lemma}

\begin{proof} Suppose $d(A)\le d(M)<\lambda$,
  $(B,N)\in {\mathfrak C}$ is a colimit of the $\lambda$-directed
  system $\tupof{\phi_{ij}:(B_i,N_i)\longrightarrow (B_j,N_j)}{i\le_I j,\sss\sss i,\sss j\in I}$ and
  $f:(A,M)\longrightarrow (B,N)$ is an $x$-morphism.
  By Proposition (\ref{directed_colimits_exist}), $|N| = \bigcup_{i\in I} \phi_i ``N_i$ and
  $|B|=\bigcup_{i\in I} \phi_i ``B_i$.

  { \begin{sublemma}  There is some $i\in I$ such that $f``M \subseteq \phi``N_i$.
  \end{sublemma}

  \begin{proof} We have
  $\card{f``M} = \card{M} \le \card{\Omega}^{d(M)} <\lambda$.
  Thus, by the $\lambda$-directedness of the system, there is some $i\in I$ such that
  $f``M \subseteq \phi``N_i$.\end{proof} }

    Fix $N_i$ as in the sublemma. By Proposition (\ref{coherence}),
  we have that $f``M$ is an submodel of $\phi_i``N_i$ and the inclusion is a morphism.

  Thus
  $\phi_i^{-1}\cdot f$ is the desired factorization of $f$ through $\phi_i$.
  Moreover, the factorization is unique: if
  $g : M\longrightarrow N_i$ with $\phi_i \cdot g = f = \phi_i \cdot (\phi_i^{-1}\cdot f)$
  then since $\phi_i$ is a monomorphism we have $g = \phi_i^{-1}\cdot f$.
  \end{proof}

  This completes the proof of the main two parts of the theorems.
\vskip6pt
  
  \begin{lemma}\label{lemma_3_of_access_thm} Let ${\mathfrak C}$ be a category as in
    Theorem (\ref{accessibility_thm_1}) or
    Theorem (\ref{accessibility_thm_2}). Every $(A,M) \in \ob({\mathfrak C})$ which is
    $\lambda$-presentable has $d(M)<\lambda$. Moroever, if ${\mathfrak C}$ a category as in
    Theorem (\ref{accessibility_thm_1}) then $d(A)\le d(M)$.
  \end{lemma}

  \begin{proof}   By Lemma (\ref{lemma_1_of_access_thm_1}) and Lemma (\ref{lemma_1_of_access_thm_2}), respectively,
    any $(A,M)\in {\mathcal K}$
  is the $\lambda$-directed colimit of the collection, ${\mathcal M}$, of its substructures $(B,N)$ with
  the property that $d(B)\le d(N)<\lambda$ in the former case and $d(N)<\lambda$ in the latter.

  Consider $id:M\longrightarrow M$. If $(A,M)$ is $\lambda$-presentable
  this map factors through the identity inclusion of some $N\in {\mathcal M}$ into $M$.
  However, the inclusions are monomorphisms, so we have $M=N$ and thus $d(M)<\lambda$.
\end{proof}

  The final paragraph of the theorem is now immediate from Lemma (\ref{lemma_2_of_access_thm_1}), or
  Lemma (\ref{lemma_2_of_access_thm_2}), and Lemma (\ref{lemma_3_of_access_thm}).
\end{proof}

\begin{remark} Let ${\mathfrak C}$ be a category and $\lambda$ a strongly inaccessible cardinal as in
  as Theorem (\ref{accessibility_thm_2}).
  Then Lemma (\ref{lemma_2_of_access_thm_2}) and Lemma (\ref{lemma_3_of_access_thm}) show that for
  any $(A,M)\in {\mathrm ob({\mathfrak C})}$ we have $d(M)<\lambda$ if and only if $(A,M)$ is $\lambda$-presentable.
  Hence inaccessible presentability rank, a category theoretic notion,  matches up with density,
  a presheaf theoretic notion, albeit at a less finely-grained scale than that at which
  presentability rank matches up with sheaf theoretic density.
  \end{remark}

\section{Presheaves and sheaves of models and independence relations}\label{indep_relns}

Kamsma (\cite{K22}) gives a framework in which he develops a theory of independ\discretionary{-}{}{}ence
relations extending that
applicable to NSOP$_1$ theories in first order model theory.

\begin{definition} (\cite{K22}, Definition 2.2) An \emph{AECat} consists of a pair $({\mathcal C},{\mathcal M})$ where
  $\mathcal C$ and $\mathcal M$ are accessible categories, $\mathcal M$ is a full subcategory of $C$,
  $\mathcal M$ has directed colimits which the inclusion functor into $\mathcal C$ preserves, and
  all morphims in $\mathcal C$ (and hence in $\mathcal M$) are monomorphism.
\end{definition}

\begin{definition}   
  We say $(\mathcal C, \mathcal M)$ has the \emph{amalgamation property} (\emph{AP}) if $\mathcal M$ has
  the amalgamation property.
\end{definition}

We have already seen
in Corollary (\ref{directed_colimits__for_models_of_theories_exist}) that the categories for which we proved
Theorem (\ref{accessibility_thm_1}) and Theorem (\ref{accessibility_thm_2}) have directed colimits.
In these categories the morphisms  are all monomorphisms. Hence, the following propositions are immediate.

\begin{proposition}\label{aecat_coroll_1}
Let $\lambda>\card{L}$, $\chi(\Omega)$ be a regular cardinal.
  Let ${\mathfrak C}$ be any of the categories
  \mbox{$(\hbox{\normalfont \bfseries bdSubShMod}(T),\hbox{\normalfont \bfseries ShMod}(T))_x$}, where
  either $x=L$-mono and the theory $T$ is a set of $\forall{\mathcal E}_{1}$ sentences, or
  $x=$elt-$L$-emb (and $T$ is arbitrary).
   Let ${\mathfrak M}$ be in each case the full subcategory of ${\mathfrak C}$ whose class of
   objects consists of all pairs of the form $(M,M)$. Then $({\mathfrak C},{\mathfrak M})$ is an AECat.
\end{proposition}

\begin{proposition}\label{aecat_coroll_2}
Let $\lambda>\card{L}$, $\card{\Omega}$ be an strongly inaccessible cardinal.
  Let ${\mathfrak C}$ be any of the categories of the form
  $(\hbox{\normalfont \bfseries SubpShMod}(T),\hbox{\normalfont \bfseries ShMod}(T))_x$,
  $(\hbox{\normalfont \bfseries SubpShMod}(T),\hbox{\normalfont \bfseries pShMod}(T))_y$ or
  \mbox{$(\hbox{\normalfont \bfseries SubShMod}(T),\hbox{\normalfont \bfseries ShMod}(T))_z$}, where for $v\in \set{x,y,z}$
  either $v=L$-mono and the theory $T$ is a set of $\forall{\mathcal E}_{1}$ sentences, or $v=$elt-$L$-emb (and $T$ is arbitrary).
Let ${\mathfrak M}$ be in each case the full subcategory of ${\mathfrak C}$ whose class of objects consists of all pairs of the
   form $(M,M)$.
   Then $({\mathfrak C},{\mathfrak M})$ is an AECat.
\end{proposition}

\vskip12pt

The main results of \cite{K22} are uniqueness theorems for independence relations on Galois types for
AECats with the amalgamation property.

The theorems are of the form that an independence relation satisfying certain conditions coincides
with non-dividing or non-forking for various notions of dividing and forking for Galois types,
in particular so called \emph{isi-dividing},
\emph{isi-forking} and \emph{long Kim-dividing}, where `isi' abbreviates \emph{initial segment invariant}.

We give the definitions adapted to our categories under the hypotheses that they have the amalgamation property.

Galois types are equivalence classes of pairs $(\tupof{a_i}{i\in I};(M,M))$ where
for each $i\in I$ we have that $a_i:(A_i,M_i)\longrightarrow (M,M)$ is a morphism.

\begin{definition} Let $(\tupof{a_i}{i\in I};(M,M))$, $(\tupof{a'_i}{i\in I};(M',M'))$ 
be pairs where 
for each $i\in I$ we have that $a_i:(A_i,M_i)\longrightarrow (M,M)$ and
$a'_i:(A'_i,M'_i)\longrightarrow (M',M')$ and are morphisms.
We say $(\tupof{a_i}{i\in I};(M,M))$ and $(\tupof{a'_i}{i\in I};(M',M'))$ are of the same \emph{Galois type}
and write $\gtp(\tupof{a_i}{i\in I};(M,M))=\gtp(\tupof{a'_i}{i\in I};(M',M'))$ if
for all $i\in I$ we have that $(A_i,M_i=(A'_i,M'_i)$ and there are morphisms 
$m:(M,M)\longrightarrow (N,N)$ and $m':(M',M')\longrightarrow (N,N)$ such that
for each $i\in I$ we have that the following diagram commutes.
\begin{equation*}
\begin{tikzcd}
(A_i,M_i) \arrow{r}{a} \arrow[swap]{d}{a'} & (M,M) \arrow{d}{m} \\
(M',M') \arrow{r}{m'} & (N,N)
\end{tikzcd}
\end{equation*}
\end{definition}

\begin{definition} Let $\kappa$ be a cardinal.
  A sequence $\tupof{b_i}{i\in \kappa}$, where for all $i<\kappa$ we have
  $b_i:(B,M)\longrightarrow (N,N)$, and a continuous chain of initial segments
  $\tupof{N_i}{i\in \kappa}$ is an \emph{initial segment invariant} (\emph{isi-}) sequence if for all $i<j<\kappa$ we have
  $\gtp(b_i,n_i;N) = \gtp(b_j,n_i;N)$ and $a_i$ factors through $m_{i+1}$.
  The sequence is an isi-sequence \emph{over} $c$ for some $c:(C,N_C)\longrightarrow N$ if
  $c$ factors through $n_0:(N_0,N_0)\longrightarrow (N,N)$.
\end{definition}

\begin{definition} Let $a:(A,M_A) \longrightarrow (M,M)$, $b:(B,M_B) \longrightarrow (M,M)$,
  $c:(C,M_C) \longrightarrow (M,M)$ be morphisms. Let $m:(M,M)\longrightarrow (N,N)$ be a morphism and
  $\tupof{b_i:(B,M_B)\longrightarrow (N,N)}{i\in I}$ a sequence of morphisms  such that for all $i\in I$ we have that
  $\gtp(b_i,c;N)=\gtp(b,c;M)$. We say $\gtp|(a,b,c;M)$ is inconsistent with $\tupof{b_i}{i\in I}$ if
  there is no morphism $n:(N,N)\longrightarrow (P,P)$ such that there is a morphism $a':(A,M_A)\longrightarrow (P,P)$
  such that for all |$i\in I$ we have $\gtp(a',b_i,c;P)=\gtp(a,b,c;M)$
  
  We say $\gtp(a,b,c;M)$ \emph{long divides} over $c$ if there is some regular
  cardinal $\mu$ such that for all regular $\lambda>\mu$
there is a morphism $m:(M,M)\longrightarrow (N,N)$ and a sequence of morphisms
$\tupof{b_i:(B_i,M_i)\longrightarrow (N,N)}{i<\lambda}$ such that for all $i\in I$ we have that
$\gtp(b_i,c;N)=\gtp(b,c;M)$ and there is some cardinal $\kappa<\lambda$ such that for every $I\in [\lambda]^{\kappa}$
we have that $\gtp(a,b,c;M)$ is inconsistent with $\tupof{b_i}{i\in I}$.

We say $\gtp(a,b,c;M)$ \emph{isi-divides} over $c$ if the $\tup-of{b_i}{i<\lambda}$ are required to be isi-sequences over $c$.
We write $(A,M_A)\forkindep[(C,M_C)]{\mathrm{isi-d}} (B,M_B)$ if $\gtp(a,b,c;M)$ does not isi-divide over $c$.

Kamsma remarks that isi-forking says $\gtp(a,b,c;M)$ implies a a collection of Galois types each isi-dividing
over $c$ -- we omit the details from \cite{K22}. We similarly write
$(A,M_A)\forkindep[(C,M_C)]{\mathrm{isi-f}} (B,M_B)$ if $\gtp(a,b,c;M)$ does not isi-fork over $c$.

Long Kim-dividing is exactly as long dividing, but with the additional requirement that the $\tupof{b_i}{i<\lambda}$
are $\forkindep[c]{\mathrm{isi-f}}$-independent, that is for all $i<\lambda$ we have $b_i\forkindep[c]{N} N_i$.
This time we write
$(A,M_A)\forkindep[(C,M_C)]{\mathrm{lK}} (B,M_B)$ if $\gtp(a,b,c;M)$ does not long Kim-divide over $c$.
\end{definition}

Kamsma's principle results as applied to our categories are the following. (Again, we eschew the
definitions of simple and NSOP$_1$-like independence relations and what exactly the
$B$-existence axiom may be in the following two propositions. The reader should refer to \cite{K22} for the details.
Our point is to note that isi-dividing and forking and long Kim-dividing play a canonical r\^ole.)

\begin{proposition}\label{applying_Kamsma's_results_1} Let $\lambda>\card{L}$, $\chi(\Omega)$ be a regular cardinal.
  Let ${\mathfrak C}$ be any of the categories
  \mbox{$(\hbox{\normalfont \bfseries bdSubShMod}(T),\hbox{\normalfont \bfseries ShMod}(T))_x$}, where
  either $x=L$-mono and the theory $T$ is a set of $\forall{\mathcal E}_{1}$ sentences, or $x=$elt-$L$-emb (and $T$ is arbitrary)
  and suppose $\hbox{\normalfont \bfseries ShMod}(T))_x$ has the amalgamation property.

  If $\forkindep[]{}$ is a \emph{simple} independence relation, $\forkindep[]{} = \forkindep[]{\mathrm{isi-d}} =
  \forkindep[]{\mathrm{isi-f}}$ over $\mathrm{base}(\forkindep[]{})$.

  Suppose $B$ is some \emph{base class}
  and \mbox{$(\hbox{\normalfont \bfseries bdSubShMod}(T),\hbox{\normalfont \bfseries ShMod}(T))_x$}
  satisfies the \emph{$B$-existence axiom}. If $\forkindep[]{}$ is an \emph{$\mathrm{NSOP}_1$-like}
  independence relation over $B$, then $\forkindep[]{} = \forkindep[]{\mathrm{lK}}$.
\end{proposition}

\begin{proposition}\label{applying_Kamsma's_results_2}
  Let $\lambda>\card{L}$, $\card{\Omega}$ be an strongly inaccessible cardinal.
  Let ${\mathfrak C}$ be any of the categories of the form
  $(\hbox{\normalfont \bfseries SubpShMod}(T),\hbox{\normalfont \bfseries ShMod}(T))_x$,
  $(\hbox{\normalfont \bfseries SubpShMod}(T),\hbox{\normalfont \bfseries pShMod}(T))_y$ or
  $(\hbox{\normalfont \bfseries SubShMod}(T),\hbox{\normalfont \bfseries ShMod}(T))_z$, where for $v\in \set{x,y,z}$
  either $v=L$-mono and the theory $T$ is a set of $\forall{\mathcal E}_{1}$ sentences, or
  $v=$elt-$L$-emb (and $T$ is arbitrary),
  and suppose $\hbox{\normalfont \bfseries ShMod}(T))_x$ or $\hbox{\normalfont \bfseries pShMod}(T))_y$, respectively,
  has the amalgamation property.

  If $\forkindep[]{}$ is a \emph{simple} independence relation, $\forkindep[]{} = \forkindep[]{\mathrm{isi-d}} =
  \forkindep[]{\mathrm{isi-f}}$ over $\mathrm{base}(\forkindep[]{})$.

  Suppose $B$ is some \emph{base class}
  and suppose
    $(\hbox{\normalfont \bfseries SubpShMod}(T),\hbox{\normalfont \bfseries ShMod}(T))_x$,
  $(\hbox{\normalfont \bfseries SubpShMod}(T),\hbox{\normalfont \bfseries pShMod}(T))_y$ or
  $(\hbox{\normalfont \bfseries SubShMod}(T),\hbox{\normalfont \bfseries ShMod}(T))_z$, respectively,
  satisfies the \emph{$B$-existence axiom}. If $\forkindep[]{}$ is a \emph{$\mathrm{NSOP}_1$-like}
  independence relation over $B$, then $\forkindep[]{} = \forkindep[]{\mathrm{lK}}$.
\end{proposition}

\begin{proof} The propositions are immediate
  from Corollary (\ref{aecat_coroll_1}) and Corollary (\ref{aecat_coroll_2}), respectively, by
  \cite{K22}, Theorems (1.2) and (1.1).
  \end{proof}

\vfill\eject


\begin{thebibliography}{299}

\bibitem{Adamek-Rosicky} Ji{\v r}\'i Ad\'amek and Ji{\v r}\'i Rosick\'y,
  \emph{Locally presentable and accessible categories}
   London Mathematical Society Lecture Note Series, 189,
   Cambridge University Press,
   (1994).

 \bibitem{AAM1} José Goudet Alvim, Caio de Andrade Mendes and Hugo Luiz Mariano,
   \emph{Q-Sets and Friends: Categorical Constructions and Categorical Properties},
   arXiv.2302.03123,
   (2023).

\bibitem{AAM2} José Goudet Alvim, Caio de Andrade Mendes and Hugo Luiz Mariano,
  \emph{Q-Sets and Friends: Regarding Singleton and Gluing Completeness},
  arXiv.2302.03691,
  (2023).
   
 \bibitem{Ara} Hisashi Aratake,
   \emph{Sheaves of structures, Heyting-valued structures, and a generalization of Łoś’s theorem},
   Mathematical Logic Quarterly,
   67,
   (2021),
   445–468.
   
 \bibitem{BBHU} Ita\"i Ben Yaacov, Alexander Berenstein, C.~Ward Henson and Alexander Usvyatsov,
   \emph{Model theory for metric structures}, in, eds.~Zoé Chatzidakis, Dugald Macpherson, Anand Pillay and Alex Wilkie,
   \emph{Model theory with applications to algebra and analysis, volume 2},
   Cambridge University Press,
   (2008).

\bibitem{BYP} Ita\"i Ben Yaacov and Arthur Pedersen,
  \emph{A proof of completeness for continuous first-order logic}
  Journal of Symbolic Logic,
  75(1),
  (2010)
  168-190.

\bibitem{BYU} Ita\"i Ben Yaacov and Alexander Usvyatsov,
  \emph{Continuous first order logic and local stability},
  Transactions of the American Mathematical Society,
  362(10),
  (2010),
  5213-5259.

 \bibitem{Becker-Kechris} Howard Becker and Alexander Kechris,
   \emph{The descriptive set theory of Polish group actions},
   London Mathematical Society Lecture Note Series, 232,
   Cambridge University Press,
   (1996).
   
 \bibitem{Bena} J. Benavides
   \emph{The logic of sheaves, sheaf forcing and the independence of the Continuum Hypothesis}
   preprint,
   arxiv.1111.5854,
   (2012)
   
 \bibitem{Beerio} Diego Berr\'io,
   \emph{Logic of sheaves of structures on a locale},
   Master's thesis,
   Universidad Nacional de Colombia,
   (2015),
   70pp.

   
 \bibitem{BHV} Alexander Berenstein, Tapani Hyttinen and Andres Villaveces,
   \emph{Hilbert spaces with generic predicates}
   Revista Colombiana de Matem\'aticas,
   52(1),
   (2018),
   107-130.

 \bibitem{Borceux} Francis Borceux,
   \emph{Handbook of Categorical Algebra: Sheaf Theory. Vol. 3},
   Encyclopedia of Mathematics and its Applications,
   Cambridge University Press,
   (1994).
   
 \bibitem{BNI} Taus Brock-Nannestad and Danko Ilik,
   \emph{An intuitionistic formula hierarchy based on high-school identities},
   Mathematical Logic Quarterly,
   65(1),
   (2019),
   57-79.

 \bibitem{B00} Andreas Brunner,
   \emph{The Method of Constants in the Model Theory of Sheaves over a Heyting Algebra},
   Ph.D. thesis, University of São Paulo,
   (2000).
   
   
 \bibitem{B16} Andreas Brunner,
   \emph{Model theory in sheaves},
   South American Journal of Logic,
   2(2),
   (2016),
   379–404.

 \bibitem{BM04} Andreas Brunner and Francisco Miraglia,
   \emph{An Omitting Types Theorem for Sheaves over Topological Spaces},
   Logic Journal of the IGPL,
   12(6),
   (2004),
   525-548.

   
 \bibitem{BM14} Andreas Brunner and Francisco Miraglia,
   \emph{The method of diagrams for presheaves of $L$-structures},
   S\'eminaire de Structures, Alg\'ebriques Ordonn\'ees, Paris,
   (2014),
   1-24.

 \bibitem{BMPIII} Andreas Brunner, Charles Morgan and Darllan Pinto,
   \emph{Some contributions to presheaf model theory, II - back and forth},
   in preparation.
   
 \bibitem{Burr} Wolfgang Burr,
   \emph{The intuitionistic arithmetical hierarchy},
   in eds.~Jan van Eijck, Vincent van Oostrom and Albert Visser, \emph{Logic Colloquium '99}, 
   Lecture Notes in Logic, 17,
   CUP,
   (2004),
   51-59.
   
 \bibitem{burwer} Stanley Burris and Heinrich Werner,
   \emph{Sheaf constructions and their elementary properties},
   Transactions of the American Mathematical Society,
   248,
   (1979),
   269-309.

 \bibitem{burwer2} Stanley Burris and Heinrich Werner,
   \emph{Remarks on Boolean products}
   Algebra Universalis,
   10,
   (1980),
   333-344.

   
\bibitem{Caicedo0} Xavier Caicedo,
\emph{L\'ogica de los haces de estructuras},
 Revista de la Academia Colombiana de Ciencias Exactas, Físicas y Naturales
 19,
 (1995),
 569-586.

\bibitem{Caicedo1} Xavier Caicedo,   
\emph{Investigaciones acerca de los conectivos intuicionistas},
Revista de la Academia Colombiana de Ciencias Exactas, Físicas y Naturales,
19, (1995),
705–716.
   

\bibitem{Caicedo2}  Xavier Caicedo,
\emph{Conectivos intuicionistas sobre espacios topol\'ogicos},
Revista de la Academia Colombiana de Ciencias Exactas,
21,
 (1997),
 521-535.

  

    
\bibitem{Carson} Andrew Carson, \emph{The model completion of the theory of commutative regular rings},
  Journal of Algebra,
  27,
  (1973),
  136-146.
  
\bibitem{Zoe} Zoe Chatzdakis,
  \emph{La repr\'sentation en termes de faisceaux des mod\`eles de la th\'eorie \'el\'ementaire
       de la multiplication des entiers naturels}, in
  eds.~Chantal Berline, Kenneth McAloon, Jean-Pierre Ressayre, \emph{Model Theory and Arithmetic},
  Lecture Notes in Mathematics 890,
  Springer,
  (1981)
  90-110.
  
\bibitem{CK} Chen Chung Chang and H. Jerome Keisler,
  \emph{Model theory} (3rd edition),
  North-Holland,
  (1990).

\bibitem{Chen} Ruiyuan Chen,
  \emph{Étale structures and the Joyal-Tierney representation theorem in countable model theory},
  preprint,
  arXiv:2310.11539,
  (2023).
  
\bibitem{com} Stephen Comer,
  \emph{Complete and model complete theories of monadic algebras},
  Colloquium Mathematicum,
  34(2),
  (1976),
  183-190.
  
\bibitem{com0} Stephen Comer,
  \emph{Elementary properties of structures of sections},
  Boletin de la Sociedad Matematica Mexicana,
  19,
  (1974),
  78-85.

\bibitem{ConMir1} Marcelo Coniglio and Francisco Miraglia,
  \emph{Non-Commutative Topology and Quantales},
  Studia Logica,
  65,
  (2000)
  223–236.

\bibitem{ConMir2} Marcelo Coniglio and Francisco Miraglia,
  \emph{Modules in the Category of Sheaves over Quantales},
  Annals of Pure and Applied Logic,
  108,
  (2001)
  103–136.



 \bibitem{EFT} Heinz-Dieter Ebbinghaus, J\"org Flum and Wolfgang Thomas,
    \emph{Mathematical Logic}, (second edition),
    Undergraduate Texts in Mathematics,
    Springer,
   (1994).

  \bibitem{d'El} Christian D’Elb\'e,
    \emph{Generic expansions by a reduct},
    Journal of Mathematical Logic,
    21(3),
    (2021).
    
  \bibitem{eli} Jonas Eliasson, \emph{Ultrasheaves as sheaves on a category of ultrafilters},
    Archive for Mathematical Logic, 43(7), (2004), 825-843.

  \bibitem{ell} David Ellerman,
    \emph{Sheaves of structures and generalized ultraproducts},
    Annals of Pure and Applied Logic,
    7,
    (1974),
    163-195.

\bibitem{Fitting} Melvin Fitting,
  \emph{Intuitionistic Logic Model Theory and Forcing},
  North Holland,
  (1969).
  
\bibitem{Fle}  Jonathan Fleischmann,
  \emph{Syntactic preservation theorems for intuitionistic predicate logic},
  Notre Dame Journal of Formal Logic,
  51(2),
  (2010),
  225–245.

\bibitem{Fra50} Roland Fra\"iss\'e,
  \emph{Sur une nouvelle classification des syst`emes de relations},
  Comptes rendus hebdomadaires des s\'eances de l’Acad\'emie des Sciences,
  230,
  (1950),
  1022-1024.


\bibitem{Fra53} Roland Fra\"iss\'e,
  \emph{Sur quelques classifications des syst\`emes de relations},
  Ph.D. thesis, University of Paris, 1953,
  published in Publications Scientifiques de l’Universit\'e d’Alger,
  series A1,
  (1954),
  35–182.

    
\bibitem{Forero}    Andr\'es Forero Cuervo,
  \emph{Una demostraci\'on alternativa del teorema de ultralı\'imites},
  Revista Colombiana de Matem\'aticas,
  43,
  115-138,
  (2009).

  
\bibitem{FS} Mike Fourman and Dana Scott, \emph{Sheaves and Logic}, in
  eds.~Mike Fourman, Chris Mulvey and Dana Scott \emph{Applications of Sheaves}, Lecture Notes in Mathematics, \textbf{753},
  Springer Verlag, (1979), 302-401.

\bibitem{FuKu} Makoto Fujiwara and Taishi Kurahashi,
  \emph{Prenex normalization and the hierarchical classification of formulas},
  Archive for Mathematical Logic,
  63, (2024), 391–403.
  

\bibitem{Haj} Petr H\'ajek, \emph{Metamathematics of fuzzy logic}, Kluwer, 1998.
    
\bibitem{hoc} Melvin Hochster,
  \emph{Prime ideal structure in commutative rings},
  Transactions of the American Mathematical Society,
  142,
  (1969),
  43-60.

\bibitem{Higgs1} Denis Higgs,
  \emph{A category approach to Boolean-valued set theory},
 Preprint,
 University of Waterloo,
 (1973).


\bibitem{Higgs2} Denis Higgs,
  \emph{Injectivity in the Topos of Complete Heyting Algebra Valued Sets},
  Canadian Journal of Mathematics,
  36(3),
  (1984),
  550-568.

  
\bibitem{hod} Wilfrid Hodges, \emph{Model Theory}, Cambridge University Press, Cambridge, 1993.

\bibitem{HK11} Ulrich H\"ohle and Tomasz Kubiak,
  \emph{A non-commutative and non-idempotent theory of quantale sets}
  Fuzzy Sets and Systems,
  166(1),
  (2011),
  1-43.
  
\bibitem{Hru} Ehud Hrushovski,
  \emph{The Manin-Mumford conjecture and the model theory of difference fields},
  Annals of Pure and Applied Logic,
  112,
  (2001),
  43–115.
  
 \bibitem{Hyland} Martin Hyland,
   \emph{Categorical Model Theory},
   lecture at Models and Fields meeting, Durham,
   www.maths.dur.ac.uk/events/Meetings/LMS/2009/Movies/hyland.wmv,
   (2009).
 
  \bibitem{johnstone} Peter Johnstone, \emph{Stone Spaces}, Cambridge
    Studies in Advanced Mathematics, Cambridge University Press, 1982.

  \bibitem{Joh} Peter Johnstone,
    \emph{Sketches of an elephant: a topos theory compendium},
    Oxford Logic Guides, 43, 44,
    Oxford University Press,
    (2002).
    
   
 \bibitem{K20} Mark Kamsma,
 \emph{The Kim–Pillay theorem for abstract elementary categories}
   Journal of Symbolic Logic,
   85(4)
   (2020),
   1717-1741. 

 \bibitem{K21} Mark Kamsma,
\emph{Independence relations in abstract elementary categories}
   Ph.D. thesis, UEA,
   (2021),
   161pp. 

   
 \bibitem{K22} Mark Kamsma,
   \emph{NSOP$_1$-like independence in AECats},
   Journal of Symbolic Logic,
   89(2),
   (2024),
   724-757.

 \bibitem{KR23} Mark Kamsma and Ji{\v r}\'i Rosick\'y,
   \emph{Unstable independence from the categorical point of view},
   preprint: arXiv:2310.15804
   (2023),
   32pp.

 \bibitem{Kirby} Jonathan Kirby,
   \emph{An invitation to model theory},
   CUP,
   (2019).


 \bibitem{KrR} Alex Kruckman and Nicholas Ramsey,
  \emph{A new Kim’s lemma},
  Model Theory,
  3(3),
   (2024),
   825-860.

 \bibitem{Kha} Karim Khanaki,
   \emph{Stability, the NIP, and the NSOP: model theoretic properties of
   formulas via topological properties of function spaces},
   Mathematical Logic Quarterly,
   66(2),
   (2020),
   136-149.

 \bibitem{LL} Ren\'e Lavendhomme and T. E. Lucas,
   \emph{A non-Boolean version of Feferman–Vaught’s theorem},
   Mathematical Logic Quarterly (Zeitschrift für mathematische Logik und Grundlagen der Mathematik),
   31(19–20),
    (1985),
   299–308.
   
   \bibitem{LR18} Michael Lieberman and Ji{\v r}\'i Rosick\'y,
     \emph{Metric abstract elementary classes as accessible categories},
     Journal of Symbolic Logic,
     82(3)
     (2017),
     1022-1040.

   \bibitem{LRZ} Michael Lieberman, J. Ji{\v r}\'i Rosick\'y  and Pedro Zambrano,
     \emph{Tameness in generalized metric structures},
     Archive for Mathematical Logic,
     62,
     (2023)
     531-558.
     
\bibitem{LRV19} Michael Lieberman, Ji{\v r}\'i Rosick\'y and Sebastien Vasey,
\emph{Forking independence from the categorical point of view},
  Advances in Mathematics,
  346
  (2019),
  719-772.

\bibitem{LRV22} Michael Lieberman, Ji{\v r}\'i Rosick\'y and Sebastien Vasey,
  \emph{Induced and higher-dimensional stable independence},
  Annals of Pure and Applied Logic,
  173(7)
  (2022),
  15pp.

\bibitem{LS} Leonard Lipschitz and Dan Saracino,
  \emph{The model companion of the theory of commutative rings without nilpotent elements},
  Proceedings of the American Mathematical Society,
  38,
  (1973),
  381-387.

  
  \bibitem{Lou79} George  Loullis,
    \emph{Sheaves and Boolean valued model theory},
     Journal of Symbolic Logic,
     44(2),
     (1979).

   \bibitem{Lou78} George Loullis,
     \emph{Infinite Forcing for Boolean valued Models},
       \textgreek{Δελτίο της Ελληνικής Μαθηματικής Εταιρίας} (Bulletin of the Hellenic Mathematical Society),
     19(B),
     (1978),
     155-182.

\bibitem{Mac} Angus Macintyre, \emph{Model-completeness for Sheaves of Structures},
  Fundamenta Mathematica, 81,
  (1973),
  73-89.
  
\bibitem{Makkai-Pare} Michael Makkai and Robert Par\'e,
  \emph{Accessible categories: the foundations of categorial model theory},
  Contemporary Mathematics,  104,
  American Mathematical Society
  (1989).


\bibitem{Marker} David Marker,
  \emph{Model theory: an introduction},
  Graduate texts in mathematics, 217,
  Springer
  (2002).

\bibitem{Mil} Haynes Miller,
  \emph{Leray in Oflag XVIIA: The origins of sheaf theory, sheaf cohomology, and spectral sequences},
  in \emph{Jean Leray (1906-1998)},
  Gazette des Mathematicens,
  84 suppl,
  (2000),
  17–34.
  
 \bibitem{M88} Francisco Miraglia,
   \emph{The downward L\"owenheim-Skolem theorem for L-structures in $\Omega$-sets},
   in Walter Carnielli and Lu\'z de Alcantara (eds),
   Methods and Applications of Mathematical Logic,
   Contemporary Mathematics, 69,
   American Mathematical Society
   (1988),
   189-208.

 \bibitem{M90} Francisco Miraglia,
   \emph{Ultraprodutos generalizados e feixes de estruturas sobre uma álgebra de Heyting completa},
   Tese de Livre Docência
   Universidade de S\~ao Paulo,
   (1990).
   
  
 \bibitem{Mir06} Francisco Miraglia,
   \emph{An introduction to partially ordered  structures and sheaves},
   Contemporary Logic Series, 1,
   Polimetrica, International Scientific Publisher,
   Milan,
   (2006).

   \bibitem{MirSol} Francisco Miraglia and Ugo Solitro,
     \emph{Sheaves and Presheaves over Right-sided Idempotent Quantales},
     Logic Journal of the IGPL,
     6(4),
     (1998),
     545–600.
   
 \bibitem{Montoya}    Andr\'es Montoya,
   \emph{Contribuciones a la teorı\'ia de modelos de haces},
   Lecturas Matem\'aticas,
   28,
   (2007),
   5–37.

   
 \bibitem{mul} Chris Mulvey, \emph{Intuitionistic algebra and the representation of rings},
   Memoirs of the American Mathematical Society, 69,
   (1988),
   189-209.

\bibitem{Mut1} Scott Mutchnik,
  \emph{On NSOP$_2$ theories},
  \emph{to appear in} Journal of the European Mathematical Society,
  preprint: arXiv:2206.08512,
   (2023),
   21pp.
  
\bibitem{Mut2} Scott Mutchnik,
  \emph{Properties of independence in NSOP$_3$ theories},
  preprint: arXiv:2305.09908,
  (2023),
  29pp.

    
\bibitem{Mut3} Scott Mutchnik,
  \emph{Conant-independence and generalized free amalgamation},
  \emph{to appear in} Journal of Mathematical Logic,
  preprint: arxiv:2210.07527,
  (2024),
  29pp.
  
\bibitem{ochvil} Maciol Ochoa and Andres Vilaveces, \emph{Sheaves of metric structures},
  Logic, Language, Information and Computation, 23rd
  Workshop WoLLIC 2016, (2016), 297-315.

\bibitem{deP} Valeria de Paiva,
  \emph{Lineales: Algebras and Categories in the Semantics of Linear Logic},
  in eds.~Dave Barker-Plummer, David Beaver, Johan van Benthem and Patrick Scotto di Luzio
  \emph{Words, Proofs,and Diagrams},
  CSLI Publications,
  (2002),
  123-142.

  \bibitem{Pappas} Peter Pappas 
    \emph{The model-theoretic structure of abelian group rings},
    Annals of Pure and Applied Logic,
    28,
    (1985),
    163-201.

    \bibitem{PV}  
      Moreno Pierobon and Mateo Viale,
      \emph{Boolean valued models, presheaves, and étalé spaces},
      preprint
      arxiv.2006.14852,
      (2023),
      
    
    \bibitem{Prest} Mike Prest,
  \emph{Model theory for sheaves of modules},
  in Logic and its Applications, Proc. 8th Indian Conf. ICLA, Delhi,
  Lecture Notes in Comp. Sci,
  11600,
  Springer,
  (2019),
  89-102.

    
\bibitem{PPR} Mike Prest, Vera Puninskaya and Alexandra Ralph
  \emph{Model theory of sheaves of modules},
  Journal of Symbolic Logic, 69(4),
  (2004),
  1187-1200.
  
\bibitem{PreRaj} Mike Prest and Ravi Rajani,
  \emph{Model-theoretic imaginaries and coherent sheaves},
  Applied Categorical Structures,
  17(6),
  (2009),
  517-559.


\bibitem{Prest-Ralph1}
  Mike Prest and Alexandra Ralph,
  \emph{On sheafification of modules},
  preprint,
  University of Manchester,
  (2001, revised 2004),
  available at personalpages.manchester.ac.uk/staff/mike.prest/ShvFP09Rev.pdf,
  11pp.
  
\bibitem{Prest-Ralph2} Mike Prest and Alexandra Ralph,
  \emph{Locally finitely presented categories of sheaves of modules},
  preprint,
  University of Manchester,
  (2001, revised 2004 and 2018),
  available at eprints.maths.manchester.ac.uk/1411/1/ShfifTorsTh.pdf,
  15pp.

\bibitem{Prest-Slavik}  Mike Prest and Alexander Slávik,
  \emph{Purity in categories of sheaves},
  Mathematische Zeitschrift,
  297(1-2),
  (2021),
  429-451.
    
\bibitem{Resende} Pedro Resende,
  \emph{Groupoid sheaves as quantale sheaves},
    Journal of Pure and Applied Algebra,
    216(1),
    (2012),
    41-70.

\bibitem{Resende-lectures} Pedro Resende,
  \emph{Lectures on étale groupoids, inverse semigroups and quantales},
  Lecture Notes for the GAMAP IP Meeting, Antwerp,
  (2006),
  115pp.

  \bibitem{Resende-talk} Pedro Resende,
    \emph{Quantale-valued sets, quantale modules, and groupoid actions},
    Lecture, International Category Thery Conference, Carvoeiro, Portugal, 17-23 June 2007.
    (2007).
  
\bibitem{RZ} David Reyes and Pedro H. Zambrano,
  \emph{Co-quantale valued logics},
  preprint,
  arXiv:2102.06067.
  (2021).
  
\bibitem{Pav} Ji{\v r}\'i Pavelka, \emph{On fuzzy logic I, II, III},
  Z. Math. Logik Grundlagen Math., 25, (1979), 45-52, 119-134, 447-464.

\bibitem{SUZ} Aleksy Schubert, Pawel Urzyczyn, and Konrad Zdanowski
  \emph{On the Mints Hierarchy in First-Order Intuitionistic Logic},
  Logical Methods in Computer Science, 12(4), (2017), 1-25.

  
\bibitem{sco} Dana Scott,
  \emph{Identity and Existence in Intuitionistic Logic},
  in eds.~Mike Fourman, Chris Mulvey and Dana Scott \emph{Applications of Sheaves},
  Lecture Notes in Mathematics, 753,
  Springer,
  (1979),
  660-696.


  \bibitem{SetCai} Antonio Sette and Xavier Caicedo,
  \emph{Equival\^encia elementar entre feixes},
  Notas de l\'ogica mathem\'atica
  38,
  129-141,
  (1993).

\bibitem{Shaf} Igor Shafarevich,
  \emph{Basic Algebraic Geometry, Volume 2}, (3rd edition),
  Springer,
  (2013).
  
\bibitem{Sureson1} Claude Sureson,
  \emph{A valuation ring analogue of von Neumann regularity},
  Annals of Pure and Applied Logic,
  145, (2007),
  204–222.
 
\bibitem{Sureson2} Claude Sureson,
  \emph{Model companion and model completion of theories of rings}
  Archive for Mathematical Logic,
  48,
  (2009),
  403–420.
  

\bibitem{Tenn} Barry Tennison,
  \emph{Sheaf Theory},
  London Mathematical Society Lecture Note Series, 20,
  CUP,
  (1975).
  
\bibitem{TAM} Ana Luiza Tenório, Caio de Andrade Mendes, and Hugo Mariano,
  \emph{On sheaves on semicartesian quantales and their truth values},
  arXiv.2204.08351,
  (2023).

\bibitem{TvD} Anne Troelstra and Dirk van Dalen.
  \emph{Constructivism in Mathematics, Volumes 1 and 2},
  North-Holland, 
  (1988).
  

\bibitem{Vaa} Jouko V\"a\"an\"anen,
  \emph{Model theory of second order logic}
  in ed.~Jose Iovino, \emph{Beyond First Order Model Theory, Volume II},
  Chapman and Hall,
  (2023),
  290-316.

\bibitem{vol1} Hugo Volger, \emph{The Feferman-Vaught Theorem, revisited},
  Colloquium Mathematicum,
  36,
  (1976),
  1-11.

\bibitem{vol2} Hugo Volger,
  \emph{Preservation Theorems for Limits of Structures and Global Sections of Sheaves of Structures},
  Mathematische Zeitschrift,
  166,
  (1979),
  27-54.

\bibitem{vol3} Hugo Volger,
  \emph{Filtered powers and stable boolean powers are relativized full boolean powers},
  Algebra Universalis,
  19
  (1984),
  399-402.

\bibitem{Weis} Volker Weispfenning,
  \emph{Model-completeness and elimination of quantifiers for subdirect products of structures},
  Journal of Algebra,
  36, (1975), 252-277.
  
\end{thebibliography}
\end{document}